\documentclass[11pt]{article}

\usepackage[margin=1in]{geometry}
\usepackage{amsmath,amssymb,amsthm}
\usepackage{mathtools}
\usepackage{enumitem}
\usepackage[numbers]{natbib}
\usepackage[colorlinks=true,linkcolor=blue,citecolor=blue]{hyperref}

\theoremstyle{plain}
\newtheorem{theorem}{Theorem}
\newtheorem{proposition}{Proposition}
\newtheorem{lemma}{Lemma}
\newtheorem{corollary}{Corollary}
\newtheorem{example}{Example}

\theoremstyle{definition}
\newtheorem{assumption}{Assumption}[]

\newtheorem{remark}{Remark}
\newtheorem{definition}{Definition}

\newcommand{\R}{\mathbb{R}}
\newcommand{\E}{\mathbb{E}}

\newcommand{\Y}{\mathbb{Y}}
\newcommand{\N}{\mathcal{N}}
\newcommand{\C}{\mathcal{C}}
\newcommand{\T}{\mathcal{T}}
\newcommand{\Ra}{\mathcal{R}}
\newcommand{\D}{\mathcal{D}}

\newcommand{\Hkx}{\mathcal{H}_{k_X}}
\newcommand{\Hkz}{\mathcal{H}_{k_Z}}
\newcommand{\HKx}{\mathcal{H}_{K_X}}
\newcommand{\HKz}{\mathcal{H}_{K_Z}}

\newcommand{\Nt}{\N(T)}
\newcommand{\Na}{\N(A)}
\newcommand{\Pna}{P_{\Na}}
\newcommand{\oX}{\overline{X}}
\newcommand{\wX}{\widetilde{X}}

\newcommand{\wY}{\widetilde{Y}}
\newcommand{\oZ}{\overline{Z}}
\newcommand{\wZ}{\widetilde{Z}}
\newcommand{\norm}[1]{\left\lVert #1 \right\rVert}

\newcommand{\Id}{\mathrm{Id}}

\title{\bf Nonparametric Instrumental Kernel Regression:\\
The Well- and Misspecified Setting\\ Under General Link Conditions}
\author{
  Guillaume Baquet\\
  \small European Center for Advanced Research in Economics and Statistics (ECARES)\\
  Solvay Brussels School of Economics and Management, Universit\'e libre de Bruxelles\\
  \texttt{guillaume.baquet@ulb.be}
  \and
  Germain Van Bever\\
  \small European Center for Advanced Research in Economics and Statistics (ECARES)\\
  Department of Mathematics, Universit\'e libre de Bruxelles and Universit\'e de Namur\\
  \texttt{germain.van.bever@ulb.be}}
\date{\today}

\begin{document}
\maketitle

\begin{abstract}
We study nonparametric instrumental variable regression when the outcome takes values in a
separable Hilbert space $\mathbb Y$, for instance a curve, a density or the kernel mean embedding
of a distribution. The structural function solves an ill-posed operator equation whose
conditional expectation operator need not be injective, and whose solutions need not belong to
the hypothesis space. We analyse a two-stage kernel estimator that extends kernel instrumental
variable regression to $\mathbb Y$-valued outcomes and allows spectral regularisation in both
stages, including Landweber iteration and Showalter's method in the second stage. We show that the
second stage is a regularisation method in the Hilbert scale generated by the inverse of the
kernel integral operator. This viewpoint identifies the solution towards which the estimator
converges: the minimum-norm solution in the reproducing kernel Hilbert space in the
well-specified case and, in the misspecified case, the unique solution in the range of the
continuous extension of the projection operator associated to the well-specified setting. 
Under source conditions and link conditions expressed
through general index functions, we derive non-asymptotic bounds in a scale of interpolation
norms that hold with high probability in terms of the first-stage error, with or without identification and
in the well- and misspecified settings. For polynomial link conditions, we recover the rate of
scalar kernel instrumental regression and extend it to targets outside the hypothesis space; for
exponential link conditions, we obtain logarithmic rates in the severely ill-posed regime.
\end{abstract}

\section{Introduction}
In many applications the outcome of interest is a whole function rather than a number:
electricity-load curves, wage or birth-weight distributions, spectrograms or growth
trajectories. Such outcomes are naturally modelled as random elements of a separable Hilbert
space~$\Y$, for instance $L^2$, a Sobolev space, or a reproducing kernel Hilbert space (RKHS)
in which distributions are embedded. When the regressor $X$ is endogenous, possibly due to omitted
variables, simultaneity or measurement error, regressing $Y$ on $X$ does not recover the
causal effect. An instrument $Z$ that shifts $X$ but affects $Y$ only through $X$ restores
identification, leading to the model
\begin{equation}
  \label{eq:model}
  Y = F_0(X) + U, \qquad \E[\,U \mid Z\,] = 0, \qquad Y,\, U \in \Y,
\end{equation}
where $F_0 \colon E_X \to \Y$ is the structural function and the conditional expectation is
taken in the Bochner sense. Consequently, $F_0$ solves the ill-posed equation
$\E[Y \mid Z] = \E[F_0(X) \mid Z]$, a $\Y$-valued version of nonparametric instrumental
variable (NPIV) regression \citep{newey2003instrumental,darolles2011nonparametric}. This
problem sits between several literatures that have developed largely apart: regression with
Hilbert-valued responses, studied mostly under exogeneity, and IV regression, studied mostly
with scalar responses.
 
\paragraph{Hilbert-valued responses under exogeneity.}
In functional data analysis, linear models with a functional response have been estimated by
functional principal components and by penalisation
\citep{ramsay2005functional,yao2005functional,sun2018optimal}. Nonparametric alternatives
include kernel smoothing with functional regressors and responses \citep{ferraty2011kernel},
additive models with general Hilbertian responses \citep{jeon2020additive}, and Fr\'echet
regression for metric-space-valued responses (either a global or local, nonparametric, version)
\citep{petersen2019frechet}. In machine learning, operator-valued kernels and their
vector-valued RKHSs \citep{micchelli2005learning,carmeli2006vector} yield nonparametric
kernel regression with functional outputs \citep{kadri2016operator,meunier2024optimal,li2024towards,mollenhauer2022learning}. 
Some of these works \citep{mollenhauer2022learning,li2024towards,meunier2024optimal} show that ridge and general spectral-regularisation
estimators attain the same optimal rates as in the scalar case, independently of the output dimension.  
None of these methods accounts for unmeasured confounding.
 
\paragraph{Classical IV regression.}
With a scalar response, NPIV is a well-understood ill-posed inverse problem. Estimators based
on sieves, on kernel density estimation with Tikhonov regularisation, or on penalisation are
available
\citep{newey2003instrumental,hall2005nonparametric,darolles2011nonparametric,florens2011identification},
together with minimax $L^2$ rates \citep{chen2011rate} and sup-norm rates with uniform
inference \citep{chen2018optimal}. 
\citet{florens2015instrumental} allow functional regressors
and instruments, still with a scalar response, but only in a linear setting.
Hilbert-valued responses under endogeneity have, to the best of our knowledge,
only been treated in linear models. \citet{benatia2017functional} estimate a linear operator
between function spaces by Tikhonov regularisation with functional instruments,
\citet{seong2024functional} develop FPCA-based functional IV estimators with asymptotic
inference, and \citet{vandijcke2026iv} extend global Fr\'echet regression to endogenous
covariates for distribution-valued outcomes in 2-Wasserstein space. In all three, the
structural map is linear in the regressor.
 
\paragraph{Machine-learning IV regression.}
Machine-learning approaches make the structural function nonparametric and flexible.
Two-stage methods generalize two-stage least squares with kernel features
\citep{singh2019kernel,meunier2024nonparametric} or neural features \citep{hartford2017deep,xu2021learning},
and both variants are now known to be minimax-optimal
\citep{meunier2024nonparametric,kim2025optimality}. Conditional-moment methods instead solve a
min--max problem over test functions of the instrument
\citep{lewis2018adversarial,bennett2019deep,muandet2020dual,dikkala2020minimax,zhang2023instrumental}.
Distributional methods target the full interventional distribution of the outcome
\citep{kook2024instrumental,holovchak2025distributional}. All of these are formulated for
scalar, or at most finite-dimensional, outcomes. Hilbert-valued regression does appear inside
kernel IV, but only on the treatment side: its first stage estimates the conditional mean
embedding of $X$ given $Z$, while the outcome enters through a scalar second-stage regression.
However, all of these works assume their problem is well-specified, i.e., they assume a solution exists in their hypothesis space.
To the best of our knowledge, no consistency results exist in the misspecified setting.
 
\paragraph{Contributions.}
We build on the works of \cite{singh2019kernel,meunier2024nonparametric}.
We extend their two-stage kernel model to Hilbert-space-valued responses and obtain convergence rates in the well- and misspecified setting.
We also assume a more general link condition than that of \cite{meunier2024nonparametric}.
Our results therefore also cover the so-called severely/exponentially ill-posed setting, besides the classical mildly ill-posed setting.
We also allow for a larger class of regularisation filters in the second stage, whereas only Tikhonov was studied in these previous works.
As a byproduct, we also discuss and characterise the solution towards which a regularisation method converges in the context of Hilbert scales when the operator in non-injective and the problem is misspecified.
To our knowledge, no result in this non-injective misspecified setting can be found in the literature of inverse problems as either well-specifiedness or injectivity of the operator is assumed.
For exemple, \cite{natterer1984error,nair2005regularization} (in a deterministic context) and \cite{rastogi2023inverse} (in a statistical context) all assume injectivity.

\paragraph{Organisation of the paper.}
Section~\ref{sec:notations} sets the notation. Section~\ref{sec:model} introduces the model, the
spectral filters and the two-stage estimator. Section~\ref{sec: problem is Hilbert scaled} shows that the
second stage is a regularisation method in the Hilbert scale generated by
$A_{\mathcal L_X^{-1}}$ (Proposition~\ref{prop:problem is hilbert scaled}). Section~\ref{sec: solution selection} identifies
the solution towards which the estimator converges. 
It corresponds to the solution of minimal $S_2([\mathcal H_{k_X}]^1,\mathbb Y)$-norm in the well-specified case. 
In the misspecified case, the limit is the unique element of $\overline\Pi(\mathcal S)$, the extension of the projection operator on $L^2(X,\mathbb{Y})$, provided it exists. We also give sufficient conditions for the existence of $\overline\Pi$.
Section~\ref{sec: convergence rates} states the assumptions, the non-asymptotic
bound (Theorem~\ref{th:main rate}) and the choice of the regularisation parameters
(Proposition~\ref{prop:rates}). It then works out polynomial and exponential link conditions
(Examples 1 and 2) and the case of a single split sample
(Corollary~\ref{cor:theta1}). Section~\ref{sec: conclusions} concludes and lists open problems.

The appendices are organised as follows. Appendix~\ref{app:hilbertscalesreg} develops
regularisation in Hilbert scales for non-injective operators and misspecified targets; these
results may be of independent interest. Appendix~\ref{app:RKHS} recalls vector-valued
reproducing kernel Hilbert spaces and their interpolation spaces. Appendix~\ref{app:AC}
studies the operators $A_C$, including their representation as $C^*\otimes\mathrm{Id}$ on
tensor products. Appendix~\ref{app: lipschitz} proves the operator-Lipschitz property of the
Landweber and Showalter filters. Appendix \ref{app: practical} gives the closed form of the
estimator in terms of Gram matrices. Appendix~\ref{app: proofs} contains the proofs of the
main results and the auxiliary lemmas.

\section{Notations}\label{sec:notations}
Let $\mathbb{Y} $ be a separable Hilbert space and $(\Omega, \mathcal{F},\mathbb{P} )$ be the underlying probability space with expectation operator $\mathbb{E} $.
Define the random variables $Y,X$ and $Z$ taking values in their respective measurable spaces $(\mathbb{Y} ,\mathcal{B}_{\mathbb{Y}}),(E_X ,\mathcal{B}_{E_X})$ and $(E_Z ,\mathcal{B}_{E_Z})$
where $E_X$ and $E_Z$ are second countable locally compact Hausdorff spaces and $\mathcal{B}_S$ denotes the Borel $\sigma$-field of $S$.
Since $E_X$ is Polish, it has a regular conditional distribution $\mathbb{P}_{X|Z}$ given $Z$.
Let $\mathbb{P} _X$, $\mathbb{P}_Z$, $\mathbb{P}_{XZ}$ and $\mathbb{P}_{ZY}$  denote the pushforward probability measures of $X$, $Z$, $(X,Z)$ and $(Z,Y)$.

In the sequel, $L_2(S,\mathcal{A} ,\mu,\mathbb{B})$ will denote the space of strongly $\mathcal{A} -\mathcal{B}_\mathbb{B} $ measurable $\mathbb{B} $-valued functions which are square Bochner-integrable\footnote{See, for example, \citep[Chapter 2]{diestel1977vector} or \cite{hytonen2016analysis} for definition and general treatment of Bochner integration.} with respect to $\mu$, where $\mathcal{B}_\mathbb{B}$ is the Borel $\sigma$-field of $\mathbb{B}$, $(S,\mathcal{A} ,\mu)$ is a measure space and $\mathbb{B} $ is a separable Banach space.
If $\mathbb{B} =\mathbb{R} $, one recovers the usual space of Lebesgue square-integrable real-valued functions, which is abbreviated by $L_2(S,\mathcal{A} ,\mu)$. In the sequel, we will also make use of the abbreviated notations 
$$L_2(X,\mathbb{Y}):=L_2(E_X,\mathcal{B}_{E_X} ,\mathbb{P} _X,\mathbb{Y})\qquad\textrm{ and }\qquad L_2(X):= L_2(E_X,\mathcal{B}_{E_X} ,\mathbb{P} _X).$$
A similar notation will be used for $Z$.

Let $\mathbb{H} $ and $\mathbb{H} '$ be two separable real Hilbert spaces with scalar products $\left\langle \cdot,\cdot\right\rangle _\mathbb{H} $ and $\left\langle \cdot,\cdot\right\rangle _{\mathbb{H}'}$.
Let $\mathcal{L} (\mathbb{H} ,\mathbb{H} ')$ and $S_2(\mathbb{H} ,\mathbb{H} ')$ denote, respectively, the space of bounded operators and Hilbert-Schmidt operators from $\mathbb{H} $ to $\mathbb{H} '$. 
In the case $\mathbb{H} '=\mathbb{H} $, the shorthand notations $\mathcal{L} (\mathbb{H})$ and $S_2(\mathbb{H})$ will be used, respectively.
For a bounded operator $A$, let $A^*,\mathcal{R} (A)$, $\mathcal{N} (A)$ and $\sigma(A)$ denote, respectively, its adjoint, range, kernel and spectrum. For some subspace $M\subseteq \mathbb{H} $, we denote by $M^\perp$ its orthogonal complement. 

We will denote by $\mathbb{H} \otimes \mathbb{H}'$ the tensor product between $\mathbb{H} $ and $\mathbb{H}'$. This is also a Hilbert space, which contains the elementary tensors $u\otimes v$ for all $(u,v)\in \mathbb{H} \times \mathbb{H}'$. The tensor $u\otimes v$ is a map from $\mathbb{H} \times \mathbb{H}'$ to $\mathbb{R} $ which satisfies the property $\left\langle u_1\otimes v_1,u_2\otimes v_2\right\rangle_{\mathbb{H} \otimes \mathbb{H}'} =\left\langle u_1,u_2\right\rangle_\mathbb{H} \left\langle v_1,v_2\right\rangle_{\mathbb{H}'}$, known as the universal property.
The tensor product Hilbert space is the completion of these elementary elements.
See \citep[Chapter 2]{murray1936rings} for more details on the construction of this Hilbert space.
Although this  construction is, in some sense, abstract, one can identify this space as the space $S_2(\mathbb{H} ,\mathbb{H}')$ since there exists an isometric isomorphism between these two spaces, i.e., $\mathbb{H} \otimes \mathbb{H}'\cong S_2(\mathbb{H} ,\mathbb{H}')$.
This isometric isomorphism, denoted $\Psi$, is identified through the relation 
$$\Psi(u\otimes v)=\left\langle \cdot,u\right\rangle_\mathbb{H} v\ :\ \mathbb{H}\rightarrow\mathbb{H}':x\mapsto \langle x,u\rangle_\mathbb{H} v.$$
The elementary tensors are therefore identified with rank-one operators. 
Consequently, the right-hand side of the latter equation is often denoted, with some abuse of notation, as $u\otimes v$ itself.
In the sequel, $u\otimes v$ will denote, interchangeably, the elementary tensor or the rank-one operators $\left\langle \cdot,u\right\rangle_\mathbb{H}  v$. The nature of $u\otimes v$ will always be clear from context.

For operators $A:\mathbb{H} \to \mathbb{H}'$ and $B:\mathbb{K} \to \mathbb{K}'$ between Hilbert spaces, we say that $A$ and $B$ are unitary equivalent, denoted $A\cong B$, if there exist isometric isomorphisms
$\Phi_1:\mathbb{H} \to \mathbb{K} $ and $\Phi_2:\mathbb{H}' \to \mathbb{K}'$ such that $A=\Phi_2^{-1} B \Phi_1 $.

Finally, let $\mathbb{H}''$ be another separable real Hilbert space. 
For a densely defined operator $C:\mathcal{D}(C)\subset \mathbb{H} \to \mathbb{H}'$, let $A_C:\mathcal{D}(A_C)\subseteq S_2(\mathbb{H}',\mathbb{H}'')\to S_2(\mathbb{H},\mathbb{H}'')$ denote the operator that acts as $K\mapsto KC$,
where $$\mathcal{D}(A_C):=\{K\in S_2(\mathbb{H}',\mathbb{H}''): KC \text{ extends uniquely to an operator in }S_2(\mathbb{H},\mathbb{H}'')\}.$$
See Appendix \ref{app:AC} for more details on this operator.

\section{Model Description}\label{sec:model}
Let $Y\in\mathbb{Y}$, $X\in E_X$ and $Z\in E_Z$ be the random variables introduced in Section~\ref{sec:notations}. Throughout, we assume the regression model
\begin{equation}\label{model 1}
   Y=F_0(X)+U,
\end{equation}
such that $\mathbb{E} [U|Z]=0$ almost surely, $F_0\in L_2(X,\mathbb{Y} )$ and $\mathbb{E}\left\lVert Y\right\rVert^2_\mathbb{Y} < +\infty $. The variable $Z$ is called the instrumental variable. Note that the conditional expectation is well-defined as $\mathbb{Y}$ is a separable Hilbert space. 

Taking conditional expectation of (\ref{model 1}) with respect to $Z$, one obtains that $F_0$ is a solution of the  inverse problem
\begin{equation}\label{initial inverse problem}
   \widetilde{r}=\mathcal{T} F,
\end{equation}
where $\widetilde{r}:=\mathbb{E} [Y|Z]$ and $\mathcal{T} :L_2(X,\mathbb{Y} )\rightarrow L_2(Z,\mathbb{Y} )\ :\ F \mapsto \mathbb{E} [F(X)|Z]$.
An application of Jensen's inequality shows that $\mathcal{T} $ is bounded with $\lVert \T\rVert=1 $.
It should be noted that if $\mathbb{Y} $ is infinite-dimensional, then $\mathcal{T} $ fails to be compact.
More precisely, letting $T :L_2(X)\to L_2(Z):f\mapsto \mathbb{E} [f(X)|Z]$, one can show that $\mathcal{T} \cong T\otimes \text{Id}_{\mathbb{Y}}$.
Therefore $\sigma(\mathcal{T^*}\mathcal{T})=\sigma(T^*T)$ and $\mathcal{T} $ is compact if and only if $T$ is compact and $\mathbb{Y} $ is finite-dimensional.
Importantly, $\mathcal{T} $ might fail to be injective and $F_0$ might therefore not be the only solution to (\ref{initial inverse problem}).
This happens, for example, when $X$ and $Z$ have finite support and the cardinality of $X$ exceeds that of $Z$ (see \cite{andrews2017examples}).
One can therefore only hope to retrieve a solution of (\ref{initial inverse problem}) which might not be $F_0$.

Before describing the estimation procedure, we recall the definition of a spectral filter function for our context, which is heavily based on the definition given by \cite{lin2020optimal}.
\begin{definition}\label{def: filter}
  A family of functions $g_\lambda : \mathbb{R}^+ \to \mathbb{R}$ indexed by $\lambda \in \Lambda\subseteq \mathbb{R}^+_0$ is called a filter with qualification $\rho \geq 1$ if there exists positive constants $E\leq 1$ and $\omega_\rho <+\infty$ such that
  $$\sup_{\alpha \in [0,1]} \sup_{t\in [0,\kappa^2]} \lambda^{1-\alpha}t^\alpha g_\lambda(t)\leq E \qquad\text{and}\qquad
  \sup_{\alpha \in [0,\rho]} \sup_{\lambda \in \Lambda} \sup_{t\in [0,\kappa^2]} \lambda^{-\alpha}t^\alpha |r_\lambda(t)|\leq \omega_\rho,$$
  with $r_\lambda(t):= 1-g_\lambda(t)t$. The constant $\kappa$ will be defined in Theorem \ref{th:main rate}.  
\end{definition}
We impose $E\leq 1$. This ensures $r_\lambda(t)\geq 0$ in the well-specified setting. We also impose the assumption below on the filter function.
\begin{assumption}\label{ass: Lipschitz}
  The filter $g_\lambda$ satisfies
  $$\left\lVert (A+\lambda Id)^{1/2}[g_\lambda(A)-g_\lambda(B)](B+\lambda Id)\right\rVert \leq c_g \lambda^{-1/2}\left\lVert A-B\right\rVert  $$
  for any bounded self-adjoint, nonnegative operators $A$ and $B$ with $\lVert A\rVert,\lVert B\rVert\leq \kappa^2  $ and some constant $c_g>0$ independent of $\lambda$.
\end{assumption}

We show in Appendix \ref{app: lipschitz} that the well-known Landweber and Showalter filters satisfy Definition \ref{def: filter} as well as Assumption~\ref{ass: Lipschitz}. This is also the case for the Tikhonov filter with $c_g=1$.

The adopted estimation procedure is that of \cite{singh2019kernel} and \cite{meunier2024nonparametric}. We generalize it here  to Hilbert-spaced valued $Y$, as well as allow for the use of a general spectral filter function for the second stage estimation.  

Let $k_X:E_X \times E_X \to \R$ and $k_Z:E_Z \times E_Z \to \R$ be positive semi-definite kernels associated with their respective RKHSs $\Hkx$ and $\Hkz$ and canonical mappings $\Phi_X:E_X\to \Hkx$ and $\Phi_Z:E_Z\to \Hkz$.
Define the embedding operator $\mathcal{I}_X:\mathcal{H}_{k_X}\to L_2(X)$ and set $\mathcal{L}_X:=\mathcal{I}_X\mathcal{I}_X^*$ and $\mathcal{C}_X:=\mathcal{I}_X^*\mathcal{I}_X$.
Furthermore, let $[\mathcal{H}_{k_X}]^\gamma$ be the interpolation spaces and let $\mathbf{I_X}:S_2(\mathcal{H}_{k_X},\mathbb{Y} )\to S_2(L_2(X),\mathbb{Y} )$ be the embedding.
Finally, let $S_2([\mathcal{H}_{k_X}]^\gamma,\mathbb{Y} )$ be its associated vector-valued interpolation spaces with norm $\lVert \cdot\rVert_{\gamma} $.
And let $\mathcal{I}_Z$, $\mathcal{L}_Z$, $\mathcal{C}_Z$, $[\mathcal{H}_{k_Z}]^\gamma$, $\mathbf{I_Z}$, $S_2([\mathcal{H}_{k_Z}]^\gamma,\mathcal{H}_{k_X})$ and $\lVert \cdot\rVert_{\gamma} $ be the similar counterparts for $k_Z$. 
It should be clear from the context which norm we are referring to.
We refer to Appendix \ref{app:RKHS} for more details on these objects.
We assume for simplicity that $\mathcal{L}_X$ is injective.
If this is not the case, then all results hold with $L^2(X)$ replaced by $[\mathcal{H}_{k_X}]^0$ (see Appendix \ref{app:RKHS}).
We assume the following concerning $k_X$ and $k_Z$:
\begin{assumption}\label{ass: RKHS}
The RKHSs $\mathcal{H} _{k_X}$ and $\mathcal{H} _{k_Z}$and their kernels satisfy
\begin{enumerate}
    \item $\mathcal{H} _{k_X}$ and $\mathcal{H} _{k_Z}$ are separable Hilbert spaces;
    \item $\Phi _X(X):=k_X(\cdot,X)$ and $\Phi _Z(Z):=k_Z(\cdot,Z)$ are measurable w.r.t. $\mathcal{B}_{E_X}$ and $\mathcal{B}_{E_Z}$, respectively;
    \item $\sup_{x\in E_X}k_X(x,x)\leq \kappa_X^2$ and $\sup_{z\in E_Z}k_Z(z,z)\leq \kappa_Z^2$.
\end{enumerate}
\end{assumption}
The assumptions above are standard and relatively weak. The separability assumption is satisfied if $k_X$ and $k_Z$ are continuous kernels and $E_X$ and $E_Z$ are separable. The measurability assumption is satisfied under the same hypothesis of continuity of the kernels. Note that Gaussian and Laplace kernels are examples of kernels satisfying these assumptions.
See Chapter 4 in \cite{steinwart2008support} for more details.
Define the vector-valued kernels $K_X:=k_X \text{Id}_{\Y}$ and $K_Z:=k_Z \text{Id}_{\Hkx}$ with respective vector RKHSs $\HKx$ and $\HKz$.
These kernels are very natural as they model the spaces $L_2(X,\mathbb{Y})\cong L_2(X)\otimes \mathbb{Y}  $ (resp. $L_2(Z,\mathcal{H}_{k_X})\cong L_2(Z)\otimes \mathcal{H}_{k_X}$) by spaces isometrically isomorphic to $\mathcal{H}_{k_X} \otimes \mathbb{Y} $ (resp. $\mathcal{H}_{k_Z} \otimes \mathcal{H}_{k_X}$) (see Theorem 4.4 in \cite{mollenhauer2020nonparametric}).

The idea is to perform the estimation within $\HKx$.
Assume we observe two samples $\{(\oX_i,\oZ_i)\}_{i=1}^m$ and $\{(\wZ_i,\wY_i)\}_{i=1}^{n}$ which are all i.i.d. copies of respectively $(X,Z)$ and $(Z,Y)$.
These two samples are also assumed to be mutually independent.
We then perform the estimation procedure in two stages for which we provide a computable closed-form expression solution in Appendix \ref{app: practical}.

\paragraph{Stage 1}
We estimate $M(\cdot):=\E[\Phi_X(X)|Z=\cdot]$ in $\HKz$ by
$$\widehat{M}(z):=\widehat{\C}_{ZX} g_\varepsilon(\widehat{\C}_{Z})\Phi_Z(z),$$
where $\widehat{\C}_{ZX}:=m^{-1}\sum_{i=1}^{m}\Phi_Z(\oZ_i)\otimes \Phi_X(\oX_i)$ and $\widehat{\C}_{Z}:=m^{-1}\sum_{i=1}^{m}\Phi_Z(\oZ_i)\otimes \Phi_Z(\oZ_i)$.

\paragraph{Stage 2}
Define our final estimator by
$$\widehat{F}_\lambda(x):=\widehat{\C}_{\hat{M}Y} g_\lambda(\widehat{\C}_{\hat{M}})\Phi_X(x)$$
where $\widehat{\C}_{\hat{M}Y}:=n^{-1}\sum_{i=1}^{n}\widehat{M}(\wZ_i)\otimes \wY_i$ and $\widehat{\C}_{\hat{M}}:=n^{-1}\sum_{i=1}^{n}\widehat{M}(\wZ_i)\otimes \widehat{M}(\wZ_i)$.

This is the the kernel instrumental variable (KIV) estimator of \cite{singh2019kernel}, here with $\mathbb{Y} $-valued outcomes and a general second-stage filter.
The idea of this estimator can be linked to the theory of regularisation in Hilbert scales as described in \cite{natterer1984error} or \cite{nair2005regularization} among many others.
We detail this in the following section.

\section{KIV as a Hilbert scale regularisation method}\label{sec: problem is Hilbert scaled}
The method estimates the regularized population candidate
\begin{equation}F_\lambda(x):=K_\lambda \Phi_X(x):=\C_{MY}g_\lambda(\C_M)\Phi_X(x),\end{equation}\label{eq:Klambda}
where $\C_{MY}=\E[M(Z)\otimes Y]$ and $\C_{M}=\E[M(Z)\otimes M(Z)]$.
Since $\mathcal{T}\cong A_{T^*}:S_2(L_2(X),\mathbb{Y} )\to S_2(L_2(Z),\mathbb{Y} )$,
Equation (\ref{initial inverse problem}) rewrites
\begin{equation}\label{eq:T}
  r=A_{T^*} C
\end{equation}
where $r:=\Xi_Z^{-1} \widetilde{r}$ and where $\Xi_Z$ is a isometric isomorphism between $S_2(L_2(Z),\mathbb{Y} )$ and $L_2(Z,\mathbb{Y})$.
It then becomes a problem of finding a solution in $S_2(L_2(X),\Y)$.
The proposed estimation method tries to find a solution to this new problem within the space $S_2(\Hkx,\Y)$.
That is, it tries to solve
\begin{equation}\label{eq:Ttilde}
  r=A_{(T\mathcal{I}_X)^*}K
\end{equation}
where we used the fact that $\mathbf{I_X}=A_{\mathcal{I_X}^*}$,  with $\mathbf{I_X}:S_2(\Hkx,\Y)\to S_2(L_2(X),\Y)$  the embedding operator described in Appendix \ref{app:RKHS}.
Setting $\widetilde{T}:=T\mathcal{I}_X$,
a regularisation candidate for (\ref{eq:Ttilde}) is $g_\lambda(A_{\tilde{T}^*\tilde{T}})A_{\widetilde{T}}r$.
Since $\widetilde{T}$ corresponds to the operator
$$\widetilde{T}:\Hkx\to L_2(Z):f \mapsto \E[\langle f,\Phi_X(X)\rangle_{\mathcal{H}_{k_X}} |Z],$$
one can verify that $K_\lambda = g_\lambda(A_{\tilde{T}^*\tilde{T}})A_{\widetilde{T}}r$.
Now equivalently, one could try to solve (\ref{eq:T}) but under the additional assumption that the solution must have a representative in $S_2(\Hkx,\Y)$.
That is, we wish to find a regularized solution $C_\lambda\in S_2(L_2(X),\Y)$ such that 
$$\left\lVert C_\lambda\right\rVert_{S_2([\Hkx]^1,\Y)}=\left\lVert C_\lambda \mathcal{L}^{-1/2}_X\right\rVert_{S_2(L_2(X),\Y)}=\left\lVert A_{\mathcal{L}^{-1/2}_X}C_\lambda \right\rVert_{S_2(L_2(X),\Y)}<+\infty.$$
This is precisely the setting of regularisation in Hilbert scales mentioned before.
Here, the Hilbert scale is generated by the operator $A_{\mathcal{L}_X^{-1}}$ which is a densely defined, self-adjoint, positive and coercive unbounded operator as $\mathcal{L}_X^{-1}$ is such an unbounded operator (see Appendix \ref{app:AC} together with Lemma \ref{lem: properties unbounded operator}).
The coercivity holds with constant $\lVert \mathcal{L}_X \rVert^{-1} $.
Define the candidate solution
\begin{equation}\label{eq:Clambda}
  C_\lambda:=A_{\mathcal{L}_X^{1/2}}g_\lambda\left(A_{\mathcal{L}_X^{1/2}}A_{T}A_{T^*}A_{\mathcal{L}_X^{1/2}}\right)A_{\mathcal{L}_X^{1/2}}A_{T}r.
\end{equation}
The following proposition shows that these two approaches coincide.
\begin{proposition}\label{prop:problem is hilbert scaled}
  Let $K_\lambda$ and $C_\lambda$ be defined as in (\ref{eq:Klambda}) and (\ref{eq:Clambda}).
  It holds that $$\mathbf{I_X}K_\lambda = C_\lambda \qquad \text{in} \qquad S_2(L_2(X),\Y).$$
\end{proposition}
This proposition will enable us to use the machinery from the theory of regularisation in Hilbert scales to study the deterministic part of our problem.

\section{Selected solution: existence and characterisation}\label{sec: solution selection}
In this section we describe the limit solution towards which  $C_\lambda$ converges in the well- and misspecified settings.
The operator $Q:=A_{\mathcal{L}_X^{1/2}T^*}$ will play an important role in the description.
Denote by $\mathcal{S}:=\{C\in S_2(L_2(X),\Y):A_{T^*}C=r\}$ the set of solutions. In the sequel, $\Xi$ denotes be the isometric isomorphism between $S_2(L_2(X),\mathbb{Y})$ and $L_2(X,\mathbb{Y})$.
Note that $\mathcal{S}$ is nonempty since it must at least contain $\Xi^{-1}(F_0)$.
We first define properly the well- and misspecified settings.

\begin{definition}
The problem in \eqref{eq:Ttilde} is \emph{well-specified} if $\mathcal S\cap S_2([\mathcal H_{k_X}]^1,\mathbb Y)\neq\emptyset$,
and \emph{misspecified} otherwise.
\end{definition}

\subsection{Well-specified setting}
Let $C\in \mathcal{S}\cap S_2([\Hkx]^1,\Y)$, then $\mathcal{S}=C + \N(A_{T^*})$ with $\lVert A_{\mathcal{L}_X^{-1/2}}C\rVert_{S_2(L_2(X),\Y)}<+\infty$.
One then obtains that
\begin{equation*}
  \begin{split}
    C_\lambda=A_{\mathcal{L}_X^{1/2}}g_\lambda(Q^*Q)Q^*QA_{\mathcal{L}_X^{-1/2}}C \underset{\lambda \to 0}{\longrightarrow }A_{\mathcal{L}_X^{1/2}}Q^\dagger Q A_{\mathcal{L}_X^{-1/2}}C&=A_{\mathcal{L}_X^{1/2}}P_{\N(Q)^\perp} A_{\mathcal{L}_X^{-1/2}}C\\
    &=A_{\big(\mathcal{L}_X^{-1/2}P_{\mathcal{N}(T\mathcal{L}_X^{1/2})^\perp}\mathcal{L}_X^{1/2}\big)}C=:\Pi C.
  \end{split}
\end{equation*}
where we used the fact that $$\mathcal{N}(A_{\mathcal{L}_X^{1/2}T^*})^\perp\cong \mathcal{N}(T\mathcal{L}_X^{1/2} \otimes Id_\mathbb{Y} )^\perp=\mathcal{N}(T\mathcal{L}_X^{1/2})^\perp \otimes \mathbb{Y}, $$
together with Lemma \ref{lem: tensor projection}.
This is the solution described in \cite{meunier2024nonparametric}. 
For any two solutions $C_1,C_2\in \mathcal{S}\cap S_2([\Hkx]^1,\Y)$, we must have $\Pi C_1=\Pi C_2$ by uniqueness of the limit.
The following can be said about the operator $\Pi$. 
\begin{proposition}\label{prop:proj_core}
$\Pi$ is an $S_2([\Hkx]^1,\Y)$-orthogonal projection with $$\N(\Pi)=\N(A_{T^*})\cap S_2([\Hkx]^1,\Y) \qquad\text{and}\qquad \mathcal{R}(\Pi)=(\N(A_{T^*})\cap S_2([\Hkx]^1,\Y))^{\perp_{S_2([\Hkx]^1,\Y)}},$$
the orthogonal complement being taken in $\langle\cdot,\cdot\rangle_{S_2([\Hkx]^1,\Y)}$.
Furthermore, for any $C\in  S_2([\Hkx]^1,\Y)$, we have that $$A_{\mathcal{L}_X^{1/2}}g_\lambda(Q^*Q)Q^*A_{T^*}C \to \Pi C$$ in $S_2(L^2(X),\Y)$ as $\lambda \to 0$.
\end{proposition}
We now characterize the solution towards which we will show convergence.
\begin{proposition}\label{prop: X_s minimal core}
  We have that $\Pi(\mathcal{S}\cap S_2([\Hkx]^1,\Y))$ contains a unique element $C^{\dagger}$ which is the solution of minimal $S_2([\Hkx]^1,\Y)$-norm.
  That is, $C^{\dagger}$ satisfies $C^{\dagger} =\arg \min \{\lVert C \rVert_{S_2([\Hkx]^1,\Y)} : C\in\mathcal{S}\cap S_2([\Hkx]^1,\Y)\}$.
\end{proposition}
We now describe the solution in the misspecified setting.
\subsection{Misspecified setting}
It becomes clear from the well-specified setting that the first step to have a well-defined solution in the misspecified setting is to ensure that $\Pi$ extends to $S_2(L_2(X),\Y)$ as a bounded operator.
\begin{assumption}\label{ass:angle core}
  The operator $\Pi:S_2([\Hkx]^1,\Y) \to S_2(L_2(X),\Y)$ extends to a bounded operator $\overline{\Pi}:S_2(L_2(X),\Y) \to S_2(L_2(X),\Y)$.
\end{assumption}
This assumption is satisfied in the following cases.
\begin{proposition}\label{prop:sufficient core}
  Each of the following implies~\ref{ass:angle core}:
  \begin{enumerate}[label=(\roman*)]
  \item \textbf{Commutation.} $P_\Nt$ commutes with the spectral measure of $\mathcal{L}_X^{-1/2}$. In that case $\overline{\Pi}=A_{P_{\Nt^\perp}}$.
  \item \textbf{Smooth finite-dimensional null space.} 
    $T$ satisfies $dim(\mathcal{N}(T))<\infty$ and $\Nt \subseteq [\Hkx]^2$.
  \item \textbf{Two-sided polynomial link condition.} $T$ and $\mathcal{L}_X^{1/2}$ satisfy
  $$\hspace{1.75cm}\lVert T x \rVert_{L_2(Z)} \ge c_1\,  \lVert \mathcal{L}_X^{b/2}x\rVert_{L_2(X)} \qquad \text{for all } x\in \mathcal{L}_X^{1/2}(\N(T \mathcal{L}_X^{1/2})^\perp)$$and
  $$\lVert T x \rVert_{L_2(Z)} \le c_2\,  \lVert \mathcal{L}_X^{b/2} x\rVert_{L_2(X)} \qquad \text{for all } x\in L_2(X)$$
  for some $b>0$ and $0 < c_1\leq c_2$.
  \end{enumerate}
\end{proposition}
Notice that the first two conditions are satisfied when $T$ is injective.
We now characterize the operator $\overline{\Pi}$.
\begin{proposition}\label{prop:pi extension core}
  Under Assumption \ref{ass:angle core}, $\overline{\Pi}$ is a bounded idempotent operator on $S_2(L_2(X),\Y)$ with
  $$\N(\overline{\Pi})=\overline{\N(A_{T^*}) \cap S_2([\Hkx]^1,\Y)}, \hspace{0.5cm} \Ra(\overline{\Pi})=\overline{(\N(A_{T^*})\cap S_2([\Hkx]^1,\Y))^{\perp_{S_2([\Hkx]^1,\Y)}}} \hspace{0.5cm} $$ and $A_{T^*}\overline{\Pi}=A_{T^*}.$
\end{proposition}
One might still wonder whether $\overline{\Pi}(\mathcal{S})$ contains a unique element.
It turns out that under general link conditions this is indeed the case.
Let $\phi:\R^+ \to \R^+$ be a continuous, strictly increasing function with $\phi(0)=0$.
We call $\phi$ an index function (see \cite{mathe2003geometry}). 
The uniqueness is then given by the following proposition.
\begin{proposition}\label{prop: uniqueness}
  Assume there exists some index fuction $\phi$  and constant $m>0$ such that
  $$\lVert Tx\rVert_{L_2(Z)}  \ge c_1\,  \lVert [\phi(\mathcal{L}_X^{2})]^{1/2}x\rVert_{L_2(X)} \qquad \text{for all } x\in \mathcal{L}_X^{1/2}(\mathcal{N}(T\mathcal{L}_X^{1/2})^\perp).$$
  Then, under assumption \ref{ass:angle core} we have that $\overline{\Pi}(\mathcal{S})$ contains a unique element.
\end{proposition} 
These link conditions are standard assumptions in the literature and are indispensable to obtain convergence rates.

\section{Convergence rates}\label{sec: convergence rates}
Let $C^\dagger$ be the solution of minimal $S_2([\mathcal{H}_{k_X}]^1,\Y)$-norm (Proposition~\ref{prop: X_s minimal core}) if the problem is well-specified, and the unique element of $\overline{\Pi}(\mathcal{S})$ (Propositions~\ref{prop: uniqueness}) otherwise. 
If the problem is well-specified and \ref{ass:angle core} holds, both definitions agree, because $\overline{\Pi}$ extends $\Pi$. We set $F^\dagger:=\Xi C^\dagger$.
Denote by $\{\mu_{X,i}\}_{i\in I_X}$ and $\{\mu_{Z,i}\}_{i\in I_Z}$ the positive eigenvalues of $\mathcal{L}_X$ and $\mathcal{L}_Z$.
The following assumptions are required to derive convergence rates of our estimator towards $F^\dagger=\Xi C^\dagger$.

\begin{assumption}\label{ass:stage1}
  \begin{enumerate}
    \item There exist constants $l_1>0$ and $p_Z\in (0,1]$ such that for all $i\in I_Z$,
    $$\mu_{Z,i}\leq l_1 i^{-1/p_Z}.$$
    \item \label{ass:embZ}There exists $\alpha_Z \in [p_Z,1]$ such that the inclusion map $\mathcal{I}_{Z}^{\alpha_Z,\infty}:[\mathcal{H}_{k_Z}]^{\alpha_Z}\to L_\infty(Z)$ satisfies, for some constant $A_Z>0$,
    $$\left\lVert \mathcal{I}_{Z}^{\alpha_Z,\infty}\right\rVert \leq A_Z. $$
    \item \label{ass:sourceZ}There exist $\beta_Z > \alpha_Z$ and a constant $B_Z\geq 0$ such that,
    $$\lVert M\rVert_{\beta_Z}\leq B_Z. $$
  \end{enumerate}
\end{assumption}
\begin{assumption}\label{ass:stage2}
  \begin{enumerate}
    \item \label{ass: px}There exist constants $l_2>0$ and $p_X\in (0,1]$ such that for all $i\in I_X$,
    $$\mu_{X,i}\leq l_2 i^{-1/p_X}.$$
    \item \label{ass:embX} There exists $\alpha_X \in [p_X,1]$ such that the inclusion map $\mathcal{I}_{X}^{\alpha_X,\infty}:[\mathcal{H}_{k_X}]^{\alpha_X}\to L_\infty(X)$ satisfies, for some constant $A_X>0$,
    $$\left\lVert \mathcal{I}_{X}^{\alpha_X,\infty}\right\rVert \leq A_X. $$
    \item \label{ass: source}There exist $\beta_X > \alpha_X$ and a constant $B_X\geq 0$ such that,
    $$\lVert F^\dagger\rVert_{\beta_X}\leq B_X. $$
    \item There are constants $\sigma,L>0$ such that for all $k \geq 2$
    $$\E\left[\left\lVert Y-\E[F^\dagger(X)|Z]\right\rVert^k_{\mathbb{Y} } |Z\right] \leq \frac{1}{2}k!\sigma^2 L^{k-2}.$$
  \end{enumerate}
\end{assumption}
\begin{assumption}\label{ass: link core}
    There exists some index function $\phi$ and constants $c_2\geq c_1>0$ such that
  $$\lVert Tx\rVert_{L_2(Z)}  \ge c_1\,  \lVert [\phi(\mathcal{L}_X^{2})]^{1/2}x\rVert_{L_2(X)} \quad \text{for all } x\in \mathcal{L}_X^{1/2}(\N(T \mathcal{L}_X^{1/2})^\perp)$$
  and
  $$\lVert Tx\rVert_{L_2(Z)}  \le c_2\,  \lVert [\phi(\mathcal{L}_X^{2})]^{1/2}x\rVert_{L_2(X)} \quad \text{for all } x\in L^2(X).$$
  Furthermore, the function $\psi_{1/2}:(0,\lVert \mathcal{L}_X\rVert ]\to \R^+:t\mapsto t\,\phi(t^{2})$ is convex.
\end{assumption}

Assumption~\ref{ass:stage1} deals with the first stage of the estimation and include standard requirements in the kernel regression literature.
See for example \cite{meunier2024optimal}, \cite{li2022optimal} and \cite{fischer2020sobolev}.
Note that assumption \ref{ass:stage1}.\ref{ass:embZ} implies that the norm of the embedding operator from $[\mathcal{H}_{K_Z}]^{\alpha_Z}$ to $L_\infty(Z,\mathcal{H}_{k_X})$ is bounded by $A_Z$ (see Lemma 4 in \cite{li2022optimal}).

Assumption~\ref{ass:stage2} deals with the second stage of the estimation.
As for assumption \ref{ass:stage1}.\ref{ass:embZ}, assumption \ref{ass:stage2}.\ref{ass:embX} ensures that the norm of the embedding operator from $[\mathcal{H}_{K_X}]^{\alpha_X}$ to $L_\infty(X,\mathbb{Y} )$ is bounded by $A_X$.

Finally, Assumption~\ref{ass: link core} imposes link conditions describing the relative smoothness of $T$ with respect to $\mathcal{L}_X$. The former (resp. latter) inequality is referred to as the lower (resp. upper) link condition.
\begin{remark}
  In \cite{meunier2024nonparametric}, in which the authors assume index functions of the form $\phi(t)=t^b$ and $\beta_X\geq 1$, the only link condition required is
  $$R \lVert \C_X^{\eta/2}f\rVert_{\mathcal{H}_{k_X}} \leq \lVert \C_M^{1/2}f\rVert_{\mathcal{H}_{k_X}},\qquad \forall f\in \mathcal{N}(\C_M)^\perp. $$
  Note that $\mathcal{I}_X(\mathcal{N}(\C_M)^\perp)=\mathcal{L}_X^{1/2}(\mathcal{N}(T\mathcal{L}_X^{1/2})^\perp)$,
  $\lVert C_M^{1/2}f\rVert _{\mathcal{H}_{k_X}}=\lVert T\mathcal{I}_X f\rVert _{L^2(Z)}$ and $\lVert \C_X^{\eta/2}f\rVert_{\mathcal{H}_{k_X}}=\mathcal \lVert \mathcal{L}_X^{(\eta-1)/2}\mathcal{I}_X f\rVert _{L^2(X)}$ $\forall f\in \mathcal{H}_{k_X}$,
so that the lower link condition in Assumption~\ref{ass: link core} corresponds to the one above if we set $\phi(t^2)=t^{\eta-1}$ and $c_1=R$.
  In our proofs below, for a general index function, the upper link condition is only required to control $\mathcal{N}_M(\lambda)=tr(\C_M(\C_M+\lambda)^{-1})$ when $\beta_X\geq 1$.
  We leave the question of whether the upper link condition can be dropped in this case for future research.
\end{remark}

In the well-specified setting we will need an additional assumption in order to exploit the additional smoothness.
\begin{assumption}\label{ass: additional convexity}
  The function $f:=\psi_{1/2}\circ (\cdot)^{1/(\beta_X-1)}$ is convex on $(0,\lVert \mathcal{L_X}\rVert^{\beta_X-1} ]$.
\end{assumption}
\begin{remark}
  For $\phi(t^2)=t^b$, this corresponds to imposing $\beta_X\leq 2 + b$.
  This assumption can already be found in Theorem 3 of \cite{meunier2024nonparametric}
  and allows one to use Heinz-Kato's theorem.
\end{remark}
The convergence result below works under a constraint on the first stage estimation. To this end, define the inequalities
\begin{equation}\label{eq: set of equations}
n \geq 8\tau c_\lambda \kappa^2_X \lambda^{-1}, \qquad \frac{J\sqrt{\tau}\lVert M-\widehat{M}\rVert_{\alpha_Z} }{\lambda \sqrt{n}}\leq \frac{1}{12},\qquad \frac{J\lVert M-\widehat{M}\rVert_{L^2(Z,\mathcal{H}_{K_X})} }{\lambda }\leq \frac{1}{12}, \end{equation}
$$
\lVert M-\widehat{M}\rVert_{L^2(Z,\mathcal{H}_{K_X})}\leq 1,\qquad \text{and}\qquad \,\lVert M-\widehat{M}\rVert_{\alpha_Z}\leq 1,$$
where $J$ is a constant only depending on $\kappa_X,\kappa_Z,\alpha_Z,\beta_Z,A_Z$ and $B_Z$ and $c_\lambda:=\ln\left(2e\mathcal{N}_M(\lambda)\frac{\lVert \C_M \rVert +\lambda}{\lVert \C_M \rVert}\right)$ and let 
$$S:=\{s \in   E_X^m  \times E_Z^m : \text{the set of inequalities in \eqref{eq: set of equations} are satisfied}\}.$$
Let $[\cdot]_X$ denote the embedding operator from $\mathcal{H}_{K_X}$ to $L_2(X,\mathbb{Y} )$.
\begin{theorem}
  \label{th:main rate}
  Let Assumptions \ref{ass: RKHS}, \ref{ass:stage1}.\ref{ass:embZ}, \ref{ass:stage2} and \ref{ass: link core} hold and set $\kappa:=\max\{\kappa_Z,\kappa_X+A_Z\}$. 
  In addition, let Assumption \ref{ass: additional convexity} (resp. \ref{ass:angle core}) hold if Assumption \ref{ass:stage2}.\ref{ass: source} holds with $\beta_X >1$ (resp. $\beta_X<1$).
  Furthermore, let $g_\lambda$ be a filter with qualification 1 satisfying Assumption \ref{ass: Lipschitz}.
  
  Then, for $\gamma\in [0,\min\{1,\beta_X\}]$ and $\tau\geq \ln(16)$, 
  we obtain with $\mathbb{P}$-probability not less than $1-10e^{-\tau}-\mathbb{P}_{XZ}^{\otimes m}(S^c)$ that
  \begin{equation*}
    \begin{split}
      \left\lVert [\widehat{F}_\lambda]_X-F^\dagger \right\rVert_{[\mathcal{H}_{K_X}]^\gamma}^2&\leq c \tau^2 \Biggl[\left[\psi^{-1}_{1/2}\left(\lambda\right)\right]^{\beta_X-\gamma}+\frac{\left[\psi^{-1}_{1/2}\left(\lambda\right)\right]^{1-\gamma-p_X}}{\lambda n}+\frac{\left[\psi^{-1}_{1/2}\left(\lambda\right)\right]^{1-\gamma}}{\lambda^2 n^2}\\
      &\hspace{1cm}+\lambda^{-2} \,\left[\psi^{-1}_{1/2}\left(\lambda\right)\right]^{a(1-\gamma)}\left([\psi_{1/2}^{-1}(\lambda)]^{(\beta_X-1)_-}+\frac{[\psi_{1/2}^{-1}(\lambda)]^{-p_X}}{\lambda n}+\frac{1}{(\lambda n)^2} \right) \times \\
      &\hspace{1cm}\left(\frac{\lVert M-\widehat{M}\rVert_{\alpha_Z}^2}{n} +\lVert M-\widehat{M}\rVert_{L_2(Z,\mathcal{H}_{k_X})}^2 \right)\Biggr],
    \end{split}
  \end{equation*}
  for some constant $c>0$ independent of $n,m, \lambda$ and $\tau$ and with $a=\mathbf{1}\{\mathcal{N}(T)=\{0\}\}$ and $(t)_-:=\min\{0,t\}$.
\end{theorem}

We now discuss an appropriate choice for $\lambda$ and $\varepsilon$ in terms of the sample sizes $n$ and $m$.
This choice is given in the next proposition.

\begin{proposition}\label{prop:rates}
Let the assumptions of Theorem~\ref{th:main rate} hold with $\gamma\in[0,\min\{1,\beta_X\})$.
In addition, let Assumption \ref{ass:stage1} hold and set $\varepsilon \varpropto (m^{-1/(\beta_Z + p_Z)})$.
Furthermore, let $g_\varepsilon$ be a filter with qualification greater than $\beta_Z/2$.
Assume $m:=\lceil \eta n^\theta \rceil $ for some constants $\eta,\theta >0$.
Define the constants
$$r(\theta):=\min\Big\{\frac{\theta\beta_Z}{\beta_Z+p_Z}, 1+\frac{\theta(\beta_Z-\alpha_Z)}{\beta_Z+p_Z}\Big\},\qquad
\nu:=a(1-\gamma)+(\beta_X-1)_-, \qquad \xi:=\beta_X-\gamma-\nu,$$
$\Lambda_n:=\ln(1/\psi_{1/2}^{-1}(1/n)),$ and functions
$$\Theta_1(t):=t^{\beta_X+p_X}\,\phi(t^2),\qquad \Theta_2(t):=t^{1+\frac{\xi}{2}}\,\phi(t^2).$$
Then, the following holds:
\begin{enumerate}
  \item Let $(0,\psi_{1/2}(\lVert \mathcal{L_X}\rVert) ]\ni \lambda_n\to 0$ be a sequence such that $n^{-r(\theta)}=o(\lambda_n^2)$ and $\ln(1/\psi_{1/2}^{-1}(\lambda_n))/\lambda_n = o(n)$, then
    $$\big\lVert [\hat F_{\lambda_n}]_X-F^\dagger\big\rVert ^2_{[\mathcal H_{K_X}]^{\gamma}}=O_{\mathbb P}\left(R_n(\lambda_n)\right),$$
    where $$R_n(\lambda_n):=\big[\psi^{-1}_{1/2}(\lambda_n)\big]^{\beta_X-\gamma}+\frac{\big[\psi^{-1}_{1/2}(\lambda_n)\big]^{1-\gamma-p_X}}{n\lambda_n}+\frac{\big[\psi^{-1}_{1/2}(\lambda_n)\big]^{\nu}}{n^{r(\theta)}\lambda_n^2}+\frac{[\psi_{1/2}^{-1}(\lambda_n)]^{\,a(1-\gamma)-p_X}}{n^{1+r(\theta)}\lambda_n^{3}}.$$
    Furthermore, for $n$ large enough, we have that
    $$R_n(\lambda_n)\gtrsim \max\left\{\big[\Theta_1^{-1}\big(n^{-1}\big)\big]^{\beta_X-\gamma},\;(\psi_{1/2}^{-1}(\Lambda_n/n))^{\beta_X-\gamma},\;\big[\Theta_2^{-1}\big(n^{-r(\theta)/2}\big)\big]^{\beta_X-\gamma}\right\}.$$
  \item  Define sequences $l_n$ and $\omega_n$ such that these satisfy $l_n:=1$ if $\xi>0$ or $l_n\to +\infty$ with $l_n=o(n^{r(\theta)})$ if $\xi=0$, and $\omega_n \to +\infty$ with $\omega_n \ln n=o(n)$.
        Furthermore, let $t_n^*$ be the unique $t\in(0,\lVert \mathcal{L_X}\rVert]$ such that $n\psi_{1/2}(t)=\omega_n \ln(1/t)$ for $n$ large enough.
        Then, for
        \begin{equation}\label{eq: lambda}
          \lambda_n^*:=\psi_{1/2}\left(\max\left\{\Theta_1^{-1}(n^{-1}),\;\; t_n^*,\;\;\Theta_2^{-1}\left(l_n^{1/2} n^{-r(\theta)/2}\right)\right\}\right),
        \end{equation}
        we obtain
        $$R_n(\lambda_n^*)\lesssim \big[\Theta_1^{-1}\big(n^{-1}\big)\big]^{\beta_X-\gamma}+(\psi_{1/2}^{-1}(\omega_n\Lambda_n/n))^{\beta_X-\gamma}+\big[\Theta_2^{-1}\big((\ell_n n^{-r(\theta)})^{1/2}\big)\big]^{\beta_X-\gamma},$$
        where $\lambda_n^*$ satisfies the conditions of point 1. above.
\end{enumerate}
\end{proposition}
  
To illustrate this result, we discuss two well-known examples.
\begin{example}[Polynomial link function: $\phi(t^2)=t^b,\;b>0$]
  Set $l_n=\omega_n=\ln n$.
  If $a=1$ and $\beta_X<1$, then $\xi=0$ and therefore $l_n=o(n^{r(\theta)})$.
  If in addition $r(\theta)>(2+2b)/(\beta_X+p_X+b)$ and $p_X+\beta_X>1$, then one can check that
  $$\big\lVert [\hat F_{\lambda_n}]_X-F^\dagger\big\rVert ^2_{[\mathcal H_{K_X}]^{\gamma}}=O_{\mathbb{P}}\Big(n^{-\frac{\beta_X-\gamma}{\beta_X+p_X+b}}\Big),$$
  which corresponds to the optimal rate in the well-specified setting presented in \cite{meunier2024nonparametric} (with $\gamma_0=\gamma_1$).
  We conjecture that this rate is still the optimal one for $\beta_X<1$.
\end{example}

\begin{example}[Exponential link condition: $\phi(t^2):=\exp(-t^{-q}),\;q>0$]
  Set $l_n=\ln n$ if $\xi =0$ and $\omega_n= \ln n$.
  For any function $h_k:\R^+\to \R^+:t\mapsto t^k\exp(-t^{-q})$ for $k\in \mathbb{R} $, one can show that $h_k^{-1}(t)=[\ln(1/t)]^{-1/q}(1+o(1))$ as $t\to 0$.
  Hence, $$\Theta^{-1}_1(n^{-1})=[\ln(n)]^{-1/q}(1+o(1)),$$and
  $$\Theta_2^{-1}\big((\ell_nn^{-r(\theta)})^{1/2}\big)=\Big[\tfrac{r(\theta)}{2}\ln n-\tfrac12\ln\ell_n\Big]^{-1/q}(1+o(1))=\Big(\tfrac{2}{r(\theta)}\Big)^{1/q}(\ln n)^{-1/q}(1+o(1)).$$

  Notice that for $n$ large, $\psi_{1/2}^{-1}(1/n)\ge(\ln n)^{-1/q}$, and hence
  $$\Lambda_n=\ln\big(1/\psi_{1/2}^{-1}(1/n)\big)\le\tfrac1q\ln\ln n.$$
  Therefore, $x_n:=\omega_n\Lambda_n/n \le(\ln n)(\ln\ln n)/(qn)\to0$.
  So $\psi_{1/2}^{-1}(x_n)$ is defined and at most $1$ for $n$ large. Moreover,
  $$\ln(1/x_n)=\ln n-\ln\ln n-\ln\Lambda_n\ \ge\ \ln n-\ln\ln n-\ln\big(\tfrac1q\ln\ln n\big)\ \ge\ \tfrac12\ln n$$
  for $n$ large. This gives
  $$\psi_{1/2}^{-1}(x_n)=\ln (1/x_n)^{-1/q}(1+o(1)) \le \Big(\frac{1}{2}\ln n\Big)^{-1/q}(1+o(1)).$$
  Plugging in $\lambda_n^\ast$ yields $\psi_{1/2}^{-1}(\lambda_n^*)=c (\ln n)^{-1/q}(1+o(1))$ for some $c>0$ and hence 
  $$\big\lVert [\hat F_{\lambda_n^*}]_X-F^\dagger\big\rVert ^2_{[\mathcal H_{K_X}]^{\gamma}}=O_{\mathbb{P}}\left((\ln n)^{(\gamma-\beta_X)/q}\right).$$
\end{example}
Now, if one only has access to one sample $\{X_i,Y_i,Z_i\}_{i=1}^N$, it has to be split into two disjoint samples of size $m$ and $n=N-m$.
For $0<\eta<1$, one could set $m=\left\lceil \eta N\right\rceil $. One can then recover the convergence rate in this setting by substituting $\theta=1$ in Proposition \ref{prop:rates}.
This gives the following rate.
\begin{corollary}\label{cor:theta1}
Let the assumptions of Proposition~\ref{prop:rates} hold with $\theta=1$, and let $\lambda_n^*$ be as in Proposition \ref{prop:rates}, with $\omega_n\ln n=o\big(n^{1-r(1)/2}\big)$. 
Then,
$$\big\lVert [\widehat{F}_\lambda]_X-F^\dagger \big\rVert_{[\mathcal{H}_{K_X}]^\gamma}^2=O_{\mathbb{P}}\left(\Big[\Theta_2^{-1}\Big(\big(\ell_n\,n^{-\frac{\beta_Z}{\beta_Z+p_Z}}\big)^{1/2}\Big)\Big]^{\beta_X-\gamma}\right).$$
\end{corollary}

\section{Conclusion}\label{sec: conclusions}
We have studied a two-stage kernel estimator for nonparametric instrumental variable regression
with outcomes in a separable Hilbert space $\mathbb Y$. The first stage estimates the conditional
mean embedding of $X$ given $Z$; the second stage regresses the $\mathbb Y$-valued outcome on this
embedding with a spectral filter. The estimator is linear in the outcomes: $\widehat F_\lambda(x)$ is
a weighted sum of $\widetilde Y_1,\dots,\widetilde Y_n$ whose weights depend only on the
observations of $X$ and $Z$. The space $\mathbb Y$ therefore enters only through the outcomes: it
does not affect how the weights are computed, and the rates take the same form as for scalar
outcomes.

Our analysis rests on the observation that the second stage is a regularisation method in the
Hilbert scale generated by $A_{\mathcal L_X^{-1}}$ (Proposition~1). This viewpoint separates two
questions that are usually treated together. 

First, we pinpoint the limit solution of the regularized problem in the set of solutions of the operator equation. 
When the conditional expectation operator $T$ is not injective and the
problem is well-specified, the target is the solution of minimal norm in
$S_2([\mathcal H_{k_X}]^1,\mathbb Y)$ (Proposition~3). When no solution lies in this space, the
target is the unique solution in the range of $\overline\Pi$, the continuous extension of the
projection that is orthogonal in $S_2([\mathcal H_{k_X}]^1,\mathbb Y)$ and whose kernel is the part
of the null space contained in this space (Propositions~5 and~6). 

Second, we provide convergence rates of the estimator to its limit.
Under a source condition and a link condition expressed through a general index function $\phi$,
Theorem~1 gives non-asymptotic bounds in the interpolation norms of $[\mathcal H_{k_X}]^\gamma$,
$0\le\gamma \le \min\{1,\beta_X\}$, holding with high probability over both samples. These bounds
apply with or without identification and in the well- and misspecified settings.

For polynomial link conditions, $\phi(t^2)=t^b$, and an injective $T$, the squared error in
$[\mathcal H_{k_X}]^\gamma$ is of order $n^{-(\beta_X-\gamma)/(\beta_X+p_X+b)}$ when the first-stage
sample is large enough (Example~1). In $L^2(X)$ this is the exponent obtained by \cite{meunier2024nonparametric} for scalar outcomes. 
Our result extends it to $\beta_X<1$ and to $\mathbb Y$-valued outcomes. For exponential
link conditions, which model severely ill-posed problems, the rate becomes
$(\ln n)^{-(\beta_X-\gamma)/q}$ (Example~2). The sample allocation matters: when a single sample is
split into two parts of proportional sizes, the first-stage error dominates the rate
(Corollary~1).

Several questions remain open. First, we have no lower bounds in the misspecified or severely
ill-posed regimes, so the optimality of our rates is established only where they match known
minimax rates; lower bounds depending on both sample sizes would also clarify how large the
first-stage sample must be. Second, our proofs use an upper link inequality, which controls the
effective dimension and the regularisation error in the misspecified case.
\cite{meunier2024nonparametric} avoid it for
$\beta_X\ge1$ through a mixed effective dimension, and extending this device to our setting would
weaken the assumptions. Finally, the regularisation parameters have to be chosen in practice. An
a posteriori rule, such as a balancing principle adapted to the two-stage structure, is left for future work.

\section*{Acknowledgements and Disclosure of Funding}
Guillaume Baquet gratefully thanks the Fonds de la recherche scientifique - FNRS for his funding through the FRIA grant ARDHIV, number FC 63567. 
Germain Van Bever gratefully acknowledges the support of the Fondation Roi Baudouin through its funding associated with the Thelam Fund.
\newpage
\bibliographystyle{plainnat}
\bibliography{references}
\newpage

\appendix
\textbf{\textit{The notations used in Appendices A, B, C and D should be seen as independent of the rest as some notations used in these sections are used elsewhere with a different meaning.}}

\section{Regularisation in Hilbert scales}\label{app:hilbertscalesreg}

In this section we recall in detail the construction of (positive) Hilbert scales based on the article \cite{nair2015role}.
Afterwards, we discuss the theory needed about regularisation in Hilbert scales.

\subsection{(Positive) Hilbert scales}
Let $H$ be a Hilbert space and $B:\D(B)\subset H \to H$ be a densely defined unbounded operator which is assumed to be positive, self-adjoint and coercive, that is, for some $c>0$
$$\langle Bx,x \rangle_H \geq c \norm{x}^2 \qquad \forall x\in \D(B).$$
Notice that coercivity implies that $\norm{Bx}_H\geq c \norm{x}_H$ for all $x\in \D(B)$.
By the spectral theorem, there exists a spectral family $\{E_\lambda\}$ satisfying (1) $E_\lambda \leq E_\mu$ for $\lambda \leq \mu$, (2) $E_\lambda E_\mu = E_{\lambda \wedge \mu}$, (3) $E_\lambda = \lim_{\mu \searrow \lambda}E_\mu$, 
and (4) $E(\lambda)=0, E(+\infty)=I$ for $\lambda < x$ such that $$B=\int_{c}^{+\infty}\lambda dE(\lambda).$$
From this representation we define the fractional powers of $B$ by
$$B^s:=\int_{c}^{+\infty}\lambda^s dE(\lambda)$$
with $\D(B^s)=\{x\in H: \int_{c}^{+\infty}\lambda^{2s}d\norm{E_\lambda x}^2 < +\infty \}$ for $s\in \mathbb{R}$.
Let us now define the scalar products $\langle u,v\rangle_s:=\langle B^s u , B^s v\rangle_H$ for $s>0$ and their associated norms $\norm{\cdot}_s$.
We start by showing the following about the spaces $(\D(B^s),\langle \cdot,\cdot \rangle_s)$.
\begin{proposition}
  For $s>0$, the space $(\D(B^s),\langle \cdot,\cdot \rangle_s)$ is a Hilbert space.
\end{proposition}
\begin{proof}
It is rather easy to prove that $\langle \cdot,\cdot \rangle_s$ satisfies the properties of a scalar product.
Hence, we will only show that the space is complete.
Since the adjoint of a densely defined operator is closed and $B$ is self-adjoint, we have that $B^s$ is closed.
Let $\{x_n\}$ be a Cauchy sequence in $(\D(B^s),\langle \cdot,\cdot \rangle_s)$, then $\{B^s x_n\}$ is a Cauchy sequence in $H$ and since $H$ is complete, there exists $y\in H$ such that $B^s x_n \to y$ in $H$.
Furthermore, since $$\lVert x_n-x_m\rVert_s^2 = \lVert B^s(x_n-x_m)\rVert_H^2 = \int_{c}^{+\infty}\lambda^{2s}d\norm{E_\lambda (x_n-x_m)}^2 \geq c^{2s}\norm{x_n-x_m}^2_H,$$
we have that $\{x_n\}$ is also a Cauchy sequence in $H$ and hence there exists $x\in H$ such that $x_n \to x$ in $H$.
By closedness of $B^s$, we have that $x\in \D(B^s)$ and $B^s x = y$.
Finally, since $$\norm{x-x_n}_s=\norm{y-B^sx_n}_H \to 0$$
as $n\to \infty$, we have that $x_n \to x$ in $(\D(B^s),\langle \cdot,\cdot \rangle_s)$ which finishes the proof.
\end{proof}
Let us now define the space $X=\bigcap_{s>0} \D(B^s)$.
This space has the following denseness property.
\begin{proposition}
  The space $X$ is dense in $\D(B^s)$ under the norm $\norm{\cdot}_s$ for all $s>0$.
\end{proposition}
\begin{proof}
Let $x\in \D(B^s)$ and $\epsilon>0$. We want to find $y\in X$ such that $\norm{x-y}_s<\epsilon$.
Since $x\in \D(B^s)$, we have that $\int_{c}^{+\infty}\lambda^{2s}d\norm{E_\lambda x}^2 < +\infty$.
Hence, there exists $R>c$ such that $\int_{R}^{+\infty}\lambda^{2s}d\norm{E_\lambda x}^2 < \epsilon^2$.
Define $y:=E_R \,x$, then
$$\norm{y}_r^2 = \int_{c}^{R}\lambda^{2r}d\norm{E_\lambda x}^2\leq R^{2r} \norm{x}_H^2 < +\infty.$$
for all $r>0$, which shows that $y\in \D(B^r)$ for all $r>0$ and therefore $y\in X$.
Finally, we have
$$\norm{x-y}_s^2 = \int_{R}^{+\infty}\lambda^{2s}d\norm{E_\lambda x}^2 < \epsilon^2.$$
which shows that $X$ is dense in $\D(B^s)$ under the norm $\norm{\cdot}_s$.
\end{proof}
Now since $B$ is bounded below, one can easily show that $B^s$ is also bounded below for all $s>0$.
Therefore, setting $X_0:= H$ and $X_s := \D(B^s)$ for $s>0$, we have that $X_s$ is dense in $X_t$ for $s > t$.
One can also define the spaces $X_s$ for $s<0$ using Gelfand triples (see \cite{nair2015role}).

\subsection{Regularisation}
Let $\{X_s\}_{s\in \mathbb{R}}$ be a Hilbert scale associated to a densely defined, self-adjoint, non-negative and coercive unbounded operator $B$ on a Hilbert space $X$.
Given an ill-posed problem $Tx=y$ with $T:X\to Y$ a bounded not necessarily injective linear operator and $y\in \Ra(T)$, we will be interested in describing the convergence of the element 
\begin{equation}\label{eq:xlambda}
x_\lambda:=B^{-s}g_\lambda(A^*A)A^*y
\end{equation}
where $A:=TB^{-s}$ for some $s>0$ and 
$g_\lambda:\R^+\to\R^+$ is a filter function, satisfying the following assumptions.
\begin{assumption}\label{ass:filter-g}
There exists positive constants $\rho_0$ and $\rho_1\leq 1$ such that 
\begin{enumerate}[label=\textup{(G\arabic*)}]
\item $\sup_{t\in [0,\norm{A^*A}]}t^{1/2}|g_\lambda(t)|\le\rho_0 \lambda^{-1/2}$;
\item $\sup_{t\in [0,\norm{A^*A}]}t|g_\lambda(t)|\le \rho_1$
\end{enumerate}
\end{assumption}
\begin{assumption}\label{ass:filter}
There exists positive constants $\gamma_0,\gamma_1$ and $\gamma_2$ such that 
\begin{enumerate}[label=\textup{(R\arabic*)}]
\item $\sup_{t\in [0,\norm{A^*A}]}|r_\lambda(t)|\le\gamma_0$;\label{F:r}
\item $\sup_{t\in [0,\norm{A^*A}]}t^{1/2}|r_\lambda(t)|\le\gamma_1\lambda^{1/2}$
\item $\sup_{t\in [0,\norm{A^*A}]}t|r_\lambda(t)|\le\gamma_2\lambda$.
\end{enumerate}
where $r_\lambda(t):=1-tg_\lambda(t)$.
\end{assumption}
In the following two sections, we will describe the convergence of $x_\lambda$ in the well- and misspecified setting.

\subsubsection{Well-specified setting}
Let $\mathcal{S}:=T^{-1}(\{y\})$ be the set of solutions.
If $\mathcal{S} \cap X_s\neq \varnothing$, then there exists $x^*\in \mathcal{S}\cap X_s$ such that $\mathcal{S}=x^*+\N(T)$.
Then, one can write $B^s x^*$ as an element in $X$ and therefore
\begin{equation}\label{eq:plain-conv}
  x_\lambda=B^{-s}g_\lambda(A^*A)A^*AB^s x^*\underset{\lambda \to 0}{\to} B^{-s}A^\dagger A B^s x^*=B^{-s}P_{\N(A)^\perp} B^s x^*=:\Pi x^*.
\end{equation}
For any two solutions $x_1^*,x_2^*\in \mathcal{S}\cap X_s$, we must have $\Pi x_1^*=\Pi x_2^*$ by uniqueness of the limit.
The following can be said about the operator $\Pi$.

\begin{proposition}\label{prop:proj}
  $\Pi$ is an $X_s$-orthogonal projection with $$\N(\Pi)=\Nt\cap X_s \qquad\text{and}\qquad \Ra(\Pi)=(\Nt\cap X_s)^{\perp_s},$$
  the orthogonal complement being taken in $\langle\cdot,\cdot\rangle_s$.
  Furthermore, if \ref{F:r} holds, then for any $x\in X_s$ we have $R_\lambda Tx:=B^{-s}g_\lambda(A^*A)A^*Tx \to \Pi x $ in $X$ as $\lambda \to 0$.
\end{proposition}

\begin{proof}
  Define the operator $U:X_s \to X:x\mapsto B^s x$.
  Since $\langle B^s x,B^s u\rangle_X = \langle x,u\rangle_s$ for all $x,u\in X_s$ we have that $U$ is an isometry which preserves the inner product.
  Also, $B^s$ is a closed operator as self-adjoint operators are closed.
  Furthermore, it is also bounded below by $c^{s}$.
  Letting $B^sx_n \to y$ in $X$ for some sequence $(x_n)\subset \D(B^s)$, we have that
  $$\norm{x_n-x_m}_X \leq c^{-s} \norm{B^sx_n - B^sx_m}_X \to 0 \quad \text{as } n,m \to \infty,$$
  since $(B^sx_n)$ is a Cauchy sequence in $X$.
  Thus, $(x_n)$ is a Cauchy sequence in $X$ and by completeness of $X$, there exists $x \in X$ such that $x_n \to x$ in $X$.
  Since $B^s$ is closed, we have that $x \in \D(B^s) = X_s$ and $B^sx = y$ in $X$. This shows that $B^s$ has closed range.
  Now since $B^s$ is injective, self-adjoint, and has closed range, we easily obtain that $B^s$ is surjective onto $X$.
  Indeed, we then obtain $$\Ra(B^s)=\overline{\Ra(B^s)}=\mathcal{N}((B^s)^*)^\perp=\mathcal{N}(B^s)^\perp=X.$$
  Hence, $U$ is an isometric isomorphism between $X_s$ and $X$ with inverse and adjoint both equal to 
  $U^*:X \to X_s:x\mapsto B^{-s}x$ such that $UU^*=I_X$ and $U^*U=I_{X_s}$. 

  Recall the notation $A=TB^{-s}$, clearly for $z\in X$ $$z\in\N(A) \iff B^{-s}z\in\N(T) \iff B^{-s}z\in\N(T)\cap X_s,$$
  since $B^{-s}z\in X_s$ for all $z\in X$. From this, one obtains that $\N(T)\cap X_s =U^*(\N(A))$.
  Since $U^*$ is a bijection and $\N(A)$ is closed in $X$, we have that $\N(T)\cap X_s$ is closed in $X_s$.
  Set $Q:=I-P_{\N(A)^\perp}$ the $X$-orthogonal projector onto $\Na$. One easily checks that
  $$\Pi^2=U^*QUU^*QU = U^*QU = \Pi \qquad\text{and}\qquad\Pi^*=(U^*QU)^*=U^*QU=\Pi.$$
  Hence, $\Pi$ is an $X_s$-orthogonal projection on $\mathcal{R}(\Pi)$.
  Now, since
  $$\mathcal{R}(I_{X_s}-\Pi) = (U^*P_{\Na}U) (X_s)=(U^*P_{\Na})(X)=U^*(\Na)=\N(T)\cap X_s,$$
  and $\mathcal{N}(\Pi) = \mathcal{R}(I_{X_s}-\Pi)$, we obtain obtain the identities concerning the kernel and the range.

  Now let $x\in \D(B^s)$, then $x=B^{-s}v$ for some $v\in X$. Then,
  \begin{equation*}
      \left\lVert R_\lambda Tx-\Pi x\right\rVert^2 = \left\lVert B^{-s}(r_\lambda(A^*A)-P_{\N(A^*A)})v\right\rVert^2 \leq c^{-2s}\int_{\sigma(A^*A)}(r_\lambda(t)-\mathbb{I}_{\{0\}}(t))^2\,d\norm{E(t)v}^2\underset{\lambda \to 0}{\to}0
  \end{equation*}
  where we applied the dominated convergence theorem since $r_\lambda(0)=1$, $r_\lambda(t)\to 0$ as $\lambda \to 0$ for all $t>0$, and
  $(r_\lambda(t)-\mathbb{I}_{\{0\}}(t))^2 \leq (\gamma_0+1)^2$ which is integrable.
\end{proof}
We now characterize the solution towards which $x_\lambda$ converges.
\begin{lemma}\label{lem: X_s minimal}
  We have that $\Pi(\mathcal{S}\cap X_s)$ contains a unique element $x^{\dagger_s}$ which is the solution of minimal $X_s$-norm.
  That is, $x^{\dagger_s}$ satisfies $x^{\dagger_s} =\arg \min \{\lVert x \rVert_s : x\in\mathcal{S}\cap X_s\}$.
\end{lemma}
\begin{proof}
  Let $x_1,x_2\in\mathcal{S}\cap X_s$. Then, $x_1-x_2\in \N(T)\cap X_s$ and hence $\Pi x_1-\Pi x_2 = \Pi(x_1-x_2)=0$ by Proposition \ref{prop:proj}.
  Therefore, $\Pi(\mathcal{S}\cap X_s)$ contains a unique element $x^{\dagger_s}$.
  Furthermore, for any $x\in \mathcal{S}\cap X_s$ one has
  $$T(I_{X_s}-\Pi)x =TB^{-s}P_{\N(A)}B^sx =AP_{\N(A)}B^s x= 0.$$
  Hence, $x^{\dagger_s}=\Pi x \in \mathcal{S}\cap X_s$.
  Now, since $\mathcal{S}=x^{\dagger_s}+\N(T)$, we obtain $\mathcal{S}\cap X_s=x^{\dagger_s}+(\N(T)\cap X_s)$.
  Therefore, for any $x\in \mathcal{S}\cap X_s$ we have by Pythagoras that
  $$\lVert x\rVert_s^2= \lVert \Pi x\rVert_s^2 +\lVert (I_{X_s}-\Pi)x\rVert_s^2 \geq \lVert \Pi x\rVert_s^2 = \lVert x^{\dagger_s}\rVert_s^2$$
\end{proof}
In order to obtain convergence rates, we will need to assume so-called link and source conditions.
The source condition describes the smoothness of $x^{\dagger_s}$, whereas
the link condition describes the relative smoothness of $T$ with respect to $B$.
The link condition uses an index function, that is, a continuous, strictly increasing function $\phi:\R^+\to\R^+$ with $\phi(0^+)=0$.
The conditions in the well-specified setting are the following
\begin{assumption}\label{ass:source-well}
  There exist constants $p\geq s$ and $E>0$ such that $\lVert B^{p}x^{\dagger_s}\rVert \leq E$.
\end{assumption}
\begin{assumption}\label{ass:link-well}
  There exists some index fuction $\phi$ with the following properties.
  \begin{enumerate}[label=(\roman*)]
  \item There exists a constant $0<m$ such that
  $$\norm{T x} \ge m\,  \lVert [\phi(B^{-2})]^{1/2}x\rVert \qquad \text{for all } x\in B^{-s}(\Na^\perp).$$
  \item The function $\psi_s:(0,c^{-2s}]\to \R^+:t\mapsto t\,\phi(t^{1/s})$ is convex.
  \end{enumerate}
\end{assumption}
 
For the case $p=s$, where $p$ is the parameter displayed in Assumption \ref{ass:source-well}, Assumptions \ref{ass:source-well} and \ref{ass:link-well} are sufficient in order to obtain convergence rates.
For the case $p>s$, we need to add the following assumption.
\begin{assumption}\label{ass:supplement}
  The function $f:(0,c^{-2(p-s)}]\to \R^+:t\mapsto t^{s/(p-s)}\,\phi(t^{1/(p-s)})$ is convex.
\end{assumption}
Before we can give the convergence rates, we will need the following results.
\begin{lemma}\label{lem:interpolation-high}
  Let $W$ be a subspace of $X$.
  Assume the following holds:
  $$\norm{T x} \ge m\,  \lVert [\phi(B^{-2})]^{1/2}x\rVert \qquad \text{for all } x\in W,$$
  with $m>0$ and $\phi$ an index function such that $\psi_a:(0,c^{-2a}]\to \R^+:t\mapsto t\,\phi(t^{1/a})$ is convex.
  Then, for any $y\in W\cap X_a$ and $\tau \in [0,a]$, we have
  $$\norm{y}_\tau \leq \norm{y}_a \left[ \psi^{-1}_a \left( \frac{\norm{Ty}^2}{m^2\norm{y}^2_a}\right)\right]^{\frac{a-\tau}{2a}},$$
  where the function $$\R^+_0 \times \R^+_0 \to \R^+:(u,t)\mapsto t \left[ \psi^{-1}_a \left( \frac{u^2}{m^2 t^2}\right)\right]^{\frac{a-\tau}{2a}}$$
  is nondecreasing in $u$ and $t$.
\end{lemma}
\begin{proof}
  Let $E_\lambda$ be the spectral family of $B$ and let $y\in W$. Then, by Jensen's inequality we obtain,
  \begin{equation*}
    \begin{split}
      \psi_a\left(\left(\frac{\norm{y}^2_\tau}{\norm{y}^2_a}\right)^{a/(a-\tau)}\right)&=\psi_a\left(\left(\frac{\int \lambda^{2\tau}\,d\norm{E_\lambda y}^2}{\int \lambda^{2a}\,d\norm{E_\lambda y}^2}\right)^{a/(a-\tau)}\right)\\
      & \leq \norm{y}_a^{-2}\int \psi_a\Big((\lambda^{2(\tau-a)})^{a/(a-\tau)}\Big) \lambda^{2a}\,d\norm{E_\lambda y}^2 \\
      &= \norm{y}_a^{-2}\int \phi(\lambda^{-2})\,d\norm{E_\lambda y}^2\\
      &=\frac{\norm{[\phi(B^{-2})]^{1/2}y}^2}{\norm{y}_a^2} \leq \frac{\norm{Ty}^2}{m^2 \norm{y}_a^2}
    \end{split}
  \end{equation*}
  The inequality now follows by applying the function $\psi_a^{-1}$, which is nondecreasing, to both sides and then rearranging the terms.
  For the second part, the nondecreasingness in $s$ is trivial. For $t$, set $u_s:=u_s(t):=\psi_a^{-1}(s^2/(m^2 t^2))$ which is decreasing in $t$ for fixed $s$.
  Then, $t=(s^2/(\psi_a(u_s)m^2))^{1/2}$ and therefore
  $$t \left[ \psi^{-1}_a \left( \frac{s^2}{m^2 t^2}\right)\right]^{\frac{a-\tau}{2a}}=\frac{s}{m}\left(\frac{u_s^{(a-\tau)/a}}{\psi_a(u_s)}\right)^{1/2}=\frac{s}{m}\left(\frac{1}{u_s^{\tau/a}\phi(u_s^{1/a})}\right)^{1/2}$$
  which is now clearly increasing in $t$ for fixed $s$.
\end{proof}

\begin{lemma}\label{lem:invariance}
  Let $f$ be a bounded Borel measurable function on the spectrum of $A^*A$.
  Then, for any $x\in \Na^\perp$, we have $f(A^*A)x \in \Na^\perp$.
\end{lemma}
\begin{proof}
  Let $v\in \Na =\N(A^*A)$, then 
  $$\langle f(A^*A)x,v \rangle =\langle f(A^*A)x,E_A(0)v \rangle= \langle E_A(0) x,f(A^*A)v \rangle = \langle P_{\Na} x,f(A^*A)v \rangle =0$$
  where we used the fact that $f(A^*A)$ commutes with $E_A$, the spectral measure of $A^*A$.
\end{proof}

\begin{lemma}\label{lem:Naperp}
  We have $\N(A)^\perp = B^s(\Ra(\Pi))$.
\end{lemma}
\begin{proof}
 By the proof of Proposition \ref{prop:proj}, we have that $\N(A)=B^s(\N(T)\cap X_s)$
 and since $B^s:X_s \to X$ is an isometric isomorphism, we obtain $\N(A)^\perp=B^s((\N(T)\cap X_s)^{\perp_s})=B^s(\Ra(\Pi))$.
\end{proof}
We can now demonstrate the convergence rates.
The proof technique is almost identical to what is presented in Proposition 5.1 in \cite{nair2005regularization}.
\begin{proposition}\label{prop: well deterministic rate}
  Let $x_\lambda$ be the element defined in \eqref{eq:xlambda} with $s>0$, where $g_\lambda$ is a filter function satisfying Assumption \ref{ass:filter}.
  Then, under Assumptions \ref{ass:source-well} and \ref{ass:link-well}, we obtain for $p=s$ the upper bound
  $$\lVert x_\lambda - x^{\dagger_s}\rVert_\tau \leq \gamma_0 E \left[\psi_s^{-1}\left(\frac{\gamma_1^2 }{m^2 \gamma_0^2}\lambda\right)\right]^{\frac{p-\tau}{2s}}\qquad \text{for }\tau \in [0,p].$$
  If in addition Assumption \ref{ass:supplement} holds, then for $p>s$ we obtain the upper bound
  $$\lVert x_\lambda - x^{\dagger_s}\rVert_\tau \leq E \left[\psi_s^{-1}\left(\frac{\gamma_2 }{m^2}\lambda\right)\right]^{\frac{p-\tau}{2s}}\qquad \text{for }\tau \in [0,p].$$
\end{proposition}
\begin{proof}
  We start by the case $p=s$.
  First notice that
  $$x_\lambda - x^{\dagger_s} = B^{-s}g_\lambda(A^*A)A^*y - x^{\dagger_s} = B^{-s}r_\lambda(A^*A)B^s  x^{\dagger_s}.$$
  Furthermore, by Lemmas \ref{lem:invariance} and \ref{lem:Naperp} we have that $r_\lambda(A^*A)B^sx^{\dagger_s}\in \Na^\perp$.
  Therefore, by Lemma \ref{lem:interpolation-high}, we obtain
  $$\lVert x_\lambda - x^{\dagger_s}\rVert_\tau \leq \lVert B^{-s}r_\lambda(A^*A)B^s x^{\dagger_s}\rVert_s \left[\psi_s^{-1}\left(\frac{\lVert Ar_\lambda(A^*A)B^s x^{\dagger_s}\rVert^2}{m^2\norm{B^{-s}r_\lambda(A^*A)B^s x^{\dagger_s}}^2_s}\right)\right]^{\frac{s-\tau}{2s}}.$$
  Since by Assumptions \ref{ass:filter} and \ref{ass:source-well}, we have that $$\lVert B^{-s}r_\lambda(A^*A)B^s x^{\dagger_s}\rVert_s^2\leq \gamma_0^2 E^2 \qquad \text{and} \qquad\lVert Ar_\lambda(A^*A)B^s x^{\dagger_s}\rVert^2\leq \gamma_1^2\lambda E^2,$$
  we obtain the desired result for $p=s$ by application of Lemma \ref{lem:interpolation-high}.

  Now for the case $p>s$, we have by the proof of Proposition 5.1 in \cite{nair2005regularization} the inequalities
  $$\lVert x_\lambda - x^{\dagger_s}\rVert_s \leq E\left[\psi^{-1}_s\left(\frac{\gamma_2}{m^2}\lambda\right)\right]^{(p-s)/(2s)} \hspace{0.2cm}\text{and}\hspace{0.4cm} \lVert T(x_\lambda - x^{\dagger_s})\rVert_s \leq \sqrt{\gamma_2 \lambda}E\left[\psi^{-1}_s\left(\frac{\gamma_2}{m^2}\lambda\right)\right]^{(p-s)/(2s)}.$$
  Again using $x_\lambda - x^{\dagger_s} \in B^{-s}(\N(A)^\perp)$ together with Lemma \ref{lem:interpolation-high} gives
  \begin{equation*}
    \begin{split}
      \lVert x_\lambda - x^{\dagger_s}\rVert_s \leq E\left[\psi^{-1}_s\left(\frac{\gamma_2}{m^2}\lambda\right)\right]^{(p-s)/(2s)} \left[\psi^{-1}_s\left(\frac{\gamma_2}{m^2}\lambda\right)\right]^{(s-\tau)/(2s)}=E\left[\psi^{-1}_s\left(\frac{\gamma_2}{m^2}\lambda\right)\right]^{(p-\tau)/(2s)}.
    \end{split}
  \end{equation*}
\end{proof}

\subsubsection{Misspecified setting}
Since in the misspecified setting $\mathcal{S}\cap X_s=\varnothing$, one needs $\Pi$ to extend to $X$ for (\ref{eq:plain-conv}) to still make sense.
This leads to the following assumption.
\begin{assumption}\label{ass:angle}
  The operator $\Pi:X_s\to X$ extends to a bounded operator $\overline{\Pi}:X \to X$. I.e.,
  \begin{equation}\label{eq:anglecond}
    \sup\Big\{\,\big\lVert B^{-s} P_{\N(A)^\perp} B^{s}x\big\rVert \ :\ x\in X_s,\ \norm{x}\le 1\,\Big\}<\infty.
\end{equation}
\end{assumption}
We give sufficient conditions for this assumption to hold.
\begin{proposition}\label{prop:sufficient}
  Each of the following implies~\ref{ass:angle}:
  \begin{enumerate}[label=(\roman*)]
  \item \textbf{Commutation.} $P_{\Nt}$ commutes with the spectral measure of $B$. In that case $\overline{\Pi}=P_{\Nt^\perp}$.
  \item \textbf{Smooth finite-dimensional null space.} $T$ satisfies $\dim\Nt<\infty$ and $\Nt\subset X_{2s}$.
  \item \textbf{Two-sided polynomial link condition.} $T$ and $B$ satisfy
  $$\hspace{1.75cm}\norm{T x} \ge m\,  \lVert B^{-a}x\rVert \qquad \text{for all } x\in B^{-s}(\Na^\perp)$$and
  $$\norm{T x} \le M\,  \lVert B^{-a} x\rVert \qquad \text{for all } x\in X$$
  for some $a>0$ and $0 < m\leq M$.
  \end{enumerate}
\end{proposition}
Notice that the case where $T$ is injective, is covered by (i) and (ii).
\begin{proof}
  \begin{enumerate}[label=(\roman*)]
  \item Let us first show that under the commutation assumption, we have $\Na=\Nt$.
  By definition of the spectral calculus, $P_{\Nt}$ must commute with $B^{-s}$ on $\D(B^{-s})=X$.
  Let $x\in \Nt$, then $$Ax=TB^{-s}x=TB^{-s}P_{\Nt}x=TP_{\Nt}B^{-s}x=0.$$ Hence, $\Nt \subseteq \Na$.
  Now let $x\in \Na$, then $B^{-s}x\in \Nt$ and therefore
  $$B^{-s}x=P_{\Nt}B^{-s}x\Longleftrightarrow B^{-s}x=B^{-s}P_{\Nt}x \Longleftrightarrow B^{-s}(x-P_{\Nt}x)=0,$$
  which implies that $x=P_{\Nt}x$ by injectivity of $B^{-s}$ and therefore $x\in \Nt$. This shows the other implication.
  Finally, for all $x\in \D(B^s)$ one obtains, again using the commutativity of $P_{\Nt}$ and $B^{-s}$, that $$(I-\overline{\Pi})x=B^{-s}P_{\Nt}B^sx=P_{\Nt}x.$$
  The final claim now follows from denseness of $\D(B^s)$ and boundedness of $P_{\Nt}$ on all of $X$.
  
  \item We start by showing that for any $e\in \Na$, we have $e \in \D(B^s)$.
  From the proof of Proposition \ref{prop:proj}, we obtain that $\N(A) = B^s(\Nt\cap X_s)=B^s(\Nt).$
  Therefore, for any $e\in \Na$, we have that $e=B^sv$ for some $v\in \Nt \subset X_{2s}$.
  From which follows that $\norm{B^se}=\norm{B^{2s}v}<+\infty$.

  Now since $\dim(\Na)=\dim(\Nt)=:d<\infty$, as $B^s$ is a bijection on $X_s$ (see proof of Proposition \ref{prop:proj}), we can find an orthonormal basis $\{e_i\}_{i=1}^{d}$ of $\Na$.
  Then for any $x\in X_s$, we have that
  $$\norm{B^{-s}\Pna B^{s}x} =\Big\lVert \sum_{i=1}^{d}\langle B^{s}x,e_i\rangle B^{-s}e_i\Big\rVert \leq c^{-s}\norm{x}\sum_{i=1}^{d}\norm{B^se_i},$$
  where $\sum_{i=1}^{d}\norm{B^se_i}$ is finite by the starting argument. This shows that Assumption \ref{ass:angle} holds.

  \item We start by showing that $$m^\nu \lVert x\rVert_{-\nu(a+s)}\leq \lVert (A^*A)^{\nu/2}x\rVert\leq M^{2\nu} \lVert x\rVert_{-\nu(a+s)}$$
  hold for all $x\in \N(A)^\perp$ and $0\leq \nu \leq 1$, where the upper bound holds on all $X$.
  Let $0\leq \nu \leq 1$. First the assumption gives $$m \norm{x}_{-(a+s)} \leq \norm{(A^*A)^{1/2}x} \leq M \norm{x}_{-(a+s)} $$
  for all $x\in \Na^\perp$ for the lower bound and for all $x\in X$ for the upper bound. 
  This gives $$m^2 P_{\Na^\perp}B^{-2(a+s)}P_{\Na^\perp} \leq P_{\Na^\perp}(A^*A)P_{\Na^\perp}=A^*A$$and $$A^*A\leq M^2 B^{-2(a+s)}.$$
  By Heinz' inequality we obtain for $0\leq \nu \leq 1$
  $$(A^*A)^\nu \leq M^{2\nu} B^{-2\nu(a+s)},$$ which gives the upper bound.
  For the lower bound, an application of Löwner-Heinz' inequality together with Hansen's inequality gives for $0\leq \nu \leq 1$
  $$m^{2\nu} P_{\Na^\perp}B^{-2\nu (a+s)} P_{\Na^\perp} \leq m^{2\nu} (P_{\Na^\perp}B^{-2(a+s)}P_{\Na^\perp})^\nu \leq (A^*A)^\nu,$$
  which in turn gives for $x\in \Na^\perp$
  $$m^{2\nu} \norm{x}_{-\nu(a+s)}^2=\langle P_{\Na^\perp}B^{-2\nu (a+s)} P_{\Na^\perp}x,x\rangle \leq \langle (A^*A)^\nu x,x\rangle =\norm{(A^*A)^{\nu/2}x}^2.$$
  Define $g_\lambda(t):=(t+\lambda)^{-1}$.
  Let $x\in \D(B^s)$ and set $\nu=s/(a+s)$. Then, using $A^*AB^sx\in \N(A)^\perp$ together with Lemma \ref{lem:invariance} and the two-sided inequality gives
  \begin{equation*}
    \begin{split}
      \norm{B^{-s}g_\lambda(A^*A)A^*Tx}=\norm{B^{-s}g_\lambda(A^*A)A^*AB^sx}&\leq m^{-\nu}\norm{(A^*A)^{\nu/2}(A^*A)g_\lambda(A^*A)B^sx}\\
      &\leq m^{-\nu}\norm{(A^*A)^{\nu/2}B^sx}\\
      &\leq \left(\frac{M}{m}\right)^\nu \norm{B^sx}_{-s}\\
      &=\left(\frac{M}{m}\right)^\nu \norm{x}.
    \end{split}
  \end{equation*}
  Since $X_s$ is dense in $X$, we obtain $$\sup_{0<\lambda }\norm{B^{-s}g_\lambda(A^*A)A^*T}<\infty.$$
  Now, since $B^{-s}g_\lambda(A^*A)A^*Tx \to \Pi x$ as $\lambda \to 0$ for all $x\in X_s$ by Proposition \ref{prop:proj}, we obtain that $\Pi$ must extend to a bounded operator on $X$ by Lemma 9.4.7 in \cite{eidelman2004functional}.
  \end{enumerate}
\end{proof}
The following can be said about $\overline{\Pi}$.
\begin{proposition}\label{prop:pi extension}
  Under Assumption \ref{ass:angle}, $\overline{\Pi}$ is a bounded idempotent operator on $X$ with
  $$\N(\overline{\Pi})=\overline{\Nt \cap X_s}, \qquad \Ra(\overline{\Pi})=\overline{(\Nt\cap X_s)^{\perp_s}} \qquad \text{and}\qquad T\overline{\Pi}=T.$$
\end{proposition}
\begin{proof}
  ($\overline{\Pi}$ is idempotent)
  For $x\in X$, let $(x_n)\subset X_s$ be a sequence such that $x_n \to x$ as $n\to \infty$. Then
  $$\overline{\Pi}^2x=\lim_{n\to \infty}\overline{\Pi}^2x_n=\lim_{n\to\infty}\Pi^2x_n=\lim_{n\to\infty}\Pi x_n=\overline{\Pi}x,$$
  which shows the idempotency of $\overline{\Pi}$.

  ($T=T\overline{\Pi}$) For $x\in X_s$, Proposition \ref{prop:proj} gives
  $$T(I-\Pi)x=TB^{-s}P_{\Na}B^sx=AP_{\Na}B^sx=0.$$
  Hence, $T=T\overline{\Pi}$ on $X_s$. Since, $T$ and $\overline{\Pi}$ are bounded, this extends on all of $X$ by continuity.

  ($\Ra(\overline{\Pi})$) By Proposition \ref{prop:proj} and continuity, we obtain $\Ra(\overline{\Pi})=\overline{\Ra(\Pi)}=\overline{(\Nt\cap X_s)^{\perp_s}}$.

  ($\N(\overline{\Pi})$) We first start by showing that $\N(\overline{\Pi})=\overline{\Nt \cap X_s}\subset \Nt$ in general.
  Let $x\in \overline{\Nt \cap X_s}$ and $(x_n)$ be a sequence in $\Nt \cap X_s$ such that $x_n\to x$ as $n\to \infty$.
  Then, by Proposition \ref{prop:proj} $$\overline{\Pi}x=\lim_{n\to \infty}\Pi x_n=0,$$ which shows $\overline{\Nt \cap X_s}\subseteq \N(\overline{\Pi})$.
  Conversely, let $x\in \N(\overline{\Pi})$. It suffices to construct a sequence $(z_n)\subset \Nt \cap X_s$ such that $z_n \to x$ as $n\to \infty$.
  Let $(x_n)\subset X_s$ be a sequence such that $x_n \to x$, then $((I-\overline{\Pi})x_n)\subset \N(\Pi)=\Nt \cap X_s$ since $\Pi(I-\overline{\Pi})x=(\Pi-\Pi^2)x=0$ for all $n$
  and $$z_n:=(I-\overline{\Pi})x_n \to (I-\overline{\Pi})x = x.$$
\end{proof}
This shows that $\overline{\Pi}(\mathcal{S})\subseteq \mathcal{S}$.
It is clear from the well-specified setting that one can only hope to converge to some element in $\overline{\Pi}(\mathcal{S})$.
However, in order for the limit to be well-defined, it needs to be unique.
In other words, $\overline{\Pi}(\mathcal{S})$ needs to contain only one element.
Fortunately for us, this is satisfied under general link conditions.
\begin{lemma}\label{lem: uniqueness}
  Assume there exists some index function $\phi$  and constant $m>0$ such that
  $\norm{T x} \ge m\,  \lVert [\phi(B^{-2})]^{1/2}x\rVert$ for all $x\in B^{-s}(\Na^\perp)$.
  Then, under Assumption \ref{ass:angle} we have that $\overline{\Pi}(\mathcal{S})$ contains a unique element, say $x^*$.
\end{lemma} 
\begin{proof}
  Since $\mathcal{S}=x+\N(T)$ for any $x\in \mathcal{S}$, it suffices to show that $\N(T)\subseteq \N(\overline{\Pi})$.
  Let $v\in \N(T)$ and $(x_n)_n\subset X_s$ a sequence such that $x_n \to v$ in $X$ as $n\to \infty$.
  Since $\Ra(\Pi)=B^{-s}(\N(A)^\perp)$ by Lemma \ref{lem:Naperp}, we obtain
  \begin{equation*}
    \begin{split}
      \lVert [\phi(B^{-2})]^{1/2}\overline{\Pi}v\rVert = \lim_{n\to \infty} \lVert [\phi(B^{-2})]^{1/2}\Pi x_n\rVert \leq m^{-1} \lim_{n\to \infty} \lVert T\Pi x_n\rVert = m^{-1} \lim_{n\to \infty} \lVert T x_n\rVert &= m^{-1}\lVert T v\rVert \\
      &=0
    \end{split}
  \end{equation*}
  Hence, $[\phi(B^{-2})]^{1/2}\overline{\Pi}v=0$ and since $[\phi(B^{-2})]^{1/2}$ is injective by injectivity of $B^{-2}$ and definition of $\phi$ we obtain that $\overline{\Pi}v=0$.
\end{proof}
We now introduce the link and source conditions needed in order to obtain convergence rates.
\begin{assumption}\label{ass:link-miss}
  There exists some index function $\phi$ with the following properties.
  \begin{enumerate}[label=(\roman*)]
  \item There exists a constant $0< M$ such that
  $$\norm{Tx} \le M\,\lVert [\phi(B^{-2})]^{1/2}x\rVert \qquad \text{for all } x\in X.$$
  \item There exists a constant $0<m\leq M$ such that
  $$\hspace{1.75cm}\norm{T x} \ge m\,\,  \lVert [\phi(B^{-2})]^{1/2}x\rVert \qquad \text{for all } x\in B^{-s}(\Na^\perp).$$
  \item The function $\psi_s:(0,c^{-2s}]\to \R^+:t\mapsto t\,\phi(t^{1/s})$ is convex.
  \end{enumerate}
\end{assumption}
\begin{assumption}\label{ass:source-miss}
  There exists constants $0<p < s$ and $R>0$ such that $\lVert B^{p}x^*\rVert \leq R$.
\end{assumption}
\begin{proposition}\label{prop: miss deterministic rate}
  Let $x_\lambda$ be the element defined in \eqref{eq:xlambda} with $s>0$, where $g_\lambda$ is a filter function satisfying Assumptions \ref{ass:filter-g} and \ref{ass:filter}.
  Then, under Assumptions \ref{ass:source-miss} and \ref{ass:link-miss}, we obtain the upper bound
  $$\lVert x_\lambda - x^*\rVert_\tau \leq (\rho_{0} M R \,c + \gamma_0 R\,c +R\norm{\overline{\Pi}}^{1-\tau/s})\,[\psi^{-1}_s(\lambda)]^{\frac{p-\tau}{2s}}\qquad \text{for }\tau \in [0,p],$$
  where $c =1\vee \frac{\rho_1}{m \rho_{0}} \vee \frac{\gamma_{1}}{m \gamma_0}$.
\end{proposition}

\begin{proof}
  Let $E$ be the spectral family of $B$ and set $v_N=\overline{\Pi}E([c,N])x^*$.
  Since $E([c,N])x^*\in X_s$, we have
  $v_N=B^{-s}P_{\mathcal N(A)^\perp}B^sE([c,N])x^*$, hence $B^sv_N\in\mathcal N(A)^\perp$ and
  $\|B^sv_N\|\le\|B^sE([c,N])x^*\|\le EN^{s-p}$. Moreover, $\overline\Pi x^*=x^*$ and
  $T\overline\Pi=T$ give
  $$v_N-x^*=\overline\Pi\,(E([c,N])-I)x^*,\qquad T(v_N-x^*)=-T(I-E([c,N]))x^*.$$
  We will use the following decomposition:
  $$\lVert x_\lambda -x^*\rVert_\tau \leq \underset{(1)}{\underbrace{\lVert B^{-s}g_\lambda(A^*A)A^*T(v_N-x^*)\rVert_\tau}}  +\underset{(2)}{\underbrace{\lVert B^{-s}g_\lambda(A^*A)A^*Tv_N-v_N\rVert_\tau}}+\underset{(3)}{\underbrace{\lVert v_N -x^*\rVert_\tau}}.$$
  The first term on the RHS is bounded by an application of Lemma \ref{lem:interpolation-high}:
  \begin{equation*}
    \begin{split}
      (1)&\leq \lVert g_\lambda(A^*A)A^*T(v_N-x^*)\rVert \left[ \psi^{-1}_s \left( \frac{\lVert Ag_\lambda(A^*A)A^*T(v_N-x^*)\rVert^2}{m^2\lVert g_\lambda(A^*A)A^*T(v_N-x^*)\rVert^2}\right)\right]^{\frac{s-\tau}{2s}}\\
      &\leq \rho_{0} M R \,\left(\psi_p(N^{-2p})\right)^{1/2} \lambda^{-1/2}\left[ \psi^{-1}_s \left( \frac{\rho_1^2}{m^2 \rho_{0}^2}\lambda\right)\right]^{\frac{s-\tau}{2s}},
    \end{split}
  \end{equation*}
  where we used the fact that 
  \begin{equation*}
    \begin{split}
      \lVert T(v_N-x^*)\rVert &\leq M \lVert [\phi(B^{-2})]^{1/2} (v_N-x^*)\rVert \\ &=M \left(\int_{(N,\infty)}\phi(\lambda^{-2})\,d\norm{E_\lambda x^*}^2\right)^{1/2} \\
      &=M \left(\int_{(N,\infty)}\psi_p(\lambda^{-2p}) \lambda^{2p}\,d\norm{E_\lambda x^*}^2\right)^{1/2}  \leq M\,R \, \left(\psi_p(N^{-2p})\right)^{1/2} .
    \end{split}
  \end{equation*}
  For the second term, again applying Lemma \ref{lem:interpolation-high} gives
  \begin{equation*}
    \begin{split}
      (2)&= \lVert B^{-s} r_\lambda(A^*A)B^s v_N\rVert_\tau \\ &\leq \lVert r_\lambda(A^*A)B^s v_N\rVert \left[ \psi^{-1}_s \left( \frac{\lVert Ar_\lambda(A^*A)B^s v_N\rVert^2}{m^2\lVert r_\lambda(A^*A)B^s v_N\rVert^2}\right)\right]^{\frac{s-\tau}{2s}}\\
      &\leq \gamma_0 R \,N^{s-p}\left[ \psi^{-1}_s \left( \frac{\gamma_{1}^2}{m^2 \gamma_0^2}\lambda\right)\right]^{\frac{s-\tau}{2s}},
    \end{split}
  \end{equation*}
  where we used $\lVert B^s v_N\rVert^2 \leq \int_{[c,N]}\lambda^{2s}d\norm{E_\lambda x^*}^2=\int_{[c,N]}\lambda^{2(s-p)}\lambda^{2p}d\norm{E_\lambda x^*}^2\leq N^{2(s-p)}\norm{x^*}_p^2$.
  For the remaining term we obtain by an application of Heinz-Kato's inequality (see e.g. Theorem 1 in \cite{roidos2019heinzkatoinequalitybanachspaces}) that
  $$(3) \leq \norm{\overline{\Pi}}^{1-\tau/s}\left(\int_{[N,\infty)}\lambda^{2\tau}d\norm{E_\lambda x^*}^2\right)^{1/2}=\norm{\overline{\Pi}}^{1-\tau/s}\left(\int_{[N,\infty)}\lambda^{2(\tau-p)}\lambda^{2p}d\norm{E_\lambda x^*}^2\right)^{1/2} \leq \norm{\overline{\Pi}}^{1-\tau/s} R N^{\tau-p}.$$
  Hence, setting $c :=1\vee \frac{\rho_1}{m \rho_{0}} \vee \frac{\gamma_{1}}{m \gamma_0}$ and using Lemma \ref{lem:nair 3.6}, we obtain
  $$\lVert x_\lambda -x^*\rVert_\tau \leq \rho_{0} M R \,c\,\left(\psi_p(N^{-2p})\right)^{1/2} \lambda^{-1/2}\left[\psi^{-1}_s \left(\lambda\right)\right]^{\frac{s-\tau}{2s}}+ \gamma_0 R\,c \,N^{s-p}\left[\psi^{-1}_s \left( \lambda\right)\right]^{\frac{s-\tau}{2s}} + \norm{\overline{\Pi}}^{1-\tau/s} R N^{\tau-p}.$$
  Now tuning the parameter $N$ by $N:=[\psi^{-1}_s(\lambda)]^{-1/(2s)}$ gives
  $$\left(\psi_p(N^{-2p})\right)^{1/2} \lambda^{-1/2}\left[\psi^{-1}_s \left(\lambda\right)\right]^{\frac{s-\tau}{2s}}=N^{s-p}\left[ \psi^{-1}_s \left( \lambda\right)\right]^{\frac{s-\tau}{2s}}=N^{\tau-p}=[\psi^{-1}_s(\lambda)]^{\frac{p-\tau}{2s}},$$
  and therefore$$\lVert x_\lambda -x^*\rVert_\tau \leq (\rho_{0} M R \,c + \gamma_0 R\,c +R\norm{\overline{\Pi}}^{1-\tau/s})\,[\psi^{-1}_s(\lambda)]^{\frac{p-\tau}{2s}}.$$
\end{proof}

\newpage
\section{Real and Vector-Valued Reproducing Kernel Hilbert Spaces}\label{app:RKHS}
Let $\mathbb{Y}  $ be a real Hilbert space, $\mathcal{X} $ a second countable locally compact topological space, $\mathbb{Y}^\mathcal{X} $ the space of functions from $\mathcal{X} $ to $\mathbb{Y} $ and $(\mathcal{X} ,\mathcal{B}_\mathcal{X} ,\mu)$ be a measure space. This section introduces the concepts surrounding vector-valued reproducing kernel Hilbert spaces (Section~\ref{subsec:RKHS}) and interpolation spaces (Section~\ref{subsec:interpolation}). It is self-sufficient for the purpose of this paper. We refer the reader to \cite{li2022optimal} for a more detailed description.

\subsection{Vector-valued reproducing kernel Hilbert spaces}
\label{subsec:RKHS}
\begin{definition}
    A Hilbert space $\mathcal{H} \subset \mathbb{Y} ^\mathcal{X} $ is a vector-valued reproducing kernel Hilbert space (RKHS) if the evaluation operator \mbox{$ev_x:\mathcal{H} \to \mathbb{Y} :f\mapsto f(x)$} is continuous for all $x\in \mathcal{X} $.
\end{definition}
In order to provide an explicit construction of such RKHSs, one requires positive semidefinite (psd) operator-valued kernels.
\begin{definition}
    A map $K:\mathcal{X} \times \mathcal{X} \to \mathcal{L} (\mathbb{Y} )$ is called a psd operator-valued kernel if $K(x,t)^*=K(t,x)$ and
    $$\sum_{i,j=1}^{N}\left\langle K(x_i,x_j)y_j,y_i\right\rangle _\mathbb{Y} \geq 0$$
    for any $x_1,\dots,x_N\in \mathcal{X} $, $y_1,\dots,y_N\in \mathbb{Y} $ and $N\geq1$.
\end{definition}
The following proposition discusses the construction of a vector-valued RKHS associated to an operator-valued kernel.
\begin{proposition}\label{vRKHS}
    Let $K:\mathcal{X} \times \mathcal{X} \to \mathcal{L} (\mathbb{Y} )$ be a psd operator-valued kernel. 
    Then, there exists a unique vector-valued RKHS, $\mathcal{H} _K$, satisfying, for all $f\in \mathcal{H} _K$ and $x\in \mathcal{X}$,
    $$K_x \in \mathcal{L} (\mathbb{Y} ,\mathcal{H} _K) \hspace{0.5cm}\text{ and }\hspace{0.5cm} f(x)=K_x^* f,$$
    where $K_x$ is defined by $(K_xy)(t):=K(t,x)y$ for all $y\in \mathbb{Y} $ and $t\in \mathcal{X} $. Furthermore, it holds $K(x,t)=K_x^*K_t$ and
    $$\mathcal{H} _K=\overline{\operatorname{span}}\{K_xy|x\in \mathcal{X} ,y\in\mathbb{Y} \}.$$
\end{proposition} 
Note that, when $\mathbb{Y} =\mathbb{R}$, one recovers real-valued RKHS as described by \cite{aronszajn1950theory}.
In that case, the operator-valued kernel is replaced by a symmetric real-valued kernel $k:\mathcal{X} \times \mathcal{X} \to \mathbb{R} $ such that $\sum_{i,j=1}^N a_i a_j k(x_i,x_j)\geq 0$ for all $\{a_i\}_{i=1}^N\subset \mathbb{R} $ and $\{x_i\}_{i=1}^N\subset \mathcal{X} $.
The kernel $k$ produces a unique RKHS 
$$\mathcal{H} _k=\overline{\operatorname{span}}\{k(\cdot,x)|x\in \mathcal{X}  \}\subset \mathbb{R}^\mathcal{X}$$ such that the reproducing property $f(x)=\left\langle f,k(\cdot,x)\right\rangle_{\mathcal{H} _k} $ holds for all $f\in \mathcal{H} _k$.
Let 
$$\Phi:\mathcal{X}\rightarrow \mathcal{H}_k:x\mapsto \Phi(x):=k(\cdot,x)$$ 
be the canonical feature map.

A popular choice of operator-valued kernels $K$ are those of the form $K(x,t)=k(x,t)\text{Id}_\mathbb{Y}$, where $k$ is a real-valued kernel and $\text{Id}_\mathbb{Y} $ is the identity operator on $\mathbb{Y} $.
In that case, one can show that $\mathcal{H} _K \cong \mathcal{H}_k \otimes \mathbb{Y} \cong S_2(\mathcal{H} _k,\mathbb{Y} )$, see \cite{carmeli2006vector}.
More precisely, this relationship is captured by the following theorem.
\begin{theorem}[\cite{li2022optimal}, Corollary 1]\label{HK vs Hk}
    For every function $F\in \mathcal{H} _K$ there exists a unique operator $C\in S_2(\mathcal{H} _k,\mathbb{Y} )$ such that $F(\cdot)=C\Phi (\cdot)$ with $\left\lVert C\right\rVert_{S_2(\mathcal{H} _k,\mathbb{Y} )} =\left\lVert F\right\rVert _{\mathcal{H} _K}$.
    Reciprocally, if $F(\cdot)=C\Phi(\cdot)$ for some $C\in S_2(\mathcal{H} _k,\mathbb{Y})$, then $F\in\mathcal{H}_K$ and the norm equality holds.
\end{theorem}
Therefore, Theorem~\ref{HK vs Hk} implies that $\mathcal{H} _K=\{F:\mathcal{X} \to \mathbb{Y} : F=C\Phi(\cdot), C\in S_2(\mathcal{H} _k,\mathbb{Y} )\}$. This representation is known as the kernel trick in the context of vector-valued RKHSs. In the sequel, $F$ and $C\Phi$ will be used interchangeably. 

A kernel $K$ (respectively, $k$) satisfies the \emph{universal property} if its RKHS is dense in $L_2(\mathcal{X} ,\mathcal{B} _\mathcal{X}  ,\mu,\mathbb{Y})$ (respectively, $L_2(\mathcal{X} ,\mathcal{B} _\mathcal{X}  ,\mu)$) for any Borel probability measure $\mu$.
Additionally, a kernel $K$ (resp., $k$) satisfies the $c_0$-universal property if its RKHS is dense in the space of continuous functions (on $\mathbb{Y}$, resp., on $\mathbb{R}$) that vanish at infinity, equipped with the supremum norm. 
Note that a $c_0$-universal kernel is universal. 
It is known (see Example 14 in \cite{carmeli2006vector}) that a kernel of the form $K=k\,\text{Id}_\mathbb{Y} $ is $c_0$-universal if and only if  $k$ is $c_0$-universal. The well-known Gaussian (i.e. $k(x,x')=\exp(-\gamma^2\|x-x'\|^2)$) and Laplace (i.e. $k(x,x')=\exp(-\lambda \|x-x'\|)$) kernels are $c_0$-universal, hence universal.

\subsection{Vector-valued interpolation spaces}
\label{subsec:interpolation}
Let $k$ be a bounded scalar-valued kernel with canonical feature map $\Phi$ and define the vector-valued kernel $K=k\,\text{Id}_\mathbb{Y} $.
Assume that $\mathcal{H} _k$ is separable and that $\Phi$ is measurable. 
This is a standard assumption in the literature. 
Let $\mathcal{I}_\mu:\mathcal{H}_k \to L_2(\mu):=L_2(\mathcal{X} ,\mathcal{B}_\mathcal{X} ,\mu)$ be the mapping that maps any $f\in \mathcal{H}_k$ to its $\mu$-equivalence class $[f]_{L_2(\mu)}$.
Under the assumption that $k$ is a bounded kernel, and noting that $\left\lVert \mathcal{I}_\mu \right\rVert_{S_2(\mathcal{H}_k,L_2(\mu))}=\left\lVert x\mapsto k(x,x)^{1/2}\right\rVert_{L_2(\mu)}$ (see \cite{steinwart2012mercer}, Lemma 2.2 and 2.3), the operator $\mathcal{I}_\mu$ is compact.

Define the self-adjoint, positive semi-definite and trace class integral operators (see \cite{steinwart2008support}, Theorem 4.27)$$\mathcal{L}_\mu := \mathcal{I}_\mu \mathcal{I}_\mu^*:L_2(\mu)\to L_2(\mu) \hspace{1cm}\text{and}\hspace{1cm}\mathcal{C}_\mu:=\mathcal{I}_\mu^* \mathcal{I}_\mu : \mathcal{H}_k \to \mathcal{H}_k.$$
One can easily show that $\mathcal{C}_\mu = \mathbb{E}_{\mu}[\Phi(X)\otimes \Phi(X)]$.
By the spectral theorem for self-adjoint compact operators, there exists a non-increasing positive sequence of real numbers $\{u_i\}_{i\in \mathbb{N} }$ and a subset $\{e_i\}_{i\in \mathbb{N} }\subseteq \mathcal{H}_k$ such that (i) $\{[e_i]_\mu\}_{i\in \mathbb{N}}$ is an orthonormal basis (ONB) of $\overline{\mathcal{R}(\mathcal{I}_\mu)}\subseteq L_2(\mu)$; (ii) $\{u_i^{1/2}e_i\}_{i\in \mathbb{N} }$ is an ONB of $\mathcal{N}(\mathcal{I}_\mu)^\perp \subseteq \mathcal{H}_k$, and (iii) the operators $\mathcal{L}_\mu$ and $\mathcal{C}_\mu$ admit the representations
$$\mathcal{L}_\mu = \sum_{i\in \mathbb{N} }u_i \left\langle \cdot, [e_i]_\mu\right\rangle_{L_2(\mu)} [e_i]_\mu \hspace{1cm}\text{and}\hspace{1cm}\mathcal{C}_\mu = \sum_{i\in \mathbb{N} }u_i\left\langle \cdot,u_i^{1/2}e_i\right\rangle_{\mathcal{H}_k} u_i^{1/2}e_i.$$

For $\gamma \geq 0$, define the $\gamma$-interpolation spaces 
$$[\mathcal{H}_k]^\gamma := \Big\{\sum_{i\in \mathbb{N}}a_i u_i^{\gamma/2}[e_i]_\mu:{a_i}\in l_2(\mathbb{N} )\Big\}\subseteq L_2(\mu),$$
equipped with the inner product 
$$\Big\langle \sum_{i\in \mathbb{N}}a_i u_i^{\gamma/2}[e_i]_\mu,\sum_{i\in \mathbb{N}}b_i u_i^{\gamma/2}[e_i]_\mu\Big\rangle_{[\mathcal{H}_k]^\gamma}=\sum_{i\in \mathbb{N}}a_i b_i,$$
for $\{a_i\}_{i\in \mathbb{N} },\{b_i\}_{i\in \mathbb{N} }\in l_2(\mathbb{N} )$.
The set $[\mathcal{H}_k]^\gamma$ is a separable Hilbert space with ONB given by $\{u_i^{\gamma/2}[e_i]_\mu\}$.

Let $\Xi$ be the isometric isomorphism between $S_2(L_2(\mu),\mathbb{Y} )$ and $L_2(\mu,\mathbb{Y})$ (which exists since $L_2(\mu,\mathbb{Y} )\cong L_2(\mu) \otimes \mathbb{Y} \cong S_2(L_2(\mu),\mathbb{Y} )$), and $\overline{\Xi}$ be the isometric isomorphism between $S_2(\mathcal{H}_k,\mathbb{Y} )$ and $\mathcal{H}_K$.
Let us also define the embedding operator $\mathbf{I}_\mu:S_2(\mathcal{H}_k,\mathbb{Y} )\to S_2(L_2(\mu),\mathbb{Y} )$ which is the unique extended operator that acts as
$$\mathbf{I}_\mu(f \otimes y) = \mathcal{I}_\mu(f) \otimes y$$ for any $f \otimes y \in S_2(\mathcal{H}_k,\mathbb{Y} ).$
One can easily see that $ \mathbf{I}_\mu = A_{\mathcal{I}^*_\mu}$.
Then, for $F\in \mathcal{H}_K$,  define $[F]_\mu := \Xi \circ \mathbf{I}_\mu \circ \overline{\Xi}^{-1}(F)$.
Let us now define the vector-valued interpolation spaces through
$$[\mathcal{H}_K]^\gamma := \Xi(S_2([\mathcal{H}_k]^\gamma,\mathbb{Y} ))\subseteq L_2(\mu,\mathbb{Y}).$$ 
One can show that $[\mathcal{H}_K]^\gamma$ is a Hilbert space when equipped with the scalar product
$$\left\langle F,G\right\rangle_{[\mathcal{H}_K]^\gamma}=\left\langle \Xi^{-1}(F),\Xi^{-1}(G)\right\rangle_{S_2([\mathcal{H}_k]^\gamma,\mathbb{Y} )}. $$
Notice that one has
$$\left\lVert F \right\rVert_{[\mathcal{H}_K]^\gamma}=\left\lVert \Xi^{-1}(F) \right\rVert_{S_2([\mathcal{H}_k]^\gamma,\mathbb{Y} )}=\left\lVert \Xi^{-1}(F) \mathcal{L}_\mu^{-\frac{\gamma}{2}}\right\rVert_{S_2(L_2(\mu),\mathbb{Y} )} .$$
Furthermore, for $\gamma=0$, we retrieve $\left\lVert F\right\rVert_{[\mathcal{H}_K]^0}=\left\lVert F\right\rVert_{L_2(\mu,\mathbb{Y} )}$.
If $F\in [\mathcal{H}_K]^1$, then there exists $\widetilde{F}\in \mathcal{H}_K$ such that $F=\widetilde{F}$ $\mu$-a.e. and $\left\lVert F\right\rVert_{[\mathcal{H}_K]^1}=\lVert \widetilde{F}\rVert_{\mathcal{H}_K}$.
\newpage

\section{$A_C$ operator}\label{app:AC}
Let $\mathbb{H} $ and $\mathbb{K} $ be separable real Hilbert spaces.
It is well-known that there exists an isometric isomorphism between the Hilbert spaces $S_2(\mathbb{H} ,\mathbb{K} )$ and $\mathbb{H} \otimes \mathbb{K} $, where $\mathbb{H} \otimes \mathbb{K}$ is the tensor Hilbert space as constructed in \cite{murray1936rings}.
Letting $\Xi:S_2(\mathbb{H},\mathbb{K})\to \mathbb{H} \otimes \mathbb{K}$ be that isometric isomorphism, one has that $A_C = \Xi^{-1}(C^* \otimes \Id_{\mathbb{K} })\Xi$, where $C:\mathcal{D}(C)\subseteq \mathbb{H} \to \mathbb{H}$ is a densely defined operator, as is shown in the next result.
We will use $\odot$ to denote the algebraic tensor.
\begin{proposition}\label{prop:AC-tensor}
  Let $C$ be densely defined. Then, $A_C$ is closed, and
  $$A_C=\Xi^{-1}\,(C^*\otimes\mathrm{Id}_{\mathbb K})\,\Xi,$$
  as an equality of unbounded operators, including domains:
  $\mathcal D(A_C)=\Xi^{-1}\big(\mathcal D(C^*\otimes\mathrm{Id}_\mathbb{K} )\big)$.
\end{proposition}

\begin{proof}
We use two elementary facts. First, for $A\in S_2(\mathbb H,\mathbb K)$ and an orthonormal basis
$(f_j)$ of $\mathbb K$,
\begin{equation}\label{eq:XiExpansion}
  \Xi A=\sum_{j\ge1}A^*f_j\otimes f_j
  \quad\text{in }\mathbb H\otimes\mathbb K,
  \qquad
  \|A\|_{S_2}^2=\sum_{j\ge1}\|A^*f_j\|_{\mathbb H}^2 .
\end{equation}
Indeed, if $P_N$ is the orthogonal projection onto $\mathrm{span}\{f_1,\dots,f_N\}$, then
$P_NA=\sum_{j\le N}\langle\cdot,A^*f_j\rangle_{\mathbb H}f_j$, so
$\Xi(P_NA)=\sum_{j\le N}A^*f_j\otimes f_j$, and
$\|A-P_NA\|_{S_2}^2=\|A^*(I-P_N)\|_{S_2}^2=\sum_{j>N}\|A^*f_j\|^2\to0$. Since
$\Xi$ is continuous, \eqref{eq:XiExpansion} follows. Second, for $a\in\mathcal D(C^*)$ and
$x\in\mathcal D(C)$ we have $\langle C^*a,x\rangle_{\mathbb H}=\langle a,Cx\rangle_{\mathbb H}$, by the
definition of the adjoint.

\emph{Step 0: $C^*\odot \Id_\mathbb{K} $ is closable.} Let $u_n\in\mathcal D(C^*)\odot\mathbb K$ with $u_n\to0$ and
$(C^*\odot\mathrm{Id}_{\mathbb K})u_n\to v$ in $\mathbb H\otimes\mathbb K$. For $x\in\mathcal D(C)$ and
$y\in\mathbb K$, we have that
$\langle(C^*\odot\mathrm{Id}_{\mathbb K})u_n,x\otimes y\rangle=\langle u_n,Cx\otimes y\rangle\to0$. Hence
$\langle v,x\otimes y\rangle=0$ for all $x\in\mathcal D(C)$ and $y\in\mathbb K$. These elementary tensors
span a dense subspace of $\mathbb H\otimes\mathbb K$, because $\mathcal D(C)$ is dense in $\mathbb H$.
Therefore $v=0$, i.e.\ $C^*\odot\mathrm{Id}_{\mathbb K}$ is closable.

\emph{Step 1: $A_C$ is closed.} Let $K_n\in\mathcal D(A_C)$ with $K_n\to K$ and $A_CK_n\to L$ in
the Hilbert--Schmidt norms. For $x\in\mathcal D(C)$, $(A_CK_n)x=K_nCx$, and convergence in
$S_2$ implies convergence in operator norm, so $K_nCx\to KCx$ and $(A_CK_n)x\to Lx$. Hence
$KCx=Lx$ for all $x\in\mathcal D(C)$, i.e.\ $KC$ extends to $L\in S_2(\mathbb H,\mathbb K)$. Thus
$K\in\mathcal D(A_C)$ and $A_CK=L$.

\emph{Step 2: the two operators agree on finite-rank elements.} Let
$u=\sum_{k=1}^ra_k\otimes y_k$ with $a_k\in\mathcal D(C^*)$ and $y_k\in\mathbb K$, and set
$K:=\Xi^{-1}u=\sum_k\langle\cdot,a_k\rangle_{\mathbb H}y_k$. For $x\in\mathcal D(C)$,
$$KCx=\sum_{k=1}^r\langle Cx,a_k\rangle_{\mathbb H}\,y_k=\sum_{k=1}^r\langle x,C^*a_k\rangle_{\mathbb H}\,y_k ,$$
which is the restriction to $\mathcal D(C)$ of the finite-rank operator
$L:=\sum_k\langle\cdot,C^*a_k\rangle_{\mathbb H}y_k\in S_2(\mathbb H,\mathbb K)$. Hence
$K\in\mathcal D(A_C)$ and
$A_CK=L=\Xi^{-1}\big(\sum_kC^*a_k\otimes y_k\big)
=\Xi^{-1}(C^*\odot\mathrm{Id}_{\mathbb K})u$.

\emph{Step 3: $\Xi^{-1}(C^*\otimes\mathrm{Id}_{\mathbb K})\Xi\subseteq A_C$.}
Let $w\in\mathcal D(C^*\otimes\mathrm{Id}_{\mathbb K})$. By definition of the closure, there are
$u_n\in\mathcal D(C^*)\odot\mathbb K$ with $u_n\to w$ and
$(C^*\odot\mathrm{Id}_{\mathbb K})u_n\to(C^*\otimes\mathrm{Id}_{\mathbb K})w$ (by step 0). By Step~2,
$K_n:=\Xi^{-1}u_n\in\mathcal D(A_C)$ and
$A_CK_n=\Xi^{-1}(C^*\odot\mathrm{Id}_{\mathbb K})u_n$. Since
$\Xi^{-1}$ is an isometry, $K_n\to\Xi^{-1}w$ and
$A_CK_n\to\Xi^{-1}(C^*\otimes\mathrm{Id}_{\mathbb K})w$. As $A_C$ is closed (Step~1),
$\Xi^{-1}w\in\mathcal D(A_C)$ and
$A_C\,\Xi^{-1}w=\Xi^{-1}(C^*\otimes\mathrm{Id}_{\mathbb K})w$.

\emph{Step 4: $A_C\subseteq\Xi^{-1}(C^*\otimes\mathrm{Id}_{\mathbb K})\Xi$.}
Let $K\in\mathcal D(A_C)$ and $L:=A_CK$. For every $j$ and $x\in\mathcal D(C)$,
$$\langle x,L^*f_j\rangle_{\mathbb H}=\langle Lx,f_j\rangle_{\mathbb K}=\langle KCx,f_j\rangle_{\mathbb K}
  =\langle Cx,K^*f_j\rangle_{\mathbb H} .$$
Thus $x\mapsto\langle Cx,K^*f_j\rangle_{\mathbb H}$ is bounded on $\mathcal D(C)$, so
$K^*f_j\in\mathcal D(C^*)$ and $C^*K^*f_j=L^*f_j$. Set
$u_N:=\sum_{j\le N}K^*f_j\otimes f_j\in\mathcal D(C^*)\odot\mathbb K$. By~\eqref{eq:XiExpansion}
applied to $K$ and to $L$,
$$ u_N\to\Xi K\qquad\text{and}\qquad
  (C^*\odot\mathrm{Id}_{\mathbb K})u_N=\sum_{j\le N}L^*f_j\otimes f_j\to\Xi L .$$
By definition of the closure, $\Xi K\in\mathcal D(C^*\otimes\mathrm{Id}_{\mathbb K})$ and
$(C^*\otimes\mathrm{Id}_{\mathbb K})\Xi K=\Xi L=\Xi A_CK$.

Steps~3 and~4 give the result.
\end{proof}
We therefore only discuss properties concerning the operator $C^*\otimes \Id_\mathbb{K} $ as these then easily extend to $A_C$ through the isometric isomorphism.
We will only discuss the properties that will be needed throughout the proofs.
We will assume throughout that $C$ is densely defined.

We start by discussing the domain of $C^* \otimes Id_\mathbb{K} $.
\begin{lemma}
  $\mathcal{D}(C^*\otimes \Id_\mathbb{K})$ is dense in $\mathbb{H}\otimes\mathbb{K}$ if and only if $\mathcal{D}(C^*)$ is dense in $\mathbb{H}$.
\end{lemma}
\begin{proof}
  Recall that $C^*$ is closed, and that $C^*\otimes \Id_\mathbb{K}$ is the closure of the operator
  $C^*\odot Id_\mathbb{K}$ defined on the algebraic tensor product $\mathcal{D}(C^*)\odot\mathbb{K}$.
  Hence
  $$\mathcal{D}(C^*)\odot\mathbb{K}\subseteq\mathcal{D}(C^*\otimes \Id_\mathbb{K}),$$
  and, by definition of the closure, every $u\in\mathcal{D}(C^*\otimes \Id_\mathbb{K})$ is the limit
  in $\mathbb{H}\otimes\mathbb{K}$ of a sequence in $\mathcal{D}(C^*)\odot\mathbb{K}$. Therefore both
  subspaces have the same closure in $\mathbb{H}\otimes\mathbb{K}$, and it suffices to show that
  $\mathcal{D}(C^*)\odot\mathbb{K}$ is dense in $\mathbb{H}\otimes\mathbb{K}$ if and only if
  $\mathcal{D}(C^*)$ is dense in $\mathbb{H}$
  
  ($\Leftarrow$) Assume $\mathcal{D}(C^*)$ is dense in $\mathbb{H}$, and let $M$ denote the closure of
  $\D(C^*)\odot\mathbb{K}$ in $\mathbb{H}\otimes\mathbb{K}$. Let $x\in\mathbb{H}$ and $y\in\mathbb{K}$.
  Choose $x_n\in \mathcal{D}(C^*)$ with $x_n\to x$. Then $x_n\otimes y\in \mathcal{D}(C^*)\odot\mathbb{K}$ and
  $$\lVert x_n\otimes y-x\otimes y\rVert_{\mathbb{H}\otimes\mathbb{K}}=\lVert (x_n-x)\otimes y\rVert_{\mathbb{H}\otimes\mathbb{K}}=\lVert x_n-x\rVert_{\mathbb{H}}\,\lVert y\rVert_{\mathbb{K}}\longrightarrow 0 .$$
  Hence every elementary tensor $x\otimes y$ belongs to $M$. Since $M$ is a closed linear subspace,
  it contains the algebraic tensor product $\mathbb{H}\odot\mathbb{K}$, which is dense in
  $\mathbb{H}\otimes\mathbb{K}$ by construction of the Hilbert tensor product.
  Thus $M=\mathbb{H}\otimes\mathbb{K}$.

  ($\Rightarrow$) Assume $\mathcal{D}(C^*)$ is not dense in $\mathbb{H}$. Since
  $\mathbb{H}=\overline{\mathcal{D}(C^*)}\oplus \mathcal{D}(C^*)^\perp$, there exists $x_0\neq 0\in \mathcal{D}(C^*)^\perp$.
  Fix $y_0\neq 0\in\mathbb{K}$. For all $x\in \mathcal{D}(C^*)$ and $y\in\mathbb{K}$,
  $$\langle x\otimes y,\,x_0\otimes y_0\rangle_{\mathbb{H}\otimes\mathbb{K}}=\langle x,x_0\rangle_{\mathbb{H}}\,\langle y,y_0\rangle_{\mathbb{K}}=0 ,$$
  so that $x_0\otimes y_0\perp \D(C^*)\odot\mathbb{K}$. Moreover
  $\lVert x_0\otimes y_0\rVert=\lVert x_0\rVert\,\lVert y_0\rVert\neq 0$. Hence
  $(\D(C^*)\odot\mathbb{K})^\perp\neq\{0\}$, and $\D(C^*)\odot\mathbb{K}$ is not dense in
  $\mathbb{H}\otimes\mathbb{K}$ which shows the contrapositive.
\end{proof}
We now characterise its adjoint.
\begin{lemma}
  One has that $(C^* \otimes \Id_\mathbb{K} )^*=\overline{C} \otimes Id_\mathbb{K} $. 
\end{lemma}
\begin{proof}
  See Proposition 7.26 in \cite{schmudgen2012unbounded}.
\end{proof}
The following can be said about its kernel and range.
\begin{lemma}\label{lem: ker tensor}
  The following properties hold
  \begin{enumerate}
    \item $\mathcal{N}(C^* \otimes \Id_\mathbb{K} )=\mathcal{N}(C^*)\otimes \mathbb{K} $,
    \item $\mathcal{R}(C^*\otimes \Id_\mathbb{K} )$ is dense in $\mathbb{H} \otimes \mathbb{K} $ iff $\mathcal{R}(C^*)$ is dense in $\mathbb{H} $.
  \end{enumerate}
\end{lemma} 
\begin{proof}
  See Lemma 3.3 in \cite{arai2018analysis}.
\end{proof}
We will also need to be able to apply the functional calculus on $C^* \otimes Id_\mathbb{K} $.
\begin{lemma}
  Let $f$ be a Borel measurable function and $C$ be self-adjoint, then $f(C \otimes \Id_\mathbb{K})=f(C)\otimes \Id_\mathbb{K}$. 
\end{lemma}
\begin{proof}
  This is a direct consequence of Theorems 3.3 and 3.5 in \cite{arai2018analysis}. 
\end{proof}

We add the following concernincg coercivity.
\begin{lemma}
  Assume $\langle C^*x,x\rangle_{\mathbb H}\ge\alpha\lVert x\rVert_{\mathbb H}^2$ for all $x\in\mathcal D(C^*)$.
  Then $\langle (C^*\otimes \Id_{\mathbb K})u,u\rangle_{\mathbb H\otimes\mathbb K}\ge\alpha\lVert u\rVert_{\mathbb H\otimes\mathbb K}^2$
  for all $u\in\mathcal D(C^*\otimes \Id_{\mathbb K})$.
\end{lemma}
\begin{proof}
  Let $u\in\mathcal D(C^*)\odot\mathbb K$, say $u=\sum_{i=1}^n x_i\otimes y_i$.
  Let $g_1,\dots,g_m$ be an orthonormal basis of $\operatorname{span}\{y_1,\dots,y_n\}$
  and write $y_i=\sum_{j}c_{ij}g_j$. Then $u=\sum_{j=1}^m \tilde x_j\otimes g_j$ with
  $\tilde x_j=\sum_i c_{ij}x_i\in\mathcal D(C^*)$. Since $\langle g_j,g_l\rangle=\delta_{jl}$,
  $$\langle (C^*\otimes \Id_{\mathbb K})u,u\rangle
    =\sum_{j,l}\langle C^*\tilde x_j,\tilde x_l\rangle_{\mathbb H}\langle g_j,g_l\rangle_{\mathbb K}
    =\sum_{j}\langle C^*\tilde x_j,\tilde x_j\rangle_{\mathbb H}
    \ge\alpha\sum_j\lVert \tilde x_j\rVert_{\mathbb H}^2=\alpha\lVert u\rVert^2 .$$
  Since $C^*\otimes \Id_{\mathbb K}$ is the closure of $C^*\odot \Id_{\mathbb K}$, so
  $\mathcal D(C^*)\odot\mathbb K$ is a core for $C^*\otimes \Id_{\mathbb K}$. 
  Given $u\in\mathcal D(C^*\otimes \Id_{\mathbb K})$,
  choose $u_n\in\mathcal D(C^*)\odot\mathbb K$ with $u_n\to u$ and
  $(C^*\otimes \Id_{\mathbb K})u_n\to(C^*\otimes \Id_{\mathbb K})u$.
  By continuity of the inner product,
  $\langle (C^*\otimes \Id_{\mathbb K})u_n,u_n\rangle\to\langle (C^*\otimes \Id_{\mathbb K})u,u\rangle$
  and $\lVert u_n\rVert^2\to\lVert u\rVert^2$, so the inequality passes to the limit.
\end{proof}
Now the following can be said if $C\in \mathcal{L}(\mathbb{H} )$.
\begin{lemma}\label{lem: properties bounded tensor op}
  Let $C_1,C_2\in \mathcal{L}(\mathbb{H} )$, $S$ be a subspace of $\mathbb{H} $ and $c>0$ some constant. Then,
  \begin{enumerate}
    \item $\lVert C_1^*\otimes \Id_\mathbb{K} \rVert = \lVert C_1^*\rVert $,
    \item $\overline{(C_1^*\otimes \Id_\mathbb{K})(S\odot \mathbb{K})}=\overline{C_1^*(S)}\otimes \mathbb{K} $,
    \item $\lVert (C^*_1\otimes \Id_\mathbb{K} )x\rVert_{\mathbb{H} \otimes \mathbb{K} } \leq c \lVert (C^*_2\otimes \Id_\mathbb{K} )x\rVert_{\mathbb{H} \otimes \mathbb{K} }$ for all $x\in \overline{S}\otimes \mathbb{K} $ if and only if
      $\lVert C^*_1 u \rVert_{\mathbb{H} } \leq c \lVert C^*_2 u\rVert_{\mathbb{H}}$ for all $u\in S$.
  \end{enumerate}
\end{lemma}
\begin{proof}
\textit{Proof of 1.}
See Theorem 3.2 in \cite{arai2018analysis}.

\textit{Proof of 2.}
  Write $D:=C_1^*(S)$.

  \emph{Step 1: $(C_1^*\otimes \Id_\mathbb K)(S\odot\mathbb K)=D\odot\mathbb K$.}
  If $x=\sum_{i=1}^n u_i\otimes y_i$ with $u_i\in S$, then
  $(C_1^*\otimes \Id_\mathbb K)x=\sum_{i=1}^n C_1^*u_i\otimes y_i\in D\odot\mathbb K$. Conversely, every
  element of $D\odot\mathbb K$ is of the form $\sum_{i=1}^n C_1^*u_i\otimes y_i$ with $u_i\in S$, which
  is the image of $\sum_{i=1}^n u_i\otimes y_i\in S\odot\mathbb K$.

  \emph{Step 2: $\overline{D\odot\mathbb K}=\overline D\otimes\mathbb K$.}
  Since $D\odot\mathbb K\subseteq\overline D\odot\mathbb K$, we have
  $\overline{D\odot\mathbb K}\subseteq\overline{\overline D\odot\mathbb K}=\overline D\otimes\mathbb K$.
  Conversely, let $z\in\overline D$ and $y\in\mathbb K$, and choose $z_n\in D$ with $z_n\to z$.
  Then $z_n\otimes y\in D\odot\mathbb K$ and
  $$\lVert z_n\otimes y-z\otimes y\rVert_{\mathbb H\otimes\mathbb K}=\lVert z_n-z\rVert_\mathbb H\,\lVert y\rVert_\mathbb K\longrightarrow 0 ,$$
  so $z\otimes y\in\overline{D\odot\mathbb K}$. Since $\overline{D\odot\mathbb K}$ is a closed
  linear subspace containing all elementary tensors $z\otimes y$ with $z\in\overline D$, it
  contains $\overline D\odot\mathbb K$ and hence its closure $\overline D\otimes\mathbb K$.

  Combining both steps, gives
  $\overline{(C_1^*\otimes \Id_\mathbb K)(S\odot\mathbb K)}=\overline{D\odot\mathbb K}
  =\overline{C_1^*(S)}\otimes\mathbb K$.

\textit{Proof of 3.}
  ($\Leftarrow$) Assume $\lVert C^*_1u\rVert_\mathbb H\le c\lVert C^*_2u\rVert_\mathbb H$ for all
  $u\in S$. Let first $x\in S\odot\mathbb K$ and write $x=\sum_{j=1}^m x_j\otimes g_j$ with
  $x_j\in S$ and wlog $\{g_j\}$ orthonormal.
  Then,
  $$\lVert (C^*_1\otimes \Id_\mathbb K)x\rVert^2=\sum_{j=1}^m\lVert C^*_1x_j\rVert^2\le c^2\sum_{j=1}^m\lVert C^*_2x_j\rVert^2=c^2\lVert (C^*_2\otimes Id_\mathbb K)x\rVert^2 .$$
  Now let $x\in\overline S\otimes\mathbb K$. Since $\overline S\otimes\mathbb K$ is the closure of
  $\overline S\odot\mathbb K$, and $\overline{S\odot\mathbb K}=\overline S\otimes\mathbb K$ by
  Step~2 of the proof of item~2, there exist $x_n\in S\odot\mathbb K$ with
  $x_n\to x$. Since $C^*_1\otimes \Id_\mathbb K$ and $C^*_2\otimes \Id_\mathbb K$ are bounded by item~1,
  $(C^*_i\otimes \Id_\mathbb K)x_n\to(C^*_i\otimes \Id_\mathbb K)x$ for $i=1,2$, and letting
  $n\to\infty$ in
  $\lVert (C^*_1\otimes \Id_\mathbb K)x_n\rVert\le c\lVert (C^*_2\otimes \Id_\mathbb K)x_n\rVert$
  gives the inequality for $x$.

  ($\Rightarrow$) Assume the inequality holds for all $x\in\overline S\otimes\mathbb K$. Fix
  $y\in\mathbb K$ with $\lVert y\rVert_\mathbb K=1$ and let $u\in S$. Then
  $u\otimes y\in S\odot\mathbb K\subseteq\overline S\otimes\mathbb K$, and for $i=1,2$,
  $$\lVert (C^*_i\otimes \Id_\mathbb K)(u\otimes y)\rVert_{\mathbb H\otimes\mathbb K}=\lVert C^*_iu\otimes y\rVert_{\mathbb H\otimes\mathbb K}=\lVert C^*_iu\rVert_\mathbb H\,\lVert y\rVert_\mathbb K=\lVert C^*_iu\rVert_\mathbb H .$$
  Hence $\lVert C^*_1u\rVert_\mathbb H\le c\lVert C^*_2u\rVert_\mathbb H$.
\end{proof}

We also have the following concerning projection operators.
\begin{lemma}\label{lem: tensor projection}
  Let $S\subseteq \mathbb{H} $ be a closed subspace. Then $P_{S\otimes \mathbb{K} }=P_S \otimes \Id_\mathbb{K} $.
\end{lemma}
\begin{proof}
  Set $Q:=P_S\otimes Id_\mathbb{K} \in\mathcal L(\mathbb H\otimes \mathbb{K} )$. 
  It is bounded with $\lVert Q\rVert=\lVert P_S\rVert\le 1$.

  \emph{$Q$ is idempotent.} For $x\in\mathbb H$ and $y\in \mathbb{K} $,
  $Q^2(x\otimes y)=P_S^2x\otimes y=P_Sx\otimes y=Q(x\otimes y)$. By linearity $Q^2=Q$ on
  $\mathbb H\odot \mathbb{K} $, and since both sides are bounded and $\mathbb H\odot \mathbb{K} $ is dense,
  $Q^2=Q$ on $\mathbb H\otimes \mathbb{K} $.

  \emph{$Q$ is self-adjoint.} For $x,x'\in\mathbb H$ and $y,y'\in \mathbb{K} $,
  $$\langle Q(x\otimes y),x'\otimes y'\rangle=\langle P_Sx,x'\rangle_\mathbb H\langle y,y'\rangle_\mathbb{K} =\langle x,P_Sx'\rangle_\mathbb H\langle y,y'\rangle_\mathbb{K} =\langle x\otimes y,Q(x'\otimes y')\rangle .$$
  By bilinearity $\langle Qu,v\rangle=\langle u,Qv\rangle$ for all $u,v\in\mathbb H\odot \mathbb{K} $,
  and by density and continuity for all $u,v\in\mathbb H\otimes \mathbb{K} $.

  Hence $Q$ is an orthogonal projection and it remains to show $\mathcal R(Q)=S\otimes \mathbb{K} $.

  \emph{$S\otimes \mathbb{K} \subseteq\mathcal R(Q)$.} For $x\in S$, $Q(x\otimes y)=x\otimes y$, so
  $S\odot \mathbb{K} \subseteq\mathcal R(Q)$. Since $\mathcal R(Q)$ is closed, it contains
  $\overline{S\odot \mathbb{K} }=S\otimes \mathbb{K} $.

  \emph{$\mathcal R(Q)\subseteq S\otimes \mathbb{K} $.} For $x\in\mathbb H$, $Q(x\otimes y)=P_Sx\otimes y\in S\odot \mathbb{K} $,
  so $Q(\mathbb H\odot \mathbb{K} )\subseteq S\odot \mathbb{K} $. Given $u\in\mathbb H\otimes \mathbb{K} $, choose
  $u_n\in\mathbb H\odot \mathbb{K} $ with $u_n\to u$; then $Qu=\lim_n Qu_n\in\overline{S\odot \mathbb{K} }=S\otimes \mathbb{K} $.
\end{proof}
\newpage

\section{Operator Lipschitz condition for Filters}\label{app: lipschitz}
Assume $A,B\geq 0$ are bounded selfadjoint operators and let $g_\lambda:\R^+\to \R$ be a bounded Borel measurable filter function.
We wish to determine well-known examples of $g_\lambda$ such that $$\left\lVert (A+\lambda \Id)^{1/2}[g_\lambda(A)-g_\lambda(B)](B+\lambda Id)\right\rVert \leq c_g \lambda^{-1/2}\left\lVert A-B\right\rVert  .$$
We show that the Landweber and Showalter filters, which are respectively defined by
\begin{equation*}
  \begin{split}
    g_{\eta,m}(t) &= \frac{1-(1-\eta t)^m}{t}\qquad \hspace{0.55cm}\text{for }m\in \mathbb{N} \text{ and }\eta \leq 1/(\left\lVert A\right\rVert \vee \left\lVert B\right\rVert  )\text{ with }\lambda=1/(\eta m),\\
    g_{\lambda}(t) &= \frac{1-e^{-t/\lambda}}{t},
  \end{split}
\end{equation*}
satisfy the operator Lipschitz condition.
We however do not claim that the obtained constants are optimal.
We first need the following lemmas.
\begin{lemma}\label{lem: bounding inequalities}
  Let $m\in \mathbb{N} $ and $l\in \{0,\dots,m\}$.
  Then, $$\sup_{0\leq x \leq m} \sqrt{1+x}\,(1-x/m)^l < \sqrt{\frac{2m}{l+1}}.$$
\end{lemma}
Recall the inequalities $1-x\leq e^{-x}$ and $1+x\leq e^x$ for all $x\in \mathbb{R^+}$.
\begin{proof}
  We obtain for $l\geq 1$ that
  \begin{equation*}
    \begin{split}
      \sqrt{1+x}\,(1-x/m)^l \leq \sqrt{1+x}\;e^{-\frac{lx}{m}} \leq \sqrt{\frac{m}{l}+x} \; e^{-\frac{lx}{m}}&\leq  \sqrt{\frac{m}{l}+x} \; e^{-\frac{lx}{2m}}\\
      &=\sqrt{m/l} \left((1+lx/m)e^{-lx/m}\right)^{1/2}\\
      &\leq \sqrt{m/l} \\
      &\leq \sqrt{\frac{2m}{l+1}}.
    \end{split}
  \end{equation*}
  for all $0\leq x \leq m$. One can easily check the case $l=0$.
\end{proof}
\begin{lemma}\label{lem: unif conv}
  The function $x\mapsto \frac{1-(1-x/m)^m}{x}$ converges uniformly to $x\mapsto (1-e^{-x})/x$ on any compact set $[0,M]$ for $M>0$.
\end{lemma}
\begin{proof}
This follows from the fact that $\frac{1-(1-x/m)^m}{x}=\int_{0}^{1}(1-sx/m)^{m-1}\,ds$, $(1-e^{-x})/x=\int_{0}^{1}e^{-sx}\,ds$ 
and the functions $s\mapsto(1-sx/m)^{m-1}$ converge uniformly to $s\mapsto e^{-sx}$ on $[0,1]$.
\end{proof}
\begin{proposition}
  The Landweber and Showalter filters satisfy the operator Lipschitz condition with $c_g=\sqrt{2}\frac{13}{3}$.
\end{proposition}
Recall the telescoping equality $A^k-B^k=\sum_{i=0}^{k-1}A^i(A-B)B^{k-1-i}$ for $k\in \mathbb{N}_0$ which can easily be shown by induction.
\begin{proof}
We start by rewriting the Landweber filter as $$g_{\eta,m}(t)=\eta \sum_{i=0}^{m-1}(1-\eta t)^{i},$$
with $\lambda=1/(\eta m ).$
Set $H_A:=(A+\lambda)g_{\eta,m}(A)$ and $H_B:=(B+\lambda)g_{\eta,m}(B)$, then
\begin{equation*}
  \begin{split}
    &\lVert (A+\lambda \Id)^{1/2}[g_{\eta,m}(A)-g_{\eta,m}(B)](B+\lambda \Id)\rVert \\ 
    &= \lVert (A+\lambda \Id)^{1/2}[(A+\lambda)^{-1}H_A-(B+\lambda)^{-1}H_B](B+\lambda)\rVert \\
    &= \lVert (A+\lambda \Id)^{-1/2}[H_A(B+\lambda)-(A+\lambda)H_B]\rVert\\
    &=\lVert (A+\lambda \Id)^{-1/2}[H_A(B-A)+(A+\lambda)H_A -(A+\lambda)H_B]\rVert\\
    &\leq \lVert (A+\lambda \Id)^{-1/2}H_A(B-A)\rVert + \lVert (A+\lambda \Id)^{1/2}[H_A -H_B]\rVert\\
    &\leq \lVert (A+\lambda \Id)^{1/2}g_{\eta,m}(A)\rVert \,\lVert B-A\rVert +\lVert (A+\lambda \Id)^{1/2}[H_A-H_B]\rVert \\
    &\leq \sqrt2\,\lambda^{-1/2}\|B-A\|+\|(A+\lambda \Id)^{1/2}[H_A-H_B]\|,
  \end{split}
\end{equation*}
where we used $0\le g_{\eta,m}(t)\le\min\{\eta m,1/t\}=\min\{\lambda^{-1},t^{-1}\}$, so that
$(t+\lambda)^{1/2}g_{\eta,m}(t)\le(2\lambda)^{1/2}\lambda^{-1}$ for $t\le\lambda$ and
$(t+\lambda)^{1/2}g_{\eta,m}(t)\le(2t)^{1/2}t^{-1}\le\sqrt2\,\lambda^{-1/2}$ for $t\ge\lambda$.

Now notice that
$$(t+\lambda)g_{\eta,m}(t)=1-(1-\eta t)^m+\frac{1}{m}\sum_{i=0}^{m-1}(1-\eta t)^i.$$
Hence, using the telescoping equality, we obtain
\begin{equation*}
  \begin{split}
    H_A-H_B 
    &=(1-\eta A)^m-(1-\eta B)^m+\frac{1}{m}\left[\sum_{i=0}^{m-1}(1-\eta A)^i-(1-\eta B)^i\right]\\
    &=\eta \sum_{i=0}^{m-1}(1-\eta A)^i(B-A)(1-\eta B)^{m-1-i}+\frac{\eta}{m}\left[\sum_{i=0}^{m-1}\sum_{j=0}^{i-1}(1-\eta A)^j(B-A)(1-\eta B)^{i-1-j}\right].
  \end{split}
\end{equation*}
Since, $\eta \leq (\left\lVert A\right\rVert \vee \left\lVert B\right\rVert )^{-1}$, we have that $\left\lVert (\Id-\eta B)^k\right\rVert \leq 1$ for all $k\in \mathbb{N}$ and by Lemma \ref{lem: bounding inequalities} we have
$$\left\lVert (A+\lambda)^{1/2}(Id - \eta A)^k\right\rVert \leq \sqrt{\lambda}\sqrt{\frac{2m}{k+1}}$$ for all $0\leq k \leq m.$
Hence, using $\sum_{k=0}^{l-1}(1+k)^{-1/2}\leq 2\sqrt{l}$ for $l\geq 1$, we obtain
\begin{equation*}
  \begin{split}
    \left\lVert (A+\lambda \Id)^{1/2}[H_A -H_B]\right\rVert &\leq \eta\sqrt{\lambda}\,\lVert A-B\rVert\left[\sum_{j=0}^{m-1}\sqrt{\frac{2m}{j+1}}+\frac{1}{m}\sum_{i=1}^{m-1}\sum_{j=0}^{i-1}\sqrt{\frac{2m}{j+1}}\right]\\
    &\leq \eta\sqrt{\lambda}\,\lVert A-B\rVert \sqrt{2m}\left[2\sqrt{m}+\frac{1}{m}\sum_{i=1}^{m-1}2\sqrt{i}\right]\\
    &\leq \eta\sqrt{\lambda}\,\lVert A-B\rVert \sqrt{2m}\left[2\sqrt{m}+\frac{4}{3}\sqrt{m}\right]\\
    &=(c_g-\sqrt{2})\lambda^{-1/2}\lVert A-B\rVert.
  \end{split}
\end{equation*}

Setting $\eta=1/(m\lambda)$ in the Landweber filter, we obtain that the Showalter filter is the uniform limit of the Landweber filter on $[0,\lVert A\rVert \vee \lVert B\rVert]$ as $m\to \infty$ by Lemma \ref{lem: unif conv}.
This gives $$g_{\eta,m}(A)\underset{m\to \infty}{\to}g_\lambda(A)\qquad \text{and}\qquad g_{\eta,m}(B)\underset{m\to \infty}{\to}g_\lambda(B)$$
in the operator norm.
This finally gives,
\begin{equation*}
  \begin{split}
    \left\lVert (A+\lambda \Id)^{1/2}[g_\lambda(A)-g_\lambda(B)](B+\lambda \Id)\right\rVert
    &= \lim_{m\to \infty}\left\lVert (A+\lambda \Id)^{1/2}[g_{\eta,m}(A)-g_{\eta,m}(B)](B+\lambda \Id)\right\rVert \\
    & \leq c_g \lambda^{-1/2}\left\lVert A-B\right\rVert.
  \end{split}
\end{equation*}
\end{proof}

\newpage

\section{Practical implementation}\label{app: practical}
Recall that the estimator of Section~3 is defined by
$$\widehat M(z)=\widehat{\mathcal C}_{ZX}\,g_\varepsilon(\widehat{\mathcal C}_Z)\,\Phi_Z(z),
  \qquad
  \widehat F_\lambda(x)=\widehat{\mathcal C}_{\widehat MY}\,g_\lambda(\widehat{\mathcal C}_{\widehat M})\,\Phi_X(x),$$
with $u\otimes v=\langle\cdot,u\rangle v$ and
$$\widehat{\mathcal C}_{ZX}=\frac1m\sum_{j=1}^m\Phi_Z(\overline Z_j)\otimes\Phi_X(\overline X_j),\qquad
  \widehat{\mathcal C}_Z=\frac1m\sum_{j=1}^m\Phi_Z(\overline Z_j)\otimes\Phi_Z(\overline Z_j),$$
$$\widehat{\mathcal C}_{\widehat MY}=\frac1n\sum_{i=1}^n\widehat M(\widetilde Z_i)\otimes\widetilde Y_i,\qquad
  \widehat{\mathcal C}_{\widehat M}=\frac1n\sum_{i=1}^n\widehat M(\widetilde Z_i)\otimes\widehat M(\widetilde Z_i).$$ 
  We use the Gram matrices and vectors
$$\mathbf K_{\overline X\overline X}:=\big(k_X(\overline X_j,\overline X_l)\big)_{j,l\le m},\qquad
  \mathbf K_{\overline Z\overline Z}:=\big(k_Z(\overline Z_j,\overline Z_l)\big)_{j,l\le m},\qquad
  \mathbf K_{\overline Z\widetilde Z}:=\big(k_Z(\overline Z_j,\widetilde Z_i)\big)_{j\le m,\,i\le n},$$
$$\mathbf k_{\overline X}(x):=\big(k_X(\overline X_j,x)\big)_{j\le m}\in\mathbb R^m,\qquad
  \mathbf k_{\overline Z}(z):=\big(k_Z(\overline Z_j,z)\big)_{j\le m}\in\mathbb R^m .$$

\begin{proposition}\label{prop:closed-form}
Let
$$\mathbf B:=\frac1m\,g_\varepsilon\Big(\frac{\mathbf K_{\overline Z\overline Z}}{m}\Big)\,
  \mathbf K_{\overline Z\widetilde Z}\in\mathbb R^{m\times n},
  \qquad
  \mathbf G:=\mathbf B^\top\mathbf K_{\overline X\overline X}\mathbf B\in\mathbb R^{n\times n}.$$
Then, for every $x\in E_X$,
$$\widehat F_\lambda(x)=\sum_{i=1}^nw_i(x)\,\widetilde Y_i,
  \qquad
  \mathbf w(x):=\frac1n\,g_\lambda\Big(\frac{\mathbf G}{n}\Big)\,\mathbf B^\top\,
  \mathbf k_{\overline X}(x)\in\mathbb R^n .$$
\end{proposition}

\begin{proof}
\emph{Stage~1.} Let $S_{\overline Z}:\mathcal H_{k_Z}\to\mathbb R^m$, $h\mapsto(h(\overline Z_j))_{j\le m}$,
and $S_{\overline X}:\mathcal H_{k_X}\to\mathbb R^m$, $f\mapsto(f(\overline X_j))_{j\le m}$. By the
reproducing property, $S_{\overline Z}^*a=\sum_ja_j\Phi_Z(\overline Z_j)$ and
$S_{\overline X}^*a=\sum_ja_j\Phi_X(\overline X_j)$, so that
$$\widehat{\mathcal C}_Z=\tfrac1mS_{\overline Z}^*S_{\overline Z},\qquad
  \widehat{\mathcal C}_{ZX}=\tfrac1mS_{\overline X}^*S_{\overline Z},\qquad
  S_{\overline Z}S_{\overline Z}^*=\mathbf K_{\overline Z\overline Z},\qquad
  S_{\overline X}S_{\overline X}^*=\mathbf K_{\overline X\overline X},$$
and $S_{\overline Z}\Phi_Z(z)=\mathbf k_{\overline Z}(z)$. Hence,
$$\widehat M(z)=\tfrac1mS_{\overline X}^*S_{\overline Z}\,g_\varepsilon\big(\tfrac1mS_{\overline Z}^*S_{\overline Z}\big)\Phi_Z(z)
  =\tfrac1mS_{\overline X}^*\,g_\varepsilon\big(\tfrac1m\mathbf K_{\overline Z\overline Z}\big)\mathbf k_{\overline Z}(z)
  =S_{\overline X}^*\,\boldsymbol\beta(z),$$
with $\boldsymbol\beta(z):=m^{-1}g_\varepsilon(\mathbf K_{\overline Z\overline Z}/m)\mathbf k_{\overline Z}(z)\in\mathbb R^m$.
The $i$-th column of $\mathbf B$ is $\boldsymbol\beta(\widetilde Z_i)$, so
$\widehat M(\widetilde Z_i)=S_{\overline X}^*\mathbf B\mathbf e_i$, where $\mathbf e_i$ is the $i$-th
canonical basis vector of $\mathbb R^n$.

\emph{Stage~2.} Let $S_{\widehat M}:\mathcal H_{k_X}\to\mathbb R^n$,
$f\mapsto(\langle f,\widehat M(\widetilde Z_i)\rangle_{\mathcal H_{k_X}})_{i\le n}$, and
$\mathbf Y:\mathbb R^n\to\mathbb Y$, $\mathbf a\mapsto\sum_ia_i\widetilde Y_i$. Then
$S_{\widehat M}^*\mathbf a=\sum_ia_i\widehat M(\widetilde Z_i)=S_{\overline X}^*\mathbf B\mathbf a$, and therefore
$$\widehat{\mathcal C}_{\widehat M}=\tfrac1nS_{\widehat M}^*S_{\widehat M},\qquad
  \widehat{\mathcal C}_{\widehat MY}=\tfrac1n\mathbf Y S_{\widehat M},\qquad
  S_{\widehat M}S_{\widehat M}^*=\mathbf B^\top S_{\overline X}S_{\overline X}^*\mathbf B
  =\mathbf B^\top\mathbf K_{\overline X\overline X}\mathbf B=\mathbf G .$$
Furthermore,
$S_{\widehat M}\Phi_X(x)=\big(\widehat M(\widetilde Z_i)(x)\big)_{i\le n}
=\mathbf B^\top S_{\overline X}\Phi_X(x)=\mathbf B^\top\mathbf k_{\overline X}(x)$. Hence,
$$\widehat F_\lambda(x)=\tfrac1n\mathbf YS_{\widehat M}\,g_\lambda\big(\tfrac1nS_{\widehat M}^*S_{\widehat M}\big)\Phi_X(x)
  =\tfrac1n\mathbf Y\,g_\lambda\big(\tfrac1n\mathbf G\big)S_{\widehat M}\Phi_X(x)
  =\mathbf Y\,\mathbf w(x)=\sum_{i=1}^nw_i(x)\widetilde Y_i .$$
\end{proof}

\newpage

\section{Proofs}\label{app: proofs}
\subsection{Proof of Proposition \ref{prop:problem is hilbert scaled}}
\begin{proof}
  Using Lemmas \ref{lem: partial isometry} and \ref{lem: unitary equivalence L_2}, we obtain
      $$\mathbf{I_X}K_\lambda = A_{\mathcal{L}_X^{1/2}}Wg_\lambda(W^*Q^*QW)W^*Q^*r
      =A_{\mathcal{L}_X^{1/2}}g_\lambda(Q^*Q)Q^*r = C_\lambda$$
\end{proof}

\subsection{Proof of Proposition \ref{prop:proj_core}}
Follows from Proposition \ref{prop:proj} with $B=A_{\mathcal{L}_X^{-1}}$, $T=A_{T^*}$ and $s=1/2$.

\subsection{Proof of Proposition \ref{prop: X_s minimal core}}
Follows from Lemma \ref{lem: X_s minimal} with $B=A_{\mathcal{L}_X^{-1}}$, $T=A_{T^*}$ and $s=1/2$.

\subsection{Proof of Proposition \ref{prop:sufficient core}}

\begin{enumerate}
  \item Follows from Proposition \ref{prop:sufficient} with $B=A_{\mathcal{L}_X^{-1}}$, $T=A_{T^*}$ and $s=1/2$, together with Lemma \ref{lem: tensor projection}.
  \item Define $\Pi_{sc}:=\mathcal{L}_X^{1/2}P_{\mathcal{N}(T\mathcal{L}_X^{1/2})^\perp}\mathcal{L}_X^{-1/2}$, then $\Pi \cong \Pi_{sc} \otimes \Id_\mathbb{Y} $.
  Hence, the boundedness of $\Pi$ on the algebraic tensor product $L_2(X)\odot \mathbb{Y} $ follows from the boundedness of $\Pi_{sc}$ on $L_2(X)$ it suffices to apply Proposition \ref{prop:sufficient}.2 to $\Pi_{sc}$.
  \item Follows from Proposition \ref{prop:sufficient} with $B=A_{\mathcal{L}_X^{-1}}$, $T=A_{T^*}$ and $s=1/2$, together with Lemma \ref{lem: properties bounded tensor op}.3.
\end{enumerate}

\subsection{Proof of Proposition \ref{prop:pi extension core}}
Follows from Proposition \ref{prop:pi extension} with $B=A_{\mathcal{L}_X^{-1}}$, $T=A_{T^*}$ and $s=1/2$.

\subsection{Proof of Proposition \ref{prop: uniqueness}}
Follows from Lemma \ref{lem: uniqueness} with $B=A_{\mathcal{L}_X^{-1}}$, $T=A_{T^*}$ and $s=1/2$ together with Lemma \ref{lem: properties bounded tensor op}.3.

\subsection{Proof of Theorem \ref{th:main rate}}
We mainly follow the steps of Appendix E.2. in \cite{meunier2024nonparametric}.
We will write $\lVert F\rVert_\gamma $ and $\lVert C\rVert_\gamma$ to denote respectively $\lVert F\rVert_{[\mathcal{H}_{K_X}]^\gamma}$ and $\lVert C\rVert_{S_2([\mathcal{H}_{k_X}]^\gamma,\mathbb{Y} )}$.
We will use the following decomposition and study each term separately:
\begin{equation}
  \lVert [\widehat{F}_\lambda]_X-F^\dagger\rVert_\gamma \leq \underset{\text{1st stage error}}{\underbrace{\lVert [\widehat{F}_\lambda]_X-[\overline{F}_\lambda]_X\rVert_\gamma}}+\underset{\text{approximation error}}{\underbrace{\lVert [\overline{F}_\lambda]_X-[F_\lambda]_X\rVert_\gamma}}+
  \underset{\text{regularisation error}}{\underbrace{\lVert F^\dagger-[F_\lambda]_X\rVert_\gamma}}
\end{equation}
where $\overline{F}_\lambda:=\overline{K}_\lambda \Phi_X(\cdot):=\widehat{\mathcal{C}}_{MY}g_\lambda(\widehat{\mathcal{C}}_M)\Phi_X(\cdot)$.
We will write $\mathcal{D}_1$ for $\{\overline{X}_i,\overline{Z}_i\}_{i=1}^m$ and $\left\lVert \cdot \right\rVert_{S_2}$ for $\left\lVert \cdot \right\rVert_{S_2(\mathcal{H}_{k_X},\mathbb{Y} )}$.
Furthermore, define
$$S:=\{u \in   E_X^m  \times E_Z^m : \text{the set of inequalities in \eqref{eq: set of equations} are satisfied}\},$$
and let $s=(x_1,\dots,x_m,z_1,\dots,z_m)\in S$.
Furthermore, let $\Omega_0(s)$ be the measurable set described in Lemma \ref{lem:25} with $\mathbb{P}\left((\wZ_1,\dots,\wZ_n)^{-1}(\Omega_0(s)^c)\right)\leq 4e^{-\tau}$.

\subsubsection{Regularisation error}
Under the setting of Theorem \ref{th:main rate},we have by an application of Propositions \ref{prop:problem is hilbert scaled}, \ref{prop: well deterministic rate} and \ref{prop: miss deterministic rate},
where the latter two are used with $p=\beta_X/2,\tau=\gamma/2$ and $s=1/2$, that
$$ \lVert F^\dagger-[F_\lambda]_X\rVert_\gamma = \lVert C^\dagger-C_\lambda\rVert_\gamma \leq c \left[\psi^{-1}_{1/2}\left(\lambda\right)\right]^{\frac{\beta_X-\gamma}{2}},$$
for some constant $c>0$.

\subsubsection{Approximation error}\label{sec:approx}
We will assume for this section and the next one (stage 1 error) that we are working on $\Omega_0(s)$.
Using Lemma \ref{lem: inequality core} together with Lemma \ref{lem: membership}, we obtain that
$$\lVert [F_\lambda-\overline{F}_\lambda]\rVert_\gamma = \lVert [K_\lambda-\overline{K}_\lambda]\rVert_\gamma \leq \lambda^{-1/2}\,\lVert (K_\lambda-\overline{K}_\lambda)\mathcal{C}_{M,\lambda}^{1/2}\rVert_{S_2} \left[\psi^{-1}_{1/2}\left(\frac{\lambda}{c_1^2}\right)\right]^{\frac{1-\gamma}{2}}.$$
Furthermore, using $\Id=\widehat{\C}_M g_\lambda(\widehat{\C}_M)+r_\lambda(\widehat{\C}_M)$ and Lemma \ref{lem:25}, we have that
\begin{equation*}
  \begin{split}
    &\left\lVert (K_\lambda-\overline{K}_\lambda)\mathcal{C}_{M,\lambda}^{1/2}\right\rVert_{S_2}\\
    &\leq \sqrt{3}\left\lVert \big(K_\lambda\widehat{\C}_M g_\lambda(\widehat{\C}_M)+K_\lambda r_\lambda(\widehat{\C}_M)-\overline{K}_\lambda\big)\widehat{\mathcal{C}}_{M,\lambda}^{1/2}\right\rVert_{S_2} \\
    &\leq \sqrt{3}\left\lVert K_\lambda r_\lambda(\widehat{\C}_M)\widehat{\mathcal{C}}_{M,\lambda}^{1/2}\right\rVert_{S_2}+\sqrt{3}\left\lVert \big(K_\lambda\widehat{\C}_M -\widehat{\C}_{MY}\big)g_\lambda(\widehat{\C}_M)\widehat{\mathcal{C}}_{M,\lambda}^{1/2}\right\rVert_{S_2} \\
    &\leq \sqrt{3}\left\lVert K_\lambda r_\lambda(\widehat{\C}_M)\widehat{\mathcal{C}}_{M,\lambda}^{1/2}\right\rVert_{S_2}+6E\left\lVert \big(K_\lambda\widehat{\C}_M -\widehat{\C}_{MY}\big)\mathcal{C}_{M,\lambda}^{-1/2}\right\rVert_{S_2}\\
    &\leq \sqrt{3}\left\lVert K_\lambda r_\lambda(\widehat{\C}_M)\widehat{\mathcal{C}}_{M,\lambda}^{1/2}\right\rVert_{S_2}+6E\left\lVert \big\{\big(K_\lambda\widehat{\C}_M -\widehat{\C}_{MY}\big)-\big(K_\lambda \C_M -\C_{MY}\big)\big\}\mathcal{C}_{M,\lambda}^{-1/2}\right\rVert_{S_2}\\
    &\hspace{1cm}+6E\left\lVert \big(K_\lambda \C_M -\C_{MY}\big)\mathcal{C}_{M,\lambda}^{-1/2}\right\rVert_{S_2}\\
    &=: \sqrt{3}(I)+6E(II)+6E(III).
  \end{split}
\end{equation*}
We now bound each of the three terms separately.
Using Lemma \ref{lem: E[eta]}, we obtain
$$(III)= \left\lVert \C_{MY}r_\lambda( \C_M)\mathcal{C}_{M,\lambda}^{-1/2}\right\rVert_{S_2}\leq c\, \lambda^{1/2} \left[\psi_{1/2}^{-1}\left(\lambda\right)\right]^{(\beta_X-1)/2},$$
and by Lemma \ref{lem: E[eta] 2}, we obtain the same bound for $(I)$.
For the second term, set $$U_i:= (\C_{M,\lambda}^{-1/2}M(\widetilde{Z}_i))\otimes (\widetilde{Y}_i-K_\lambda M(\widetilde{Z}_i)).$$
Then, $$(II)=\left\lVert n^{-1}\sum_{i=1}^{n}U_i - \mathbb{E}[U]\right\rVert_{S_2}, $$where
\begin{equation*}
  \begin{split}
    &\mathbb{E}[\lVert U \rVert^k ]\\ &=\mathbb{E}\left[\lVert Y-K_\lambda M(Z)\rVert^k \big\lVert \C_{M,\lambda}^{-1/2}M(Z) \big\rVert^k \right]\\
    &\leq 2^{k-1} \mathbb{E}\left[\lVert Y-\mathbb{E} [F^\dagger(X)|Z]\rVert^k \big\lVert \C_{M,\lambda}^{-1/2}M(Z) \big\rVert^k \right]+2^{k-1} \mathbb{E}\left[\lVert K_\lambda M(Z)-\mathbb{E} [F^\dagger(X)|Z]\rVert^k \big\lVert \C_{M,\lambda}^{-1/2}M(Z) \big\rVert^k \right]\\
    &\leq \frac{2^{k-1}}{2}k!\sigma^2 L^{k-2}\mathbb{E} \left[\big\lVert \C_{M,\lambda}^{-1/2}M(Z) \big\rVert^k \right]+2^{k-1}\mathbb{E} \left[\lVert F^\dagger(X)-K_\lambda \Phi_X(X)\rVert^k \,\,\big\lVert \C_{M,\lambda}^{-1/2}M(Z) \big\rVert^k \right]\\
    &\leq \frac{2^{k-1}}{2}k!\sigma^2 L^{k-2}\mathbb{E} \left[\big\lVert \C_{M,\lambda}^{-1/2}M(Z) \big\rVert^k \right]+2^{k-1} A_X^k \left\lVert F^\dagger-[F_\lambda]\right\rVert_{\alpha_X}^k \mathbb{E} \left[\big\lVert \C_{M,\lambda}^{-1/2}M(Z) \big\rVert^k \right]\\
  \end{split}
\end{equation*}
with
\begin{equation*}
  \begin{split}
    \mathbb{E} \left[\big\lVert \C_{M,\lambda}^{-1/2}M(Z) \big\rVert^k \right]&\leq (\kappa_X \lambda^{-1/2})^{k-2}\mathbb{E} \left[\big\lVert \C_{M,\lambda}^{-1/2}M(Z) \big\rVert^2 \right]\\
    &=(\kappa_X \lambda^{-1/2})^{k-2} \, tr\left(\C_{M,\lambda}^{-1/2}\C_M \C_{M,\lambda}^{-1/2}\right)\\
    &=(\kappa_X \lambda^{-1/2})^{k-2} \mathcal{N}_M(\lambda).
  \end{split}
\end{equation*}
Hence,
\begin{equation*}
  \begin{split}
    \mathbb{E}[\lVert U \rVert^k ]&\leq \frac{2^{k-1}}{2}k!\sigma^2 L^{k-2}\mathbb{E} \left[\big\lVert \C_{M,\lambda}^{-1/2}M(Z) \big\rVert^k \right]+2^{k-1} k! A_X^k \left\lVert F^\dagger-[F_\lambda]\right\rVert_{\alpha_X}^k \mathbb{E} \left[\big\lVert \C_{M,\lambda}^{-1/2}M(Z) \big\rVert^k \right]\\
    &\leq \frac{k!}{2}\left((\sqrt{2}\sigma +2\sqrt{2}A_X \left\lVert F^\dagger-[F_\lambda]\right\rVert_{\alpha_X})\sqrt{\mathcal{N}_M(\lambda)} \right)^2 \left(\frac{2L\kappa_X+2A_X \kappa_X\left\lVert F^\dagger-[F_\lambda]\right\rVert_{\alpha_X} }{\sqrt{\lambda}}\right)^{k-2}
  \end{split}
\end{equation*}
Since, $\alpha_X < \beta_X$ and therefore $\lVert F^\dagger-[F_\lambda]\rVert_{\alpha_X}\to 0 $ as $\lambda \to 0$ by the regularisation error, we may assume for ease of writing
$$\max\left(2\sigma +\sqrt{2A_X} \left\lVert F^\dagger-[F_\lambda]\right\rVert_{\alpha_X},2L\kappa_X+2\kappa_X\left\lVert F^\dagger-[F_\lambda]\right\rVert_{\alpha_X}\right) \leq c$$
for some constant $c>0$.
Bernstein (see Theorem \ref{thm:bernstein}) now gives us, with $\mathbb{P}_{YZ}^{\otimes n}$-probability not less than $1-2e^{-\tau}$, that
$$(II)\leq \sqrt{32}c\,\frac{\tau}{\sqrt{n}}\left(\sqrt{\mathcal{N}_M(\lambda)}+\frac{\lambda^{-1/2}}{\sqrt{n}}\right)\leq \sqrt{32}c\,\frac{\tau}{\sqrt{n}}\left(\left[\psi^{-1}_{1/2}\left(\lambda\right)\right]^{-p_X/2}+\frac{\lambda^{-1/2}}{\sqrt{n}}\right),$$
where we used Lemma \ref{lem: effective dimension}.
Hence, denoting by $\Omega_1\in \mathbb{Y}^n \times E_Z^n$ the set on which this event holds, we have that
$$\lVert [F_\lambda-\overline{F}_\lambda]\rVert_\gamma \leq c \left(\left[\psi^{-1}_{1/2}\left(\frac{\lambda}{c_1^2}\right)\right]^{\frac{\beta_X-\gamma}{2}}+\tau \;\frac{\left[\psi^{-1}_{1/2}\left(\lambda\right)\right]^{\frac{1-\gamma-p_X}{2}}}{\sqrt{\lambda n}}+\tau \;\frac{\left[\psi^{-1}_{1/2}\left(\lambda\right)\right]^{\frac{1-\gamma}{2}}}{\lambda n}\right), $$
for some constant $c>0$ on $(\wZ_1,\dots,\wZ_n)^{-1}(\Omega_0(s)) \cap (\wY_1,\dots,\wY_n,\wZ_1,\dots,\wZ_n)^{-1}(\Omega_1)$.

\newpage
\subsubsection{Stage 1 error}\label{sec:stage 1}
We start by noticing that if $T$ is injective, then $A_{\mathcal{L}_X^{1/2}}(\N(Q)^\perp)$ is dense in $S_2(L_2(X),\mathbb{Y})$.
Indeed, one then obtains that
$$\N(A_{\mathcal{L}_X^{1/2}T^*})=\N(A_{\mathcal{L}_X^{1/2}})=\N(A_{\mathcal{L}_X})=0$$
by injectivity of $\mathcal{L}_X$.
And therefore $$\overline{A_{\mathcal{L}_X^{1/2}}(\N(A_{\mathcal{L}_X^{1/2}T^*})^\perp)}=\overline{\mathcal{R}(A_{\mathcal{L}_X^{1/2}})}=\overline{\mathcal{R}(A_{\mathcal{L}_X})}=S_2(L_2(X),\mathbb{Y} ).$$
Hence, if $T$ is injective then the first inequality displayed in Assumption \ref{ass: link core} holds by continuity of the norm on all of $S_2(L_2(X),\mathbb{Y})$
and one can then use Lemma \ref{lem: inequality core} to obtain 
$$\lVert [\overline{F}_\lambda-\widehat{F}_\lambda]\rVert_\gamma = \lVert [\overline{K}_\lambda-\widehat{K}_\lambda]\rVert_\gamma \leq \lambda^{-1/2}\,\lVert (\overline{K}_\lambda-\widehat{K}_\lambda)\mathcal{C}_{M,\lambda}^{1/2}\rVert_{S_2} \left[\psi^{-1}_{1/2}\left(\frac{\lambda}{c_1^2}\right)\right]^{\frac{1-\gamma}{2}}.$$
However, if $T$ is not injective, then $[\widehat{K}_\lambda]$ is not guarenteed to be in $A_{\mathcal{L}_X^{1/2}}(\N(Q)^\perp)$ or its closure.
In that case, the best we can do is
$$\lVert [\overline{F}_\lambda-\widehat{F}_\lambda]\rVert_\gamma = \lVert (\overline{K}_\lambda-\widehat{K}_\lambda)\mathcal{C}_X^{\frac{1-\gamma}{2}}\rVert_{S_2} \leq \kappa_X^{1-\gamma}\lambda^{-1/2}\,\lVert (\overline{K}_\lambda-\widehat{K}_\lambda)\mathcal{C}_{M,\lambda}^{1/2}\rVert_{S_2}.$$
Hence, setting $a:=\mathbb{I} [\N(T)=\{0\}]$, we obtain $$\lVert [\overline{F}_\lambda-\widehat{F}_\lambda]\rVert_\gamma \leq c \,\lambda^{-1/2} \,\lVert (\overline{K}_\lambda-\widehat{K}_\lambda)\mathcal{C}_{M,\lambda}^{1/2}\rVert_{S_2}\,\left[\psi^{-1}_{1/2}\left(\frac{\lambda}{c_1^2}\right)\right]^{a\frac{1-\gamma}{2}}.$$
It now remains to bound the term $\lVert (\overline{K}_\lambda-\widehat{K}_\lambda)\mathcal{C}_{M,\lambda}^{1/2}\rVert_{S_2}$.
We decompose this term the following way:
\begin{equation*}
  \begin{split}
    \lVert (\overline{K}_\lambda-\widehat{K}_\lambda)\mathcal{C}_{M,\lambda}^{1/2}\rVert_{S_2} &\leq \lVert \widehat{C}_{MY}(g_\lambda(\widehat{C}_{M})-g_\lambda(\widehat{C}_{\hat{M}}))\mathcal{C}_{M,\lambda}^{1/2}\rVert_{S_2} + \lVert( \widehat{C}_{MY}-\widehat{C}_{\hat{M}Y})g_\lambda(\widehat{C}_{\hat{M}})\mathcal{C}_{M,\lambda}^{1/2}\rVert_{S_2} =:S_0+S_{-1}.\\
  \end{split}
\end{equation*}
For the second term on the RHS, we obtain
\begin{equation*}
  \begin{split}
    S_{-1}&\leq \lVert( \widehat{C}_{MY}-\widehat{C}_{\hat{M}Y})\widehat{\mathcal{C}}_{M,\lambda}^{-1/2}\rVert_{S_2} \,\lVert \widehat{\mathcal{C}}_{M,\lambda}^{1/2}\widehat{\mathcal{C}}_{\hat{M},\lambda}^{-1/2}\rVert \,\lVert \widehat{\mathcal{C}}_{\hat{M},\lambda}g_\lambda(\widehat{\C}_{\hat{M}})\rVert \, \lVert \widehat{\mathcal{C}}_{\hat{M},\lambda}^{-1/2}\mathcal{C}_{M,\lambda}^{1/2}\rVert \\
    &\leq 2E\,\lVert( \widehat{C}_{MY}-\widehat{C}_{\hat{M}Y})\widehat{\mathcal{C}}_{M,\lambda}^{-1/2}\rVert_{S_2} \,\lVert \widehat{\mathcal{C}}_{M,\lambda}^{1/2}\widehat{\mathcal{C}}_{\hat{M},\lambda}^{-1/2}\rVert \, \lVert \widehat{\mathcal{C}}_{\hat{M},\lambda}^{-1/2}\mathcal{C}_{M,\lambda}^{1/2}\rVert \\
    &\leq 8E\,\lambda^{-1/2} \lVert\widehat{C}_{MY}-\widehat{C}_{\hat{M}Y}\rVert_{S_2}
  \end{split}
\end{equation*}
by Lemma \ref{lem:25}.
Furthermore, setting $W^s_i:=(\widehat{M}_s(\wZ_i)-M(\wZ_i))\otimes (\wY_i-\mathbb{E}[F^\dagger(\wX_i)|\wZ_i])$, where $\widehat{M}_s$ is assumed to be observed with observation $s\in S$, we obtain
\begin{equation*}
  \begin{split}
    \lVert\widehat{C}_{MY}-\widehat{C}_{\hat{M}_sY}\rVert_{S_2} \leq \Big\lVert n^{-1}\sum_{i=1}^{n}W^s_i \Big\rVert_{S_2} + \Big\lVert n^{-1}\sum_{i=1}^{n}(\widehat{M}_s(\wZ_i)-M(\wZ_i))\otimes \mathbb{E}[F^\dagger(\wX_i)|\wZ_i] \Big\rVert_{S_2}.
  \end{split}
\end{equation*}
Since
$$\mathbb{E}[W^s_i]=\mathbb{E}[\mathbb{E} [W^s_i|\wZ_i]]=\mathbb{E} \left[(\widehat{M}_s(\wZ_i)-M(\wZ_i))\otimes (\mathbb{E} [\wY_i|\wZ_i]-\mathbb{E}[F^\dagger(\wX_i)|\wZ_i])\right] =0$$and
\begin{equation*}
  \begin{split}
    \mathbb{E}[\lVert W^s_i\rVert_{S_2}^m ]&=\int_{E_Z}\lVert \widehat{M}_s(z)-M(z)\rVert_{\mathcal{H}_{k_X}}^m \mathbb{E}\left[ \big\lVert \wY_i-\mathbb{E}[F^\dagger(\wX_i)|\wZ_i=z]\big\rVert_\mathbb{Y} ^m\big. \big|\wZ_i=z \right] \mathbb{P}_Z(dz)\\
    &\leq \frac{1}{2}m!(LA_Z \lVert \widehat{M}_s-M\rVert_{\alpha_Z})^{m-2} (\sigma A_Z\lVert \widehat{M}_s-M\rVert_{\alpha_Z})^2 
  \end{split}
\end{equation*}
we have by Bernstein's theorem that there exists a set $\Omega_2(s)$ with $\mathbb{P}_{YZ}^n(\Omega_2(s))\geq 1-2e^{\tau}$ such that
$$\Big\lVert n^{-1}\sum_{i=1}^{n}W^s_i \Big\rVert_{S_2} \leq \sqrt{32}A_Z\frac{\tau}{\sqrt{n}}\left(\sigma\lVert \widehat{M}_s-M\rVert_{\alpha_Z}+\frac{L\lVert \widehat{M}_s-M\rVert_{\alpha_Z}}{\sqrt{n}}\right),$$
on that set.
Furthermore, setting $U_i:=\lVert (\widehat{M}_s(\wZ_i)-M(\wZ_i))\otimes \mathbb{E}[F^\dagger(\wX_i)|\wZ_i]\rVert_{S_2} $, we obtain another set $\Omega_3(s)$ with $\mathbb{P}_{Z}^n(\Omega_3(s))\geq 1-2e^{\tau}$, such that
\begin{equation*}
  \begin{split}
     &\Big\lVert n^{-1}\sum_{i=1}^{n}(\widehat{M}_s(\wZ_i)-M(\wZ_i))\otimes \mathbb{E}[F^\dagger(\wX_i)|\wZ_i] \Big\rVert_{S_2} \\&\leq n^{-1}\sum_{i=1}^{n}(U_i-\mathbb{E}[U])+\mathbb{E} [U]\\
     &\leq A_X \lVert F^\dagger\rVert_{\beta_X} \left(A_Z \kappa_X^{\beta_X-\alpha_X}\lVert M-\widehat{M}_s\rVert_{\alpha_Z}\sqrt{\frac{\tau}{n}}+\lVert M-\widehat{M}_s\rVert_{L_2(Z,\mathcal{H}_{k_X})}\right)
  \end{split}
\end{equation*}
on that set by Hoeffding's theorem, where we used the fact that almost surely
$$ U_i \leq \lVert \widehat{M}_s-M\rVert_\infty \big\lVert F^\dagger\big\rVert_\infty \leq A_Z A_X \lVert F^\dagger\rVert_{\alpha_X} \lVert M-\widehat{M}_s\rVert_{\alpha_Z}\leq \kappa_X^{\beta_X-\alpha_X}A_Z A_X \lVert F^\dagger\rVert_{\beta_X} \lVert M-\widehat{M}_s\rVert_{\alpha_Z} $$and
$$\mathbb{E} [U]\leq A_X \lVert F^\dagger\rVert_{\beta_X}\lVert M-\widehat{M}_s\rVert_{L_2(Z,\mathcal{H}_{k_X})}  $$
Hence, we have on
$$(\wZ_1,\dots,\wZ_n)^{-1}(\Omega_0(s)) \cap (\wY_1,\dots,\wY_n,\wZ_1,\dots,\wZ_n)^{-1}(\Omega_2(s))\cap (\wZ_1,\dots,\wZ_n)^{-1}(\Omega_3(s)),$$that
$$\lVert S_{-1}\rVert_{S_2} \leq J\lambda^{-1/2} \tau\left(\frac{\lVert M-\widehat{M}_s\rVert_{\alpha_Z}}{\sqrt{n}} +\lVert M-\widehat{M}_s\rVert_{L_2(Z,\mathcal{H}_{k_X})} \right)$$
for some constant $J>0$ and $\tau \geq 1$.

Let us now turn to the term $S_0$ for which we have on
$$(\wZ_1,\dots,\wZ_n)^{-1}(\Omega_0(s)) \cap (\wY_1,\dots,\wY_n,\wZ_1,\dots,\wZ_n)^{-1}(\Omega_1)$$that
\begin{equation*}
  \begin{split}
    \lVert S_0\rVert_{S_2} &\leq \lVert \widehat{\C}_{MY}\widehat{\C}_{M,\lambda}^{-1}\rVert_{S_2} \lVert \widehat{\C}_{M,\lambda}(g_\lambda(\widehat{\C}_M)-g_\lambda(\widehat{\C}_{\hat{M}}))\widehat{\C}_{\hat{M},\lambda}^{1/2}\rVert\, \lVert \widehat{\C}_{\hat{M},\lambda}^{-1/2}\C_{M,\lambda}^{1/2}\rVert \\
    &\leq 3c_g \lambda^{-1/2}\lVert \widehat{\C}_M-\widehat{\C}_{\hat{M}}\rVert \left(\left\lVert \widehat{\C}_{MY}\widehat{\C}_{M,\lambda}^{-1}-\C_{MY}\C_{M,\lambda}^{-1} \right\rVert_{S_2} +\left\lVert \widehat{\C}_{MY}\widehat{\C}_{M,\lambda}^{-1} \right\rVert_{S_2}\right)\\
    &\lesssim  3c_g \lambda^{-1/2}\lVert \widehat{\C}_M-\widehat{\C}_{\hat{M}}\rVert \left(\left\lVert \mathbf{I_X}(\widehat{\C}_{MY}\widehat{\C}_{M,\lambda}^{-1}-\C_{MY}\C_{M,\lambda}^{-1}) \right\rVert_{1} +\left[\psi_{1/2}^{-1}(\lambda)\right]^{((\beta_X-1)_-/2)}\right)\\
    &\leq c \tau \lambda^{-1/2}\lVert \widehat{\C}_M-\widehat{\C}_{\hat{M}}\rVert \left([\psi_{1/2}^{-1}(\lambda)]^{(\beta_X-1)_-/2}+\frac{[\psi_{1/2}^{-1}(\lambda)]^{-p_X/2}}{\sqrt{\lambda n}}+\frac{1}{\lambda n} \right),
  \end{split}
\end{equation*}
where we used Lemma \ref{lem: C lambda control}, \ref{lem:25} and the rate for the Approximation term with the Tikhonov filter and $\gamma=1$.
Using again Lemma \ref{lem:25} and putting the bounds together gives
\begin{equation*}
  \begin{split}
    &\lVert [\overline{F}_\lambda-\widehat{F}_\lambda]\rVert_\gamma \leq c \tau \,\lambda^{-1} \,\left[\psi^{-1}_{1/2}\left(\lambda\right)\right]^{a\frac{1-\gamma}{2}}\left([\psi_{1/2}^{-1}(\lambda)]^{(\beta_X-1)_-/2}+\frac{[\psi_{1/2}^{-1}(\lambda)]^{-p_X/2}}{\sqrt{\lambda n}}+\frac{1}{\lambda n} \right) \times \\ 
    &\hspace{5cm}\left(\frac{\lVert M-\widehat{M}\rVert_{\alpha_Z}}{\sqrt{n}} +\lVert M-\widehat{M}\rVert_{L_2(Z,\mathcal{H}_{k_X})} \right)
  \end{split}
\end{equation*}

for some constant $c>0$ on $$\Omega_4(s):=(\wZ_1,\dots,\wZ_n)^{-1}(\Omega_0(s)) \cap (\wY_1,\dots,\wY_n,\wZ_1,\dots,\wZ_n)^{-1}(\Omega_1\cap \Omega_2(s))\cap (\wZ_1,\dots,\wZ_n)^{-1}(\Omega_3(s)).$$
Notice that (given $s$) this event includes the event on which the bound for the approximation error holds.

\subsubsection{Probability of the bound}
Let $B\subseteq (E_X\times E_Z)^m \times (\mathbb{Y}\times E_Z )^n$ denote the event that the upperbound in Theorem \ref{th:main rate} is satisfied. 
For $d_1\in (E_X\times E_Z)^m$, define the sets $B_{d_1}:=\{d_2\in (\mathbb{Y} \times E_Z)^n:(d_1,d_2)\in B\}$.
Notice that sections \ref{sec:approx} and \ref{sec:stage 1} show that $(\wY_1,\dots,\wY_n,\wZ_1,\dots,\wZ_n)(\Omega_4(s))\subseteq B_s$ for any $s\in S$ with $\mathbb{P}(\Omega_4(s)^c)\leq 10 e^{-\tau} $.
We then obtain, by an application of Tonelli's theorem, that
\begin{equation*}
  \begin{split}
    \mathbb{P}(B^c)&=\mathbb{P} \left(B^c\cap (\oX_1,\dots,\oX_m,\oZ_1,\dots,\oZ_m)^{-1}(S)\right) + \mathbb{P} \left(B^c\cap (\oX_1,\dots,\oX_m,\oZ_1,\dots,\oZ_m)^{-1}(S)^c\right)\\
    &\leq \int_{S}\mathbb{P}_{YZ}^{\otimes n}(B_{d_1})\mathbb{P}_{XZ}^{\otimes m}(dd_1)+\mathbb{P}_{XZ}^{\otimes m}(S^c)\\
    &\leq \int_{S}\mathbb{P}_{YZ}^{\otimes n}(\Omega_4(d_1))\mathbb{P}_{XZ}^{\otimes m}(dd_1)+\mathbb{P}_{XZ}^{\otimes m}(S^c)\\
    &\leq 10 e^{-\tau}+\mathbb{P}_{XZ}^{\otimes m}(S^c).
  \end{split}
\end{equation*}

\newpage
\subsection{Proof of Proposition \ref{prop:rates}}
\begin{proof}
  Notice that by Theorem \ref{th: meunier optimal} there exists a set with $\mathbb{P}$-probability not less than $1-6e^{-\tau}$ on which
    $$\lVert M-\widehat{M}\rVert_{L^2(Z,\mathcal{H}_{K_X})}^2\leq J\tau^2 m^{-\beta_Z/(\beta_Z+p_Z)}\qquad \text{and}\qquad \lVert M-\widehat{M}\rVert_{\alpha_Z}^2\leq J' \tau^2 m^{(\alpha_Z-\beta_Z)/(\beta_Z+p_Z)}.$$
    We will show in Part (i) below that under the assumptions of the proposition, the set of inequalities in \eqref{eq: set of equations} will eventually be satisfied and therefore
    $\mathbb{P}_{XZ}^{\otimes m}(S^c)\leq 6e^{-\tau}$ for large enough $m$.

  For $n$ large enough, let $s^*_n$ be the unique $s\in (0,\lVert \mathcal{L_X}\rVert ]$ such that $n\psi_{1/2}(s)=\ln(1/s)$.

  We start by showing that for $n$ large enough $$\left(1-\frac{\ln \Lambda_n}{\Lambda_n}\right)\psi_{1/2}^{-1}(\Lambda_n/n)\leq s_n^* ,$$
  where $\Lambda_n:=\ln (1/\psi_{1/2}^{-1}(1/n))$. Since $s_n^*\to 0$, we have for large enough $n$ that $s_n^*\leq 1/e$ so that $\ln(1/s_n^*)\geq 1$.
  Then, $\psi_{1/2}(s_n^*)\geq 1/n$ and therefore $n\psi_{1/2}(s^*_n)=\ln(1/s_n^*)\leq \Lambda_n$.
  Using equation Lemma \ref{lem:nair 3.6}, together with upper bound, gives
  $$\ln(1/s_n^*)\geq \ln(1/\psi_{1/2}^{-1}(\Lambda_n/n))\geq \ln\left(\frac{1}{\Lambda_n\psi_{1/2}^{-1}(1/n)}\right)=\Lambda_n-\ln \Lambda_n.$$
  Therefore,
  $$n\psi_{1/2}(s_n^*)\geq \Lambda_n-\ln \Lambda_n \Longleftrightarrow s_n^*\geq\psi_{1/2}^{-1}\left(\frac{\Lambda_n-\ln \Lambda_n}{n}\right)\implies s^*_n \geq \left(1-\frac{\ln \Lambda_n}{\Lambda_n}\right)\psi_{1/2}^{-1}\left(\frac{\Lambda_n}{n}\right),$$
  for $n$ large enough.

  \textbf{(Part (i))}
  For the first part, since the condition $\ln(1/\psi_{1/2}^{-1}(\lambda_n))/\lambda_n = o(n)$ implies $n\lambda_n \geq 1$ eventually, it suffices to show (\ref{eq: set of equations}) is satisfied.
  The rate is then an application of Theorem \ref{th:main rate}.
  By Lemma \ref{lem: effective dimension}, we have that
  $$\frac{c_{\lambda_n}}{n\lambda_n}\lesssim \frac{\ln\big(1/\psi_{1/2}^{-1}(\lambda_n)\big)}{n\lambda_n}\to 0.$$
  The first inequality of \eqref{eq: set of equations} will therefore be satisfied eventually.
  On the event of Theorem \ref{th: meunier optimal}, we have that
  $$\frac{\lVert M-\widehat{M}\rVert_{L^2(Z,\mathcal{H}_{K_X})}}{\lambda_n}\lesssim \frac{m^{-\frac{\beta_Z}{2(\beta_Z+p_Z)}}}{\lambda_n}\lesssim \frac{n^{-\frac{\theta \beta_Z}{2(\beta_Z+p_Z)}}}{\lambda_n}\leq \frac{n^{-r(\theta)/2}}{\lambda_n} \to 0$$ and similarly
  $$\frac{\lVert M-\widehat{M}\rVert_{\alpha_Z}}{\lambda_n \sqrt{n}}\lesssim \frac{m^{-\frac{\beta_Z-\alpha_Z}{2(\beta_Z+p_Z)}}}{\sqrt{n}\lambda_n}\lesssim \frac{n^{-\theta\frac{\beta_Z-\alpha_Z}{2(\beta_Z+p_Z)}-\frac{1}{2}}}{\lambda_n}\leq \frac{n^{-r(\theta)/2}}{\lambda_n} \to 0.$$
  Hence, inequalities 2 and 3 of \eqref{eq: set of equations} will eventually be satisfied.
  The last condition is obtained as a consequence of the second and third. We now turn to the lower bound. Let us define the terms
  $$w(\lambda):=\big[\psi_{1/2}^{-1}(\lambda)\big]^{\beta_X-\gamma},\qquad v_n(\lambda):=\frac{[\psi_{1/2}^{-1}(\lambda)]^{1-\gamma-p_X}}{\lambda n},\qquad d_n(\lambda):=\frac{[\psi_{1/2}^{-1}(\lambda)]^{\nu}\,n^{-r(\theta)}}{\lambda^2},$$
so that $R_n\geq w+v_n+d_n$. Since $\psi_{1/2}^{-1}$ is nondecreasing, $w$ is nondecreasing. With $t=\psi^{-1}_{1/2}(\lambda)$ we have $v_n(\lambda)=t^{-\gamma-p_X}/(n\phi(t^2))$ and $d_n(\lambda)=n^{-r(\theta)} t^{\nu-2}/\phi(t^2)^2$. As $\gamma+p_X>0$, $\nu\le1$ and $\phi$ is increasing, both are nonincreasing in $t$, hence in $\lambda$. 
We will also use the following elementary fact.

\smallskip
\noindent($\ast$) \emph{Let $f,g\ge0$ on $(0,\psi_{1/2}(\lVert \mathcal{L}_X\rVert )]$, with $f$ nondecreasing and $g$ nonincreasing, and let $\lambda^*$ satisfy $f(\lambda^*)=g(\lambda^*)$. Then $f(\lambda)+g(\lambda)\ge f(\lambda^*)$ for every $\lambda\in(0,\psi_{1/2}(\lVert \mathcal{L}_X\rVert )]$.}
Indeed, if $\lambda\ge\lambda^*$, then $f(\lambda)+g(\lambda)\ge f(\lambda)\ge f(\lambda^*)$. If $\lambda<\lambda^*$, then $f(\lambda)+g(\lambda)\ge g(\lambda)\ge g(\lambda^*)=f(\lambda^*)$.
\smallskip

Let $n$ be large enough that $n^{-1}$ and $n^{-r(\theta)/2}$ lie in the ranges of $\Theta_1$ and $\Theta_2$, and that $\ln n/n\le \psi_{1/2}(\lVert \mathcal{L}_X\rVert )$.

\emph{(a)} Set $t_n^{\mathrm v}:=\Theta_1^{-1}(n^{-1})$ and $\lambda^{\mathrm v}:=\psi_{1/2}(t_n^{\mathrm v})$. A direct computation gives $v_n(\psi_{1/2}(t))/w(\psi_{1/2}(t))=1/(n\Theta_1(t))$. Since $n\Theta_1(t_n^{\mathrm v})=1$, we get $w(\lambda^{\mathrm v})=v_n(\lambda^{\mathrm v})$. By ($\ast$) with $f=w$ and $g=v_n$, for every $\lambda\in(0,\psi_{1/2}(\lVert \mathcal{L}_X\rVert )]$,
$$R_n(\lambda)\ \ge\ w(\lambda)+v_n(\lambda)\ \ge\ w(\lambda^{\mathrm v})=\big[\Theta_1^{-1}(n^{-1})\big]^{\beta_X-\gamma}.$$

\emph{(b)} Let $\lambda^{0}:=\psi_{1/2}\big(\Theta_2^{-1}(n^{-r(\theta)/2})\big)$. Similarly, $d_n(\psi_{1/2}(t))/w(\psi_{1/2}(t))=n^{-r(\theta)}/\Theta_2(t)^2$, so $w(\lambda^0)=d_n(\lambda^0)$. 
By ($\ast$) with $f=w$ and $g=d_n$, for every $\lambda\in(0,\psi_{1/2}(\lVert \mathcal{L}_X \rVert )]$,
$$R_n(\lambda)\ \ge\ w(\lambda)+d_n(\lambda)\ \ge\ w(\lambda^0)=\big[\Theta_2^{-1}(n^{-r(\theta)/2})\big]^{\beta_X-\gamma}.$$

\emph{(c)} Let $(\lambda_n)$ be as in (i) and $t_n:=\psi_{1/2}^{-1}(\lambda_n)$. We claim that $t_n>s_n^*$ for all $n$ large enough. Otherwise, $t_n\le s_n^*$ for infinitely many $n$. For each such $n$, since $\psi$ is increasing,
$$n\lambda_n=n\psi(t_n)\le n\psi(s_n^*)=\ln(1/s_n^*)\le\ln(1/t_n),$$
so $\ln(1/t_n)/(n\lambda_n)\ge1$. This contradicts $\ln(1/t_n)=o(n\lambda_n)$. Since $w$ is increasing, for $n$ large,
$$R_n(\lambda_n)\ \ge\ w(\lambda_n)=t_n^{\beta_X-\gamma}\ \ge\ (s_n^*)^{\beta_X-\gamma}.$$
Taking the maximum of the bounds in (a), (b) and (c) gives the lower bound in (i).

\textbf{(Part (ii))}
We start by checking that $\lambda_n^*$ is well-defined and satisfies the conditions of (i).
Under the condition $\omega_n \ln n=o(n)$, one can easily check that $t_n^*$ exists and is unique for $n$ large enough and satisfies $t_n^*\to 0$.
The other two terms in the maximum in \eqref{eq: lambda} can easily be shown to tend towards $0$.
Hence, $\lambda_n^*\to 0$.
Next, since $\Theta_2(t)=t^{\xi/2}\psi_{1/2}(t)$, we have that
$$\frac{n^{-r(\theta)}}{(\lambda_n^*)^2}\leq \frac{n^{-r(\theta)}}{\psi_{1/2}\big(\Theta_2^{-1}\big(l_n^{1/2} n^{-r(\theta)/2}\big)\big)^2}=\frac{\Theta_2^{-1}\big(l_n^{1/2} n^{-r(\theta)/2}\big)^\xi}{l_n}\to 0.$$
For the remaining condition we obtain
$$\frac{\ln(1/\psi_{1/2}^{-1}(\lambda_n^*))}{n\lambda_n^*}\leq \frac{\ln(1/t_n^*)}{n\psi_{1/2}(t_n^*)}=\frac{1}{\omega_n}\to 0.$$
We now derive the upper-bound.
First notice that from (a) we have that
$$v_n\circ \psi_{1/2}(\Theta_1^{-1}(n^{-1}))=w\circ \psi_{1/2}(\Theta_1^{-1}(n^{-1}))$$and similarly as in (b) one can show
$$l_n (d_n\circ \psi_{1/2})\big(\Theta_2^{-1}((l_n n^{-r(\theta)})^{1/2})\big)=w\circ \psi_{1/2}\big(\Theta_2^{-1}((l_n n^{-r(\theta)})^{1/2})\big).$$
Since, $v_n$ and $d_n$ are nonincreasing, $w$ nondecreasing and $l_n\geq 1$, we have that
\begin{equation*}
  \begin{split}
    R_n(\lambda_n^*)&\leq w(\lambda_n^*)+v_n\circ \psi_{1/2}\big(\Theta_1^{-1}(n^{-1})\big)+d_n\circ \psi_{1/2}\big(\Theta^{-1}_2((l_n n^{-r(\theta)})^{1/2})\big)+\frac{[\psi_{1/2}^{-1}(\lambda^*_n)]^{\,a(1-\gamma)-p_X}}{n^{1+r(\theta)}(\lambda^*_n)^{3}}\\
    &\leq w(\lambda_n^*)+w\circ \psi_{1/2}\big(\Theta_1^{-1}(n^{-1})\big)+w\circ \psi_{1/2}\big(\Theta^{-1}_2((l_n n^{-r(\theta)})^{1/2})\big)+\frac{[\psi_{1/2}^{-1}(\lambda^*_n)]^{\,a(1-\gamma)-p_X}}{n^{1+r(\theta)}(\lambda^*_n)^{3}}\\
    &\leq 3\;w(\lambda_n^*)+\frac{[\psi_{1/2}^{-1}(\lambda^*_n)]^{\,a(1-\gamma)-p_X}}{n^{1+r(\theta)}(\lambda^*_n)^{3}}.
  \end{split}
\end{equation*}
Now since,
\begin{equation*}
  \begin{split}
    \frac{[\psi_{1/2}^{-1}(\lambda^*_n)]^{\,a(1-\gamma)-p_X}}{n^{1+r(\theta)}(\lambda^*_n)^{3}}/d_n(\lambda_n^*)=[\psi_{1/2}^{-1}(\lambda^*_n)]^{(\beta_X-1)_+}\cdot \frac{v_n(\lambda_n^*)}{w(\lambda_n^*)}&\leq[\psi_{1/2}^{-1}(\lambda^*_n)]^{(\beta_X-1)_+}\cdot \frac{v_n(\lambda^{\mathrm v})}{w(\lambda^{\mathrm v})}\\
    &=[\psi_{1/2}^{-1}(\lambda^*_n)]^{(\beta_X-1)_+}\\
    &\leq 1\vee \lVert \mathcal{L}_X\rVert^{\beta_X-1} 
  \end{split}
\end{equation*}and
$$d_n(\lambda_n^*)\leq w(\lambda_n^*),$$
we obtain
$$R_n(\lambda_n^*)\lesssim w(\lambda_n^*)$$
The upperbound now follows easily once we have shown $t_n^*\leq \psi_{1/2}^{-1}\Big(\frac{\omega_n\Lambda_n}{n}\Big)$, since then
$$w(\psi_{1/2}(t_n^*)) = (t_n^*)^{\beta_X-\gamma} \leq \psi_{1/2}^{-1}\Big(\frac{\omega_n\Lambda_n}{n}\Big).$$

For $n$ large, $\omega_n\ge1$, so
$n\psi_{1/2}(t_n^*)=\omega_n\ln(1/t_n^*)\ge\ln(1/t_n^*)$. Since
$t\mapsto n\psi_{1/2}(t)-\ln(1/t)$ is increasing and vanishes at $s_n^*$, this gives
$t_n^*\ge s_n^*$, hence $\ln(1/t_n^*)\le\ln(1/s_n^*)\le\Lambda_n$. Therefore
$$t_n^*=\psi_{1/2}^{-1}\Big(\frac{\omega_n\ln(1/t_n^*)}{n}\Big)\le\psi_{1/2}^{-1}\Big(\frac{\omega_n\Lambda_n}{n}\Big).$$
\end{proof}

\subsection{Proof of Corollary \ref{cor:theta1}}
\begin{proof}
Write $\bar t:=\lVert \mathcal{L}_X \rVert $, $\bar\phi:=\phi(\bar t^{\,2})$ and notice that $r:=r(1)=\beta_Z/(\beta_Z+p_Z)<1$. For $t\in(0,\bar t\,]$ let
$$A(t):=t^{\beta_X-\gamma},\qquad V_n(t):=\frac{t^{1-\gamma-p_X}}{n\,\psi(t)},\qquad D_n(t):=\frac{t^{\nu}}{n^{r}\,\psi(t)^2}.$$
Notice that $A/V_n=n\,\Theta_1$ and $A/D_n=n^{r}\,\Theta_2^2$. Let
$$t_n^{1}:=\Theta_1^{-1}(n^{-1})\qquad\text{and}\qquad t_n^{2}:=\Theta_2^{-1}\big(\ell_n^{1/2}n^{-r/2}\big),$$
so that $\lambda_n^*=\psi\big(\max\{t_n^{1},t_n^*,t_n^{2}\}\big)$.

For every $t\in(0,\bar t\,]$, using $\psi(t)=t\,\phi(t^2)$,
\begin{equation}\label{eq:corollary}
  \frac{D_n(t)}{V_n(t)}=\frac{n^{1-r}\,t^{\nu}}{t^{1-\gamma-p_X}\,\psi(t)}=\frac{n^{1-r}\,t^{\nu+\gamma+p_X-2}}{\phi(t^2)}\ \ge\ c\,n^{1-r},\qquad c:=\frac{\bar t^{\,\nu+\gamma+p_X-2}}{\bar\phi}.
\end{equation}
This holds because $\nu\le1-\gamma$ gives $\nu+\gamma+p_X-2\le p_X-1\le0$, so that $t^{\nu+\gamma+p_X-2}\ge\bar t^{\,\nu+\gamma+p_X-2}$, and because $\phi(t^2)\le\bar\phi$.

Let $n$ be large enough such that $\ell_n\ge1$. Since $\Theta_2$ is increasing and $A/D_n=n^{r}\Theta_2^2$, for every $t\in(0,\bar t\,]$, we have
$$t<t_n^{2}\iff n^{r}\,\Theta_2(t)^2<\ell_n\iff\frac{A(t)}{D_n(t)}<\ell_n.$$
In particular, $A(t)<D_n(t)$ implies $t<t_n^{2}$.
This will now be used to show $t_n^1<t_n^2$ and $t_n^*<t_n^2$.

By $A/V_n=n \Theta_1$, we have that $A(t_n^{1})=V_n(t_n^{1})$. Equation \eqref{eq:corollary} then gives
$$\frac{A(t_n^{1})}{D_n(t_n^{1})}=\frac{V_n(t_n^{1})}{D_n(t_n^{1})}\le\frac{n^{r-1}}{c}\longrightarrow0,$$
so $t_n^{1}<t_n^{2}$ for $n$ large enough.

Since $\Theta_2(t)^2=t^{2+\xi}\phi(t^2)^2=t^{\xi}\psi_{1/2}(t)^2$ and $n\psi_{1/2}(t_n^*)=\omega_n\ln(1/t_n^*)$,
$$\frac{A(t_n^*)}{D_n(t_n^*)}=n^{r}\,(t_n^*)^{\xi}\,\psi_{1/2}(t_n^*)^2=n^{r-2}\,(t_n^*)^{\xi}\,\big(\omega_n\ln(1/t_n^*)\big)^2.$$
By the proof of Proposition~\ref{prop:rates}, $\ln(1/t_n^*)\le\Lambda_n$. Since $\psi_{1/2}(t)\le\bar\phi\,t$, we have (by setting $t=1/(\overline{\phi}n)$) that $\psi_{1/2}^{-1}(1/n)\ge1/(\bar\phi n)$, hence $\Lambda_n\le\ln(\bar\phi n)$. 
It follows that $\ln(1/t_n^*)\leq \Lambda_n \leq \ln(\overline{\phi}n)$, and since $\xi\geq 0$ we have
$$\frac{A(t_n^*)}{D_n(t_n^*)}\le\bar t^{\,\xi}\Big(\frac{\omega_n\ln(\bar\phi n)}{n^{1-r/2}}\Big)^2\longrightarrow0$$
by the assumption on $\omega_n$. So $t_n^* \leq t_n^{2}$ for $n$ large enough.

Hence, $\max\{t_n^{1},t_n^*,t_n^{2}\}=t_n^{2}$ for $n$ large, which gives the formula for $\lambda_n^*$.
And since, from the proof of Proposition \ref{prop:rates}, $R_n(\lambda_n^*)\leq 3 \big[\psi_{1/2}^{-1}(\lambda_n^*)\big]^{\beta_X-\gamma} $ the result follows.
\end{proof}

\newpage

\subsection{Useful Lemmas and Propositions}
\begin{lemma}\label{lem: properties unbounded operator}
  Let $C:\mathbb{H}\to \mathbb{H} $ be a bounded, injective, self-adjoint, positive operator on a separable real Hilbert space $\mathbb{H} $.
  Then $C^{-1}:\mathcal{R}(C)\to \mathbb{H} $ is a densely defined, self-adjoint, positive and coercive operator.
  Furthermore, the coercivity holds with constant $\lVert C\rVert^{-1} $.
\end{lemma}
\begin{proof}
  Let $x=Cu$ and $y=Cv$ be both in $\mathcal{R}(C)=\mathcal{D}(C^{-1})$.

  \textit{Densely defined.}
  This follows easily, since $$\overline{\mathcal{D}(C^{-1})}=\overline{\mathcal{R}(C)}=\mathcal{N}(C)^\perp = \mathbb{H} .$$

  \textit{Self-adjoint.}
  First notice that
  $$\langle C^{-1}x,y\rangle_\mathbb{H}=\langle u,Cv\rangle_\mathbb{H}=\langle Cu,v\rangle_\mathbb{H}=\langle x,C^{-1}y\rangle_\mathbb{H}.$$
  This shows that $C^{-1}$ is symmetric. It now suffices to show that $\mathcal{D}((C^{-1})^*)\subseteq \mathcal{D}(C^{-1})$.
  Let $z\in \mathcal{D}((C^{-1})^*)$ and set $h=(C^{-1})^*z$. For every $w\in \mathbb{H}$ we have that
  $$\langle w,z\rangle_\mathbb{H} = \langle C^{-1}C w,z\rangle_\mathbb{H}=\langle Cw,h\rangle_\mathbb{H}=\langle w,Ch\rangle_\mathbb{H}. $$
  Hence, $z=Ch \in \mathcal{R}(C)=\mathcal{D}(C^{-1})$.

  \textit{Positive.}
  $$\langle C^{-1}x,x\rangle_\mathbb{H}=\langle u,Cu\rangle_\mathbb{H}=\lVert C^{1/2}u\rVert_\mathbb{H}^2\geq 0.$$

  \textit{Coercive.}
  $$\lVert x\rVert^2_\mathbb{H}\leq \lVert C^{1/2}\rVert^2 \lVert C^{1/2}u\rVert^2_\mathbb{H} = \lVert C\rVert  \langle Cu,u\rangle_\mathbb{H} =\lVert C\rVert  \langle x,C^{-1}x\rangle_\mathbb{H}.$$
\end{proof}

\begin{lemma}\label{lem: inclusion empirical cov}
  Let $X_1,\dots,X_n$ be centered i.i.d. random elements defined on a separable Hilbert $H$ space with covariance operator $\mathcal{C}_X$ and let $\widehat{\mathcal{C}}_X:=(1/n)\sum_{i=1}^{n}X_i \otimes X_i$.
  Then, $\overline{\mathcal{R}(\widehat{\mathcal{C}}_X)}\subseteq \overline{\mathcal{R}(\mathcal{C}_X)}$ almost surely.
\end{lemma}
\begin{proof}
  Let $u\in \mathcal{N}(\mathcal{C}_X)$, it then suffices to show that $u\in \mathcal{N}(\widehat{\mathcal{C}}_X)$ almost surely.
  Since
  $$\mathbb{E}[\lVert \widehat{\mathcal{C}}_X^{1/2} u\rVert_H^2 ]=\mathbb{E}[\langle u,\widehat{\mathcal{C}}_Xu\rangle_H ]=\langle u,\mathcal{C}_Xu\rangle_H=0$$
  and $\mathcal{N}(\widehat{\mathcal{C}}_X^{1/2})=\mathcal{N}(\widehat{\mathcal{C}}_X)$, the claim follows.
\end{proof}

\begin{lemma}\label{lem: partial isometry}
  There exists a partial isometry $W:S_2(\mathcal{H}_{k_X},\mathbb{Y} )\to S_2(L_2(X),\mathbb{Y})$ with initial space $\N(\mathbf{I_X})^\perp$ and final space $S_2(L_2(X),\mathbb{Y})$, by injectivity of $\mathcal{L}_X$, such that $\mathbf{I_X}=A_{\mathcal{L}_X^{1/2}}W$.
  Furthermore, one has that $$\mathbf{I_X}(\N(A_{\widetilde{T}^*})^\perp)\subseteq A_{\mathcal{L}_X^{1/2}}(\N(Q)^\perp).$$
\end{lemma}
\begin{proof}
  The first part of the lemma is a direct application of the polar decomposition theorem and using the fact that $\mathbf{I_X}\mathbf{I_X}^*=A_{\mathcal{L}_X}$.
  For the second part, by the first part, it suffices to show that $W(\N(A_{\widetilde{T}^*})^\perp)\subseteq \N(Q)^\perp$.
  For this, let $f\in \N(A_{\widetilde{T}^*})^\perp$ and $g\in \N(Q)$.
  First notice that $\N(\mathbf{I_X})\subseteq \N(A_{\widetilde{T}^*})$, and therefore $f\in \N(A_{\widetilde{T}^*})^\perp \subseteq \N(\mathbf{I_X})^\perp$.  
  Using the fact that $W^*W=P_{\N(\mathbf{I_X})^\perp}$, $WW^*=P_{\overline{\mathcal{R}(\mathbf{I_X})}}$ and $Wf\in \overline{\Ra(\mathbf{I_X})}$ we obtain
  $$\langle Wf,g\rangle_{S_2(L_2(X),\mathbb{Y})} = \langle W f,WW^*g\rangle_{S_2(L_2(X),\mathbb{Y})}=\langle f,W^* g\rangle_{S_2(L_2(X),\mathbb{Y})}.$$
  It now remains to show that $W^* g \in \N(A_{\tilde{T}^*})$, which holds since by $\mathbf{I_X}= A_{\mathcal{L}_X^{1/2}}W$ we obtain
  $$A_{\widetilde{T}^*}W^*g = A_{T^*} \mathbf{I_X}W^* g=A_{T^*} A_{\mathcal{L}_X^{1/2}}P_{\overline{\Ra(\mathbf{I_X})}} g=A_{T^*} A_{\mathcal{L}_X^{1/2}}P_{\N(A_{\mathcal{L}_X^{1/2}})^\perp} g=A_{T^*} A_{\mathcal{L}_X^{1/2}} g = Qg=0.$$
\end{proof}
\begin{lemma}\label{lem: membership}
  One has that $\mathbf{I_X}\C_{MY}g_\lambda(\C_M)\in A_{\mathcal{L}_X^{1/2}}(\N(Q)^\perp)$ and $\mathbf{I_X}\widehat{\C}_{MY}g_\lambda(\widehat{\C}_M)\in A_{\mathcal{L}_X^{1/2}}(\N(Q)^\perp)$ $\mathbb{P} $-almost surely.
\end{lemma}
We first recall the following notation: for a Hilbert space $\mathbb{H}$ and natural number $n$, we denote by $\mathbb{H}^n$ the Hilbert space equipped with the inner product $\langle \cdot, \cdot \rangle_{\mathbb{H}^n} = n^{-1} \sum_{i=1}^n \langle \cdot, \cdot \rangle_{\mathbb{H}}$.
\begin{proof}
  By Proposition \ref{prop:problem is hilbert scaled} we have that
  $$\mathbf{I_X}\C_{MY}g_\lambda(\C_M)=A_{\mathcal{L}_X^{1/2}}g_\lambda\left(A_{\mathcal{L}_X^{1/2}}A_T A_{T^*}A_{\mathcal{L}_X^{1/2}}\right)A_{\mathcal{L}_X^{1/2}}A_Tr=A_{\mathcal{L}_X^{1/2}}g_\lambda\left(Q^*Q\right)Q^*r=A_{\mathcal{L}_X^{1/2}}Q^*g_\lambda\left(QQ^*\right)r.$$
  Clearly, $Q^*g_\lambda(QQ^*)r\in \overline{\Ra(Q^*)}=\N(Q)^\perp$, showing the first part.
  For the second part, we introduce the operator $$\widehat{T}:S_2(\mathcal{H}_{k_X},\mathbb{Y} )\to \mathbb{Y}^{n}:C\mapsto (CM(\widetilde{Z}_i))_{i=1}^{n}.$$
  One can then verify that $$\widehat{T}^* :\mathbb{Y}^{n}\to S_2(\mathcal{H}_{k_X},\mathbb{Y} ):(y_i)_{i=1}^{n}\mapsto n^{-1}\sum_{i=1}^{n} M(\widetilde{Z}_i)\otimes y_i,$$
  and therefore $\widehat{T}^*\widehat{T}=A_{\widehat{\C}_M}$ and $\widehat{T}^*(\widetilde{Y}_i)_{i=1}^{n}=\widehat{\C}_{YM}$.
  As a consequence, this allows us to write$$\mathbf{I_X}\widehat{\C}_{YM}g_\lambda(\widehat{\C}_M)=\mathbf{I_X}g_\lambda(\widehat{T}^* \widehat{T})\widehat{T}^*(\widetilde{Y}_i)_{i=1}^{n}=\mathbf{I_X}\widehat{T}^*g_\lambda(\widehat{T} \widehat{T}^*)(\widetilde{Y}_i)_{i=1}^{n}.$$
  By Lemma \ref{lem: partial isometry}, it now remains to show that $\overline{\Ra(\widehat{T})}\subseteq \N(A_{\widetilde{T}}^*)^\perp$, which is equivalent to showing that $\overline{\Ra(A_{\widehat{\C}_M})}\subseteq \overline{\Ra(A_{\C_M})}$ $\mathbb{P}-$a.e since $\widehat{T}^*\widehat{T}=A_{\widehat{\C}_M}$ and $A_{\widetilde{T}}A_{\widetilde{T}^*}=A_{\C_M}$.
  Since $A_{\widehat{\C}_M}\cong \widehat{\C}_M \otimes \Id_\mathbb{Y} $ and $A_{\C_M}\cong \C_M \otimes \Id_\mathbb{Y} $, it suffices to show that $\overline{\Ra(\widehat{\C}_M)}\subseteq \overline{\Ra(\C_M)}$ $\mathbb{P}-$a.e. (see Lemma \ref{lem: properties bounded tensor op}.2). 
  This is shown by Lemma \ref{lem: inclusion empirical cov}.
\end{proof}

\begin{lemma}\label{lem: N IX membership}
  One has that $K_\lambda \in \mathcal{N}(\mathbf{I_X})^\perp$ and $\overline{K}_\lambda,\widehat{K}_\lambda \in \mathcal{N}(\mathbf{I_X})^\perp$ almost surely.
\end{lemma}
\begin{proof}
  For $K_\lambda$ and $\overline{K}_\lambda$, by the proof of the previous lemma, it suffices to show that $$\mathbf{I_X}(\mathcal{N}(A_{\widetilde{T}^*})^\perp)\subseteq \mathbf{I_X}(\mathcal{N}(\mathbf{I_X})^\perp),$$
  which is true since $A_{\widetilde{T}^*}=A_{T^*}\mathbf{I_X}$.
  
We now turn to the proof of $\widehat{K}_\lambda$ which is divided into five steps.

  \emph{Step 0: characterisation of $\mathcal{N}(\mathbf{I_X})^\perp$.} Let $P$ denote the orthogonal projection onto
$\mathcal{N}(\mathcal{I}_X)^\perp\subseteq\mathcal{H}_{k_X}$. The operator $A_{P}:K\mapsto KP$ is idempotent and self-adjoint on
$S_2(\mathcal{H}_{k_X},\mathbb{Y})$, since $\langle KP,L\rangle_{S_2}=\operatorname{tr}(PK^*L)=\operatorname{tr}(K^*LP)=\langle K,LP\rangle_{S_2}$;
hence it is an orthogonal projection. Moreover, $\mathbf{I_X}K=K\mathcal{I}_X^*=0$ if and only if $K$ vanishes on
$\overline{\mathcal{R}(\mathcal{I}_X^*)}=\mathcal{N}(\mathcal{I}_X)^\perp$, i.e.\ $\mathcal{N}(\mathbf{I_X})=\mathcal{N}(A_{P})$. Therefore
$$\mathcal{N}(\mathbf{I_X})^\perp=\mathcal{R}(A_{P})=\{K\in S_2(\mathcal{H}_{k_X},\mathbb{Y}):\ Kh=0\ \text{for all } h\in\mathcal{N}(\mathcal{I}_X)\}.$$

\emph{Step 1: $\Phi_X(\overline{X}_j)\in\mathcal{N}(\mathcal{I}_X)^\perp$ for all $j=1,\dots,m$.} If $h\in\mathcal{N}(\mathcal{I}_X)$, then $h=0$
$\mathbb{P}_X$-almost everywhere, so that, by the reproducing property, $\langle h,\Phi_X(\overline{X}_j)\rangle_{\mathcal{H}_{k_X}}=h(\overline{X}_j)=0$
almost surely for each $j$. Since $\mathcal{N}(\mathcal{I}_X)$ is a closed subspace of the separable space $\mathcal{H}_{k_X}$, it admits a
countable dense subset $\{h_k\}_{k\in\mathbb{N}}$, and the event
$$\Omega_0:=\bigcap_{j=1}^{m}\bigcap_{k\in\mathbb{N}}\{h_k(\overline{X}_j)=0\}$$
has probability one as a countable intersection of events of probability one. On $\Omega_0$, for any $h\in\mathcal{N}(\mathcal{I}_X)$ and
$h_{k_l}\to h$, we have $|h(\overline{X}_j)|=|\langle h-h_{k_l},\Phi_X(\overline{X}_j)\rangle|\leq\lVert h-h_{k_l}\rVert\,k_X(\overline{X}_j,\overline{X}_j)^{1/2}\to 0$.
Hence, on $\Omega_0$, $\Phi_X(\overline{X}_j)\in\mathcal{N}(\mathcal{I}_X)^\perp$ for all $j=1,\dots,m$ since $h(\overline{X}_j)=\left\langle \Phi_X(\overline{X}_j),h\right\rangle$.

\emph{Step 2: $\widehat{M}(z)\in \mathcal{N}(\mathcal{I}_X)^\perp$ for all $z\in E_Z$.} For every $z\in E_Z$, setting
$u:=g_\varepsilon(\widehat{\C}_Z)\Phi_Z(z)\in\mathcal{H}_{k_Z}$, we have
$$\widehat{M}(z)=\widehat{\C}_{ZX}u=\frac{1}{m}\sum_{j=1}^m\langle u,\Phi_Z(\overline{Z}_j)\rangle_{\mathcal{H}_{k_Z}}\Phi_X(\overline{X}_j)
\in\operatorname{span}\{\Phi_X(\overline{X}_j)\}_{j=1}^m .$$
Hence, on $\Omega_0$, $\widehat{M}(z)\in\mathcal{N}(\mathcal{I}_X)^\perp$ for all $z\in E_Z$, and in particular for $z=\widetilde{Z}_1,\dots,\widetilde{Z}_{n}$.

\emph{Step 3: $\mathcal{N}(\mathcal{I}_X)$ is annihilated by $\widehat{\C}_{\hat{M}}$ and $\widehat{\C}_{\hat{M}Y}$.} Let $h\in\mathcal{N}(\mathcal{I}_X)$. On $\Omega_0$,
$\langle h,\widehat{M}(\widetilde{Z}_i)\rangle=0$ for all $i$, hence
$$\widehat{\C}_{\hat{M}}h=\frac{1}{n}\sum_{i=1}^{n}\langle h,\widehat{M}(\widetilde{Z}_i)\rangle\widehat{M}(\widetilde{Z}_i)=0
\qquad\text{and}\qquad
\widehat{\C}_{\hat{M}Y}h=\frac{1}{n}\sum_{i=1}^{n}\langle h,\widehat{M}(\widetilde{Z}_i)\rangle\widetilde{Y}_i=0 .$$

\emph{Step 4: conclusion.} The operator $\widehat{\C}_{\hat{M}}$ is positive, self-adjoint and of finite rank, so that
$\widehat{\C}_{\hat{M}}=\sum_{k=1}^r\sigma_k\, e_k\otimes e_k$ with $\sigma_k>0$ and orthonormal $e_k$, and
$$g_\lambda(\widehat{\C}_{\hat{M}})=\sum_{k=1}^r g_\lambda(\sigma_k)\, e_k\otimes e_k+g_\lambda(0)\,P_{\mathcal{N}(\widehat{\C}_{\hat{M}})} .$$
By Step 3, $h\in\mathcal{N}(\widehat{\C}_{\hat{M}})$, so that $g_\lambda(\widehat{\C}_{\hat{M}})h=g_\lambda(0)h$ and
$$\widehat{K}_\lambda h=g_\lambda(0)\,\widehat{\C}_{\hat{M}Y}h=0 .$$
Since $h\in\mathcal{N}(\mathcal{I}_X)$ was arbitrary, Step 0 yields $\widehat{K}_\lambda\in\mathcal{N}(\mathbf{I_X})^\perp$ on $\Omega_0$.
\end{proof}

\begin{lemma}\label{lem: Li extended}
  For any $K\in S_2(\mathcal{H}_{k_X},\mathbb{Y} )$ and $0\leq \gamma \leq 1$, one has that $$\left\lVert \mathbf{I_X}K\right\rVert_{S_2([\mathcal{H}_{k_X}]^\gamma,\mathbb{Y} )} \leq \left\lVert K\C_X^{(1-\gamma)/2}\right\rVert_{S_2(\mathcal{H}_{k_X},\mathbb{Y} )},$$
  where the inequality holds if $\gamma<1$ or if $\gamma=1$ and $K\in \N(\mathbf{I_X})^\perp$.
\end{lemma}
\begin{proof}
  By Lemma 2 in \cite{li2022optimal}, it only remains to prove the equality for $\gamma=1$.
  This follows from Lemma \ref{lem: partial isometry}, which gives
  $$\left\lVert \mathbf{I_X}K\right\rVert_{S_2([\mathcal{H}_{k_X}]^1,\mathbb{Y} )} = \left\lVert A_{\mathcal{L}_X^{1/2}}WK\right\rVert_{S_2([\mathcal{H}_{k_X}]^1,\mathbb{Y} )}=\left\lVert WK\right\rVert_{S_2(\mathcal{H}_{k_X},\mathbb{Y} )}=\left\lVert K\right\rVert_{S_2(\mathcal{H}_{k_X},\mathbb{Y} )},$$
  where we used the fact that $(W^*W)^{1/2}=P_{\mathcal{N}(\mathbf{I_X})^\perp}^{1/2}=P_{\mathcal{N}(\mathbf{I_X})^\perp}$.
\end{proof}

\begin{lemma}\label{lem: unitary equivalence Hkx}
Let $P:=W^*W$ be the orthogonal projection onto $\mathcal{N}(\mathbf{I_X})^\perp$ and let $A$ be a bounded self-adjoint operator
on $S_2(\mathcal{H}_{k_X},\mathbb{Y})$ such that $PA=A$. Set $S:=WAW^*$. Then $A=W^*SW$ and, for every bounded Borel function
$g$ on $\mathbb{R}$,
$$g(A)=W^*g(S)W+g(0)(Id-P)\qquad\text{and}\qquad g(A)W^*=W^*g(S).$$
\end{lemma}

\begin{proof}
Set $\mathcal{H}_0:=\mathcal{R}(P)=\mathcal{N}(\mathbf{I_X})^\perp$. By Lemma \ref{lem: partial isometry}, $W$ is a partial isometry with initial
space $\mathcal{H}_0$ and final space $S_2(L^2(X),\mathbb{Y})$.
Hence $U:=W|_{\mathcal{H}_0}:\mathcal{H}_0\to S_2(L^2(X),\mathbb{Y})$ is unitary, $WW^*=\Id$, $\mathcal{R}(W^*)=\mathcal{H}_0$, $Wh=WPh+W(Id-P)h=WPh$, and
$W^*x=U^{-1}x=U^*x$ for all $x\in S_2(L^2(X),\mathbb{Y})$.

Since $A$ is self-adjoint, $AP=(PA)^*=A$, so that $PA=AP=A$. Hence $\mathcal{H}_0$ and $\mathcal{H}_0^\perp$ are both
invariant under $A$, and $Ah=APh=0$ for $h\in\mathcal{H}_0^\perp$. With respect to $\mathcal{H}_0\oplus\mathcal{H}_0^\perp$ we may thus write
$A=A_0\oplus 0$, where $A_0:=A|_{\mathcal{H}_0}$. Moreover, for $h\in\mathcal{H}_0$,
$$U^*SUh=W^*WAW^*Wh=PAPh=Ah=A_0h,$$
i.e.\ $A_0=U^*SU$.

By the functional calculus for self-adjoint operators, we obtain that
$$g(A_0)=g(U^*SU)=U^*g(S)U \qquad \text{and}\qquad g(A)=g(A_0\oplus 0)=g(A_0)\oplus g(0)\,Id_{\mathcal{H}_0^\perp}.$$
Hence, for $h=h_0+h_1\in\mathcal{H}_0\oplus\mathcal{H}_0^\perp$, we obtain
$$g(A)h=U^*g(S)Uh_0+g(0)h_1=W^*g(S)Wh+g(0)(Id-P)h .$$

Finally, since $WW^*=\Id$ and $(\Id-P)W^*=0$, we get $g(A)W^*=W^*g(S)$. Taking $g(t)=t$ gives $A=W^*SW$.
\end{proof}

\begin{lemma}\label{lem: unitary equivalence L_2}
Let $A$ be a bounded self-adjoint operator on $S_2(L_2(X),\mathbb{Y})$ and set $S:=W^*AW$. 
Then $A=WSW^*$ and, for every bounded Borel function
$g$ on $\mathbb{R}$,
$$g(A)=Wg(S)W^*\qquad\text{and}\qquad g(A)W=W g(S).$$
\end{lemma}

\begin{proof}
  The proof is similar (but shorter and simpler) to the proof of Lemma \ref{lem: unitary equivalence Hkx} and is therefore omitted.
\end{proof}

The following inequality will be the key to obtain convergence rates.
\begin{lemma}\label{lem: inequality core}
  Under Assumption \ref{ass: link core}, we have for $\gamma \in [0,1]$ that
  $$\left\lVert [K]_X \right\rVert_\gamma \leq \lambda^{-1/2}\left\lVert K\mathcal{C}_{M,\lambda}^{1/2}\right\rVert_{S_2}\left[\psi^{-1}_{1/2}\left(\frac{\lambda}{m^2}\right)\right]^{\frac{1-\gamma}{2}} $$
  for any $[K]_X\in A_{\mathcal{L}_X^{1/2}}(\N(Q)^\perp)$ with $K\in S_2(\mathcal{H}_{k_X},\mathbb{Y} )$.
\end{lemma}
\begin{proof}
  Since $Q^*=A_{T\mathcal{L}_X^{1/2}}\cong (\mathcal{L}_X^{1/2}T^*)\otimes \Id_\mathbb{Y} $,
  we have by Lemma \ref{lem: properties bounded tensor op}.2 that $$\mathcal{N}(Q)^\perp = \overline{\mathcal{R}(Q^*)} \cong \overline{\mathcal{R}(\mathcal{L}_X^{1/2}T^*)}\otimes \mathbb{Y}=\N(T\mathcal{L}_X^{1/2})^\perp \otimes \mathbb{Y} .$$
  This together with Assumption \ref{ass: link core} and Lemma \ref{lem: properties bounded tensor op}.3 gives
  $$\left\lVert A_{T^*}A_{\mathcal{L}_X^{1/2}}C\right\rVert_{L^2(Z,\mathbb{Y} )}\geq m \left\lVert [\phi(A_{\mathcal{L}_X}^2)]^{1/2}A_{\mathcal{L}_X^{1/2}}C\right\rVert_{S_2}$$
  for all $C\in \N(Q)^\perp$ and therefore
  $$\left\lVert A_{T^*}C\right\rVert_{L^2(Z,\mathbb{Y} )}\geq m \left\lVert [\phi(A_{\mathcal{L}_X}^2)]^{1/2}C\right\rVert_{S_2}$$
  for all $C\in A_{\mathcal{L}_X^{1/2}}(\N(Q)^\perp)$.
  An application of Lemma \ref{lem:interpolation-high} with $W=A_{\mathcal{L}_X^{1/2}}(\N(Q)^\perp)$, $B=A_{\mathcal{L}_X^{-1}}$ and $T=A_{T^*}$ gives
  $$\norm{C}_\tau \leq \norm{C}_1 \left[ \psi^{-1}_{1/2} \left( \frac{\norm{A_{T^*}C}^2_{S_2(L_2(Z),\mathbb{Y} )}}{m^2\norm{C}^2_1}\right)\right]^{\frac{1-\gamma}{2}}$$
  for all $C\in A_{\mathcal{L}_X^{1/2}}(\N(Q)^\perp)$.
  Furthermore, since $[K]_X=A_{\mathcal{I}_X^*}K$, $A_{\widetilde{T}^*}=A_{T^*} A_{\mathcal{I}_X^*}$ and $A_{\widetilde{T}}A_{\widetilde{T}^*}=A_{\mathcal{C}_M}$ we obtain the upperbounds
  $$\lVert [K]_X\rVert_1=\lVert K\rVert_{S_2} \leq \lambda^{-1/2}\lVert K\mathcal{C}_{M,\lambda}^{1/2}\rVert_{S_2}  $$and
  $$\lVert A_{T^*}[K]_X\rVert_{S_2(L_2(Z),\mathbb{Y} )} =\lVert A_{\widetilde{T}^*}K\rVert_{S_2(L_2(Z),\mathbb{Y} )}=\lVert K\mathcal{C}_M^{1/2}\rVert_{S_2}\leq \lVert K\mathcal{C}_{M,\lambda}^{1/2}\rVert_{S_2} \,\lVert \mathcal{C}_{M,\lambda}^{-1/2}\mathcal{C}_{M}^{1/2}\rVert \leq \lVert K\mathcal{C}_{M,\lambda}^{1/2}\rVert_{S_2}.$$
  It now suffices to plug these upperbounds into the inequality to obtain the final result.
\end{proof}

\begin{lemma}[{\cite[Proposition~7]{rudi2015less}}]\label{lem:rudi}
Let $A,B$ be bounded, self-adjoint, positive semi-definite operators on a Hilbert space and $\lambda>0$.
If $\beta:=\|(B+\lambda\mathrm{Id})^{-1/2}(B-A)(B+\lambda\mathrm{Id})^{-1/2}\|<1$, then
$\|(A+\lambda\mathrm{Id})^{-1/2}(B+\lambda\mathrm{Id})^{1/2}\|\le(1-\beta)^{-1/2}$.
\end{lemma}

\begin{lemma}[{\cite[Lemma~6]{meunier2024nonparametric}}]\label{lem:FS17}
Let $\tau\ge1$, $\lambda\in(0,1]$ and $n\ge8\tau c_\lambda\kappa_X^2\lambda^{-1}$. Then, with
$\mathbb P_Z^{\otimes n}$-probability at least $1-2e^{-\tau}$,
$\|\mathcal C_{M,\lambda}^{-1/2}(\mathcal C_M-\widehat{\mathcal C}_M)\mathcal C_{M,\lambda}^{-1/2}\|\le\frac23$.
\end{lemma}

\begin{lemma}\label{lem:25}
Let Assumption \ref{ass:stage1}.\ref{ass:embZ} hold, and let $\tau\ge1$, $\lambda\in(0,1]$ and $n\ge8\tau c_\lambda\kappa_X^2\lambda^{-1}$.
Then with $\mathbb P^{\otimes n}_Z$-probability at least $1-2e^{-\tau}$ 
\begin{equation}\label{eq: sqrt 3 bound}
  \big\|\widehat{\mathcal C}_{M,\lambda}^{-1/2}\mathcal C_{M,\lambda}^{1/2}\big\|\le\sqrt3.
\end{equation}
Furthermore, let $u\in  E_X^m \times E_Z^m$ be such
that the constraints in \eqref{eq: set of equations} hold for $\widehat M=\widehat M(u)$, with a constant $J\ge J_0$, where
$$J_0:=\max\big\{2\kappa_X+1,\ A_Z(2\kappa_X+A_Z)\big\}.$$
Then there exists a measurable set $\Omega_0(u)\subseteq E_Z^n$
with $\mathbb P^{\otimes n}_Z\big(\Omega_0(u)^c\big)\le4e^{-\tau}$ on which the following bounds hold
simultaneously:
\begin{align}
  &\big\|\widehat{\mathcal C}_{M,\lambda}^{-1/2}\mathcal C_{M,\lambda}^{1/2}\big\|\le\sqrt3,\\
  &\big\|\widehat{\mathcal C}_M-\widehat{\mathcal C}_{\widehat M}\big\|
   \le J_0\Big(\sqrt{\tfrac{\tau}{n}}\,\|M-\widehat{M}\|_{\alpha_Z}+\|M-\widehat{M}\|_{L^2(Z,\mathcal H_{k_X})}\Big),
   \label{eq:L25-i}\\
  &\big\|\widehat{\mathcal C}_{\widehat M,\lambda}^{-1/2}\widehat{\mathcal C}_{M,\lambda}^{1/2}\big\|
   \le\sqrt{6/5}, \label{eq:L25-iii}\\
  &\big\|\mathcal C_{M,\lambda}^{1/2}\widehat{\mathcal C}_{\widehat M,\lambda}^{-1/2}\big\|\le\sqrt6,
   \label{eq:L25-iv}
\end{align}
\end{lemma}

\begin{proof}
Set $\Delta(z):=M(z)-\widehat{M}(z)$.

\emph{Step 1: bound~\eqref{eq:L25-i}.} Since $\widehat M=M-\Delta$, for every $z\in E_Z$,
$$M(z)\otimes M(z)-\widehat M(z)\otimes\widehat M(z)=M(z)\otimes\Delta(z)+\Delta(z)\otimes M(z)-\Delta(z)\otimes\Delta(z).$$
Using $\|u\otimes v\|=\|u\|\|v\|$ and $\|M(z)\|_{\mathcal H_{k_X}}\le\kappa_X$, and setting
$\xi_i:=\|\Delta(\widetilde Z_i)\|_{\mathcal H_{k_X}}$, we get
\begin{equation}\label{eq:L25-split}
  \big\|\widehat{\mathcal C}_M-\widehat{\mathcal C}_{\widehat M}\big\|
  \le\frac{2\kappa_X}{n}\sum_{i=1}^n\xi_i+\frac1n\sum_{i=1}^n\xi_i^2 .
\end{equation}
By Assumption \ref{ass:stage1}.\ref{ass:embZ} and \cite[Lemma~4]{li2022optimal},
$\|\Delta\|_{L_\infty(Z,\mathcal H_{k_X})}\le A_Z\|\Delta\|_{\alpha_Z}=:b$. Since
$\widetilde Z_1,\dots,\widetilde Z_n$ are i.i.d.\ with law $\mathbb P_Z$ and $u$ is fixed, the
$\xi_i$ are i.i.d.\ with $0\le\xi_i\le b$ almost surely. Moreover,
$\mathbb E\xi_1^2=\|\Delta\|^2_{L^2(Z,\mathcal H_{k_X})}$ and, by Jensen's inequality,
$\mathbb E\xi_1\le\|\Delta\|_{L^2(Z,\mathcal H_{k_X})}$. The one-sided Hoeffding inequality,
applied to $(\xi_i)$ and to $(\xi_i^2)$, gives an event $E_0$ with
$\mathbb P^{\otimes n}_Z(E_0^c)\le2e^{-\tau}$ on which
$$\frac1n\sum_{i=1}^n\xi_i\le\|\Delta\|_{L^2(Z,\mathcal H_{k_X})}+b\sqrt{\tfrac{\tau}{2n}}, \qquad \frac1n\sum_{i=1}^n\xi_i^2\le\|\Delta\|^2_{L^2(Z,\mathcal H_{k_X})}+b^2\sqrt{\tfrac{\tau}{2n}} .$$
By \eqref{eq: set of equations}, $\|\Delta\|_{L^2(Z,\mathcal H_{k_X})}\le1$ and $\|\Delta\|_{\alpha_Z}\le1$, so that squares
are bounded by first powers. Inserting the two bounds in~\eqref{eq:L25-split} yields, on
$E_0$,
$$\big\|\widehat{\mathcal C}_M-\widehat{\mathcal C}_{\widehat M}\big\|\le(2\kappa_X+1)\|\Delta\|_{L^2(Z,\mathcal H_{k_X})}+A_Z(2\kappa_X+A_Z)\sqrt{\tfrac{\tau}{n}}\,\|\Delta\|_{\alpha_Z},$$
which implies~\eqref{eq:L25-i}.

\emph{Step 2: a consequence of \eqref{eq: set of equations}.} Since $\|\widehat{\mathcal C}_{M,\lambda}^{-1/2}\|^2\le\lambda^{-1}$
and $\|\mathcal C_{M,\lambda}^{-1/2}\|^2\le\lambda^{-1}$, bound~\eqref{eq:L25-i} and the second and
third constraints in \eqref{eq: set of equations} with $J\ge J_0$ give
\begin{equation}\label{eq:L25-sixth}
  \lambda^{-1}\big\|\widehat{\mathcal C}_M-\widehat{\mathcal C}_{\widehat M}\big\|
  \le\frac{J_0\sqrt\tau\,\|\Delta\|_{\alpha_Z}}{\lambda\sqrt n}
  +\frac{J_0\|\Delta\|_{L^2(Z,\mathcal H_{k_X})}}{\lambda}
  \le\frac1{12}+\frac1{12}=\frac16 .
\end{equation}

\emph{Step 3: bound~\eqref{eq:L25-iii}.} Apply Lemma~\ref{lem:rudi} with
$A=\widehat{\mathcal C}_{\widehat M}$ and $B=\widehat{\mathcal C}_M$. By~\eqref{eq:L25-sixth},
$\beta\le\lambda^{-1}\|\widehat{\mathcal C}_M-\widehat{\mathcal C}_{\widehat M}\|\le\frac16$, hence
$\|\widehat{\mathcal C}_{\widehat M,\lambda}^{-1/2}\widehat{\mathcal C}_{M,\lambda}^{1/2}\|
\le(5/6)^{-1/2}=\sqrt{6/5}$.

\emph{Step 4: bound~\eqref{eq:L25-iv}.} Apply Lemma~\ref{lem:rudi} with
$A=\widehat{\mathcal C}_{\widehat M}$ and $B=\mathcal C_M$. By the triangle inequality,
Lemma \ref{lem:FS17} and~\eqref{eq:L25-sixth},
$$\beta\le\big\|\mathcal C_{M,\lambda}^{-1/2}(\mathcal C_M-\widehat{\mathcal C}_M)\mathcal C_{M,\lambda}^{-1/2}\big\|
  +\lambda^{-1}\big\|\widehat{\mathcal C}_M-\widehat{\mathcal C}_{\widehat M}\big\|
  \le\frac23+\frac16=\frac56,$$
hence $\|\widehat{\mathcal C}_{\widehat M,\lambda}^{-1/2}\mathcal C_{M,\lambda}^{1/2}\|\le\sqrt6$.
\end{proof}

\begin{lemma}\label{lem: C lambda control}
  Assume \ref{ass:stage2}.\ref{ass: source} with $\beta_X\leq1$ and the upper link condition in \ref{ass: link core}.
  Then, $$\left\lVert A_{\mathcal{L}_X^{-1/2}}C_\lambda  \right\rVert_{S_2(L_2(X),\mathbb{Y} )}\leq E (1+M)B_X \left[\psi_{1/2}^{-1}(\lambda)\right]^{\frac{\beta_X-1}{2}}.$$
\end{lemma}
\begin{proof}
  Using Proposition \ref{prop:problem is hilbert scaled}, we have that
  $$\left\lVert A_{\mathcal{L}_X^{-1/2}}C_\lambda  \right\rVert_{S_2(L_2(X),\mathbb{Y} )}=\left\lVert g_\lambda(Q^*Q)Q^*A_{T^*}C^\dagger \right\rVert_{S_2(L_2(X),\mathbb{Y} )}.$$
  Let $E_A$ be the spectral measure of $A_{\mathcal{L}_X^{-1}}$ and set $C_N:=E_A([0,N])C^\dagger$.
  Then,
  \begin{equation*}
    \begin{split}
      &\left\lVert g_\lambda(Q^*Q)Q^*A_{T^*}C^\dagger \right\rVert_{S_2(L_2(X),\mathbb{Y} )} \\ &\leq \left\lVert g_\lambda(Q^*Q)Q^*A_{T^*}C_N \right\rVert_{S_2(L_2(X),\mathbb{Y} )}+\left\lVert g_\lambda(Q^*Q)Q^*A_{T^*}(C^\dagger-C_N) \right\rVert_{S_2(L_2(X),\mathbb{Y} )}\\
    \end{split}
  \end{equation*}
  where
  \begin{equation*}
    \begin{split}
      \left\lVert g_\lambda(Q^*Q)Q^*A_{T^*}(C^\dagger-C_N) \right\rVert_{S_2(L_2(X),\mathbb{Y} )}^2&\leq E^2\lambda^{-1}\left\lVert A_{T^*}(C^\dagger-C_N) \right\rVert_{L^2(Z,\mathbb{Y} )}^2\\
      &\leq E^2M^2\lambda^{-1}\left\lVert [\phi(A_{\mathcal{L}_X^2})]^{1/2}(C^\dagger-C_N) \right\rVert_{S^2(L_2(X),\mathbb{Y} )}^2\\
      &=E^2M^2\lambda^{-1}\int_{(N,\infty)}\phi(t^{-2})\,d\lVert E_A(t)C^\dagger\rVert^2 \\
      &\leq E^2M^2\lambda^{-1}\phi(N^{-2})N^{-\beta_X}\int_{(N,\infty)}t^{\beta_X}\,d\lVert E_A(t)C^\dagger\rVert^2 \\
      &\leq E^2 M^2 \lVert C^\dagger\rVert_{\beta_X}^2\lambda^{-1}\phi(N^{-2})N^{-\beta_X}
    \end{split}
  \end{equation*}
  and
  \begin{equation*}
    \begin{split}
      \left\lVert g_\lambda(Q^*Q)Q^*A_{T^*}C_N \right\rVert_{S_2(L_2(X),\mathbb{Y} )}^2 &\leq E^2\left\lVert A_{\mathcal{L}_X^{-1/2}}C_N \right\rVert_{S_2(L_2(X),\mathbb{Y} )}^2\\
      &= E^2 \int_{[0,N]}t\,d\lVert E_A(t)C^\dagger\rVert^2\\
      &\leq E^2 \lVert C^\dagger\rVert_{\beta_X}^2 N^{1-\beta_X}.
    \end{split}
  \end{equation*}
  Hence, 
  $$\left\lVert A_{\mathcal{L}_X^{-1/2}}C_\lambda  \right\rVert_{S_2(L_2(X),\mathbb{Y} )}\leq EB_XM\lambda^{-1/2}[\phi(N^{-2})N^{-\beta_X}]^{1/2}+EB_X N^{\frac{1-\beta_X}{2}}.$$
  We can now balance both terms by choosing $N^{-1}=\psi^{-1}_{1/2}(\lambda)$ which gives the desired bound.
\end{proof}
\begin{lemma}\label{lem: E[eta]}
  Assume \ref{ass:stage2}.\ref{ass: source} and \ref{ass: link core} together with \ref{ass: additional convexity} if $\beta_X>1$. Then,
  $$\left\lVert \C_{MY}r_\lambda(\C_M)(\C_M+\lambda)^{-1/2}\right\rVert_{S_2} \leq c\, \lambda^{1/2} \left[\psi_{1/2}^{-1}\left(\lambda\right)\right]^{(\beta_X-1)/2} $$
  for some constant $c>0$ only depending on $B_X$, $\omega_\rho$, $\beta_X$ and $c_1$.
\end{lemma}
\begin{proof}
  We only treat the case $\beta_X>1$ as the case $\beta_X=1$ is easier and the case $\beta_X<1$ is treated by a similar truncation argument as in Lemma \ref{lem: C lambda control}.
  Set $V:=A_{\mathcal{L}_X^{-1/2}}C^\dagger$, $Z:=A_{\mathcal{L}_X^{1/2}}r_\lambda(Q^*Q)V$ and $\Delta:=\lVert r_\lambda(Q^*Q)^{1/2}V\rVert_{S_2(L^2(X),\mathbb{Y} )}^2 $.
  Set $h_\lambda:=r_\lambda(t)(t+\lambda)^{-1/2}$ and notice that
  \begin{equation*}
    \begin{split}
      \C_{MY}r_\lambda(\C_M)(\C_M+\lambda)^{-1/2}&= h_\lambda(\mathbf{I_X}^*A_{T}A_{T^*}\mathbf{I_X})\mathbf{I_X}^*A_T A_{T^*}C^\dagger\\
      &= \mathbf{I_X}^*A_{T}h_\lambda(A_{T^*}\mathbf{I_X}\mathbf{I_X}^*A_T) A_{T^*}C^\dagger\in\overline{\mathcal{R}(\mathbf{I_X^*})}= \mathcal{N}(\mathbf{I_X})^\perp.
    \end{split}
  \end{equation*}
  Hence, by Lemma \ref{lem: Li extended} (with $\gamma=1$), we have that
  \begin{equation*}
    \begin{split}
      (I):=\left\lVert \C_{MY}r_\lambda(\C_M)(\C_M+\lambda)^{-1/2}\right\rVert_{S_2}&=\left\lVert A_{\mathcal{L}_X^{-1/2}}\mathbf{I_X}\C_{MY}r_\lambda(\C_M)(\C_M+\lambda)^{-1/2}\right\rVert_{S_2(L^2(X),\mathbb{Y})}\\
      &=\left\lVert (Q^*Q+\lambda)^{-1/2}r_\lambda(Q^*Q)Q^*QV\right\rVert_{S_2(L^2(X),\mathbb{Y} )}
    \end{split}
  \end{equation*}
  where we have by Lemma \ref{lem: unitary equivalence L_2} with $A=Q^*Q$ that
  \begin{equation*}
    \begin{split}
      A_{\mathcal{L}_X^{-1/2}}\mathbf{I_X}\C_{MY}r_\lambda(\C_M)(\C_M+\lambda)^{-1/2}&=Wr_\lambda(W^*Q^*QW)(W^*Q^*QW+\lambda)^{-1/2}W^*Q^*QV\\
      &=(Q^*Q+\lambda)^{-1/2}r_\lambda(Q^*Q)Q^*QV.
    \end{split}
  \end{equation*}

  Furthermore, we obtain the following bounds,
  $$(I)^2\leq \omega_\rho^2 \lambda \Delta,\qquad \left\lVert TZ\right\rVert_{L^2(Z,\mathbb{Y} )}^2 \leq \omega_\rho^2 \lambda \Delta \qquad \text{and}\qquad \lVert Z\rVert_{1}^2\leq \Delta, $$
  as well as,
  $$\Delta =  \langle A_{\mathcal{L}_X^{-1/2}}Z,V\rangle_{S_2(L^2(X),\mathbb{Y} )}=\langle A_{\mathcal{L}_X^{(\beta_X-2)/2}}Z,A_{\mathcal{L}_X^{-\beta_X/2}}C^\dagger\rangle_{S_2(L^2(X),\mathbb{Y} )}\leq B_X \lVert Z\rVert_{2-\beta_X}. $$
  Letting $E$ be the spectral measure of $A_{\mathcal{L}_X^{-1}}$, we obtain in a similar fashion as in the proof of Proposition 5.1 in \cite{nair2005regularization}, that
  \begin{equation*}
    \begin{split}
      f\left(\frac{\left\lVert Z\right\rVert_{2-\beta_X}^2 }{\left\lVert Z\right\rVert_1^2 }\right)=f\left(\frac{\int t^{2-\beta_X}\,d\left\lVert E(t)Z\right\rVert^2}{\int t\,d\left\lVert E(t)Z\right\rVert^2}\right)
      \leq \frac{\int f(t^{1-\beta_X})t\,d\left\lVert E(t)Z\right\rVert^2}{\int t\,d\left\lVert E(t)Z\right\rVert^2}
      &=\frac{\int \phi(t^{-2})\,d\left\lVert E(t)Z\right\rVert^2}{\int t\,d\left\lVert E(t)Z\right\rVert^2}\\
      &=\frac{\lVert [\phi(A_{\mathcal{L}_X^2})]^{1/2}Z\rVert^2 }{\left\lVert Z\right\rVert_1^2}\\
      &\leq \frac{\left\lVert A_{T^*}Z\right\rVert^2 }{m^2\left\lVert Z\right\rVert_1^2}\\
      &\leq \frac{\omega_\rho^2 \lambda}{m^2}\frac{\Delta}{\left\lVert Z\right\rVert_1^2 }\\
      &\leq \frac{B_X \omega_\rho^2 \lambda}{m^2}\frac{\left\lVert Z\right\rVert_{2-\beta_X} }{\left\lVert Z\right\rVert_1^2 }
    \end{split}
  \end{equation*}
  Since, $f$ is convex with $f(0^+)=0$, we have that $f(t)/t$ is nondecreasing.
  Hence, since
  $$\frac{\left\lVert Z\right\rVert_1^2 }{B_X}\leq \frac{\Delta}{B_X}\leq \left\lVert Z\right\rVert_{2-\beta_X} $$
  and therefore $$\frac{\left\lVert Z\right\rVert_{2-\beta_X} }{B_X}\leq \frac{\left\lVert Z\right\rVert^2_{2-\beta_X}}{\left\lVert Z\right\rVert^2_1 },$$
  we obtain
  $$f\left(\frac{\left\lVert Z\right\rVert_{2-\beta_X} }{B_X}\right)\leq f\left(\frac{\left\lVert Z\right\rVert^2_{2-\beta_X}}{\left\lVert Z\right\rVert^2_1 }\right)\frac{\left\lVert Z\right\rVert^2_1 }{\left\lVert Z\right\rVert_{2-\beta_X}}\leq \frac{\omega_\rho^2}{m^2}\lambda.$$
  Finally, using $f^{-1}=[\psi_{1/2}^{-1}]^{\beta_X-1}$ which is increasing, we have
  \begin{equation*}
    \begin{split}
      (I)^2 \leq \omega_\rho^2 B_X \lambda \left\lVert Z\right\rVert_{2-\beta_X}&\leq \omega_\rho^2 B_X^2 \lambda \left[\psi_{1/2}^{-1}\left(\frac{\omega_\rho^2}{m^2}\lambda\right)\right]^{\beta_X-1}\\
      &\leq \omega_\rho^2 B_X^2 \left(\frac{\omega_\rho^2}{m^2}\vee 1\right)^{\beta_X-1} \lambda \left[\psi_{1/2}^{-1}\left(\lambda\right)\right]^{\beta_X-1},
    \end{split}
  \end{equation*}
  where we used the fact that $\psi_{1/2}^{-1}(ct)\leq (c\vee 1)\,\psi_{1/2}^{-1}(t)$ (see Lemma \ref{lem:nair 3.6}).
\end{proof}

\begin{lemma}\label{lem: E[eta] 2}
  Assume \ref{ass:stage2}.\ref{ass: source} and \ref{ass: link core} together with \ref{ass: additional convexity} if $\beta_X>1$. Then, on the event of Equation \eqref{eq: sqrt 3 bound}, we have that
  $$\left\lVert K_{\lambda}r_\lambda(\widehat{\C}_M)(\widehat{\C}_M+\lambda)^{1/2}\right\rVert_{S_2} \leq c\, \lambda^{1/2} \left[\psi_{1/2}^{-1}\left(\lambda\right)\right]^{(\beta_X-1)/2} $$

  for some constant $c>0$ only depending on $B_X$, $\omega_\rho$, $\beta_X$ and $m$.
\end{lemma}
\begin{proof}
  Let $\beta_X<1$. Then, by Lemmas \ref{lem: N IX membership}, \ref{lem: Li extended} and \ref{lem: C lambda control} we have that
  \begin{equation*}
    \begin{split}
      \left\lVert K_{\lambda}r_\lambda(\widehat{\C}_M)(\widehat{\C}_M+\lambda)^{1/2}\right\rVert &\leq \left\lVert r_\lambda(\widehat{\C}_M)^{1/2}\right\rVert \left\lVert r_\lambda(\widehat{\C}_M)^{1/2}(\widehat{\C}_M+\lambda)^{1/2}\right\rVert_{S_2}\left\lVert K_{\lambda}\right\rVert_{S_2}\\
      &\leq \omega_\rho \lambda^{1/2}  \left[\psi_{1/2}^{-1}\left(\lambda\right)\right]^{(\beta_X-1)/2}.
    \end{split}
  \end{equation*}
  Now let $\beta_X>1$ and set $\widehat{S}:=WA_{\hat{\C}_M}W^*$ and $V:=A_{\mathcal{L}_X^{-1/2}}C^\dagger$. Then, since $\overline{\mathcal{R}(\widehat{\C}_M)}\subseteq \overline{\mathcal{R}(\C_M)}$ almost surely (see Lemma \ref{lem: inclusion empirical cov}), we have that 
  $$\overline{\mathcal{R}(A_{\widehat{\C}_M})}\subseteq \overline{\mathcal{R} (A_{\C_M})}=\mathcal{N}(A_{T^*}\mathbf{I_X})^\perp \subseteq \mathcal{N}(\mathcal{\mathbf{I_X}})^\perp,$$
  and therefore by Lemma \ref{lem: unitary equivalence Hkx}
  $$K_\lambda r_\lambda(\widehat{\C}_M)(\widehat{\C}_M+\lambda)^{1/2}=W^*r_\lambda(\widehat{S})\widehat{S}_\lambda^{1/2}g_\lambda(Q^*Q)Q^*QV\in \mathcal{N}(\mathbf{I_X})^\perp,$$
  where we used the fact that $K_\lambda = W^*g_\lambda(Q^*Q)Q^*QV$.

  Hence, by Lemma \ref{lem: Li extended}, we obtain
  \begin{equation*}
    \begin{split}
      \left\lVert K_{\lambda}r_\lambda(\widehat{\C}_M)(\widehat{\C}_M+\lambda)^{1/2}\right\rVert_{1} 
      &=\left\lVert r_\lambda(\widehat{S})\widehat{S}_\lambda^{1/2}g_\lambda(Q^*Q)Q^*QV\right\rVert_1\\
      &\leq \left\lVert r_\lambda(\widehat{S})\widehat{S}_\lambda^{1/2}V\right\rVert_1 + \left\lVert r_\lambda(\widehat{S})\widehat{S}_\lambda^{1/2}r_\lambda(Q^*Q)V\right\rVert_1=:(I)+(II),\\
    \end{split}
  \end{equation*}
  where $V:=A_{\mathcal{L}_X^{-1/2}}C^\dagger$.
  For the second term, we have that
  \begin{equation*}
    \begin{split}
      (II)&\leq \lVert r_\lambda(\widehat{S})^{1/2}\rVert \lVert r_\lambda(Q^*Q)^{1/2}\rVert \lVert (r_\lambda(\widehat{S})\widehat{S}_\lambda)^{1/2}\rVert \lVert r_\lambda(Q^*Q)^{1/2}V\rVert_1 \\
      &\leq 2\omega_\rho B_X \left(\frac{\omega_\rho^2}{m^2}\vee 1\right)^{( \beta_X-1 )/2}\lambda^{1/2} \left[\psi_{1/2}^{-1}\left(\lambda\right)\right]^{(\beta_X-1)/2}
    \end{split}
  \end{equation*}
  where we used the bound obtained in the proof of Lemma \ref{lem: E[eta]} for $\lVert r_\lambda(Q^*Q)^{1/2}V\rVert $.
  For the first term, we have
  $$(I)\leq 2\omega_\rho \lambda^{1/2}\lVert r_\lambda(\widehat{S})^{1/2}V\rVert $$
  and we can now use similar arguments as in the proof of Lemma \ref{lem: E[eta]} to bound $\widehat{\Delta}:=\lVert r_\lambda(\widehat{S})^{1/2}V\rVert^2$.
  Set $\widehat{Z}:=A_{\mathcal{L}_X^{1/2}}r_\lambda(\widehat{S})V$ and one can check that $\widehat{Z}\in A_{\mathcal{L}_X^{1/2}}(\mathcal{N}(Q)^\perp)$ since $\widehat{S}\subset W(\mathcal{N}(A_{\widetilde{T}^*})^\perp)\subseteq \mathcal{N}(Q)^\perp$ and $V\in \mathcal{N}(Q)^\perp$. Then,
  $$\lVert \widehat{Z} \rVert^2_1= \lVert r_\lambda(\widehat{S})V \rVert^2\leq \widehat{\Delta}, $$
  $$\widehat{\Delta}=\big\langle V,r_\lambda(\widehat{S})V\big\rangle =\big\langle V,A_{\mathcal{L}_X^{-1/2}}\widehat{Z}\big\rangle =\big\langle A_{\mathcal{L}_X^{(\beta_X-1)/2}}A_{\mathcal{L}_X^{-\beta_X/2}}C^\dagger,A_{\mathcal{L}_X^{-1/2}}\widehat{Z}\big\rangle \leq B_X \lVert \widehat{Z}\rVert_{2-\beta_X} .$$ and
  $$\lVert A_{T^*}\widehat{Z} \rVert = \lVert Qr_\lambda(\widehat{S})V \rVert \leq \lVert Q(Q^*Q)_\lambda^{-1/2}\rVert \lVert  (Q^*Q)_\lambda^{1/2} \widehat{S}_\lambda^{-1/2}\rVert \lVert \widehat{S}_\lambda^{1/2} r_\lambda(\widehat{S})V \rVert \leq 2\sqrt{3}\omega_\rho \lambda^{1/2}\Delta^{1/2},$$
  with $\mathbb{P}_Z^{\otimes n}$-probability at least $1-4e^{-\tau}$, where we used the Lemmas \ref{lem:25} and \ref{lem: unitary equivalence Hkx} to obtain that
  $$\lVert  (Q^*Q)_\lambda^{1/2} \widehat{S}_\lambda^{-1/2}\rVert = \lVert  W A_{\C_{M,\lambda}^{1/2}} A_{\hat{\C}_{M,\lambda}^{-1/2}}W^*\rVert \leq \lVert  \C_{M,\lambda}^{1/2} \widehat{\C}_{M,\lambda}^{-1/2}\rVert \leq \sqrt{3}.$$
  One can now easily copy the arguments of the proof of Lemma \ref{lem: E[eta]} to obtain the desired bound.

  The proof for the case $\beta_X=1$ is easier and therefore omitted.
\end{proof}

\begin{lemma}\label{lem:nair 3.6}
One has that
$$\min\{c,1\}\,\psi_{1/2}^{-1}(\lambda)\;\le\;\psi_{1/2}^{-1}(c\lambda)\;\le\;\max\{c,1\}\,\psi_{1/2}^{-1}(\lambda),$$
for all $c,\lambda>0$.
\end{lemma}

\begin{proof}
The upperbound is shown by equation (3.6) in \cite{nair2005regularization}.
Let $c\geq 1$, then since $\psi_{1/2}^{-1}$ is nondecreasing we have that $\psi_{1/2}^{-1}(\lambda)\le \psi_{1/2}^{-1}(c\lambda)$.
Now let $c<1$ and $t>0$. Due to monotonicity of $\phi$, we have that $c \psi_{1/2}(t)=ct \phi(t^2)\geq c t \phi((ct)^2)=\psi_{1/2}(ct)$, and therefore
$\psi_{1/2}^{-1}\big(c\psi_{1/2}(t))\geq ct$. Choosing $\lambda =\psi_{1/2}(t)$ gives the final result.
\end{proof}

\begin{lemma}\label{lem: effective dimension}
Let $\mathcal{N}_M(\lambda):=\operatorname{tr}\big(\mathcal{C}_M(\mathcal{C}_M+\lambda)^{-1}\big)$. Under Assumptions \ref{ass:stage2}.\ref{ass: px} and \ref{ass: link core}
we have for sufficiently small $\lambda$ that
$$\mathcal{N}_M(\lambda)\leq c[\psi_{1/2}^{-1}(\lambda)]^{-p_X},$$
for some constant $c>0$ only depending on $c_2,l_2$ and $p_X$.
\end{lemma}

\begin{proof}
  Assume wlog $l_2=1$.
Setting $\widetilde{T}=T \mathcal{I}_X$, we have that $\C_M=\widetilde{T}^*\widetilde{T}$.
And therefore, by the upper inequality of Assumption \ref{ass: link core}, we have that
\begin{equation*}
  \begin{split}
    \langle \mathcal{C}_Mh,h\rangle =\lVert T\mathcal{I}_Xh\rVert_{L_2(Z)}^2
    &\leq c_2^2\big\lVert [\phi(\mathcal{L}_X^2)]^{1/2}\mathcal{I}_Xh\big\rVert_{L_2(X)}^2\\
    &=c_2^2\big\lVert \mathcal{I}_X[\phi(\mathcal{C}_X^2)]^{1/2}h\big\rVert_{L_2(X)}^2\\
    &=c_2^2\langle \psi_{1/2}(\mathcal{C}_X)h,h\rangle ,
  \end{split}
\end{equation*}
that is, $\mathcal{C}_M\leq c_2^2\psi_{1/2}(\mathcal{C}_X)$, and therefore also $(\mathcal{C}_M+\lambda)^{-1}\geq (c_2^2\psi_{1/2}(\mathcal{C}_X)+\lambda)^{-1}$.

Since $t\mapsto -t^{-1}$ is operator monotone on $(0,\infty)$, we obtain
$$\mathcal{C}_M(\mathcal{C}_M+\lambda)^{-1}=\Id-\lambda(\mathcal{C}_M+\lambda)^{-1}
\leq \Id-\lambda\big(c_2^2\psi_{1/2}(\mathcal{C}_X)+\lambda\big)^{-1}
=c_2^2\psi_{1/2}(\mathcal{C}_X)\big(c_2^2\psi_{1/2}(\mathcal{C}_X)+\lambda\big)^{-1}.$$
Hence,
$$\mathcal{N}_M(\lambda)\leq tr\left(c_2^2\psi_{1/2}(\mathcal{C}_X)\big(c_2^2\psi_{1/2}(\mathcal{C}_X)+\lambda\big)^{-1}\right)=\sum_i \frac{c_2^2\psi_{1/2}(\mu_{X,i})}{c_2^2\psi_{1/2}(\mu_{X,i})+\lambda}.$$

Set $t_\lambda:=\psi_{1/2}^{-1}(\lambda/c_2^2)\in(0,\lVert\mathcal{L}_X\rVert]$, which holds for sufficiently small $\lambda$. If $\mu_{X,i}\geq t_\lambda$, the $i$-th term
above is bounded by $1$. If $\mu_{X,i}<t_\lambda$, then, since the map
$t\mapsto\psi_{1/2}(t)/t=\phi(t^2)$ is nondecreasing, we have that
$c_2^2\psi_{1/2}(\mu_{X,i})\leq c_2^2\psi_{1/2}(t_\lambda)\,\mu_{X,i}/t_\lambda=\lambda\,\mu_{X,i}/t_\lambda$ and the $i$-th term is bounded by
$\mu_{X,i}/t_\lambda$. Therefore
$$\mathcal{N}_M(\lambda)\leq \sum_{i\geq 1}\min\left(1,\frac{\mu_{X,i}}{t_\lambda}\right).$$
On the one hand, $\sum_i\mu_{X,i}/t_\lambda$ is bounded by $\operatorname{tr}(\mathcal{L}_X)/t_\lambda$ which shows the result for $p_X=1$.
On the other hand, we have $\min(1,\mu_{X,i}/t_\lambda)\leq\min\{1,(c_2/t_\lambda)\,i^{-1/p_X}\}$ which is nonincreasing in $i$.
Hence, if $p_X<1$, then
\begin{equation*}
  \begin{split}
    \sum_{i\geq 1}\min\left(1,\frac{\mu_{X,i}}{t_\lambda}\right)&\leq\int_0^\infty\min\left\{1,\frac{c_2}{t_\lambda}u^{-1/p_X}\right\}du\\
    &=\int_0^{(c_2/t_\lambda)^{p_X}}\min\left\{1,\frac{c_2}{t_\lambda}u^{-1/p_X}\right\}du + \int_{(c_2/t_\lambda)^{p_X}}^\infty\min\left\{1,\frac{c_2}{t_\lambda}u^{-1/p_X}\right\}du\\
    &\leq \left(\frac{c_2}{t_\lambda}\right)^{p_X}+\frac{c_2}{t_\lambda}\int_{(c_2/t_\lambda)^{p_X}}^\infty u^{-1/p_X}du\\
    &=\left(1+\frac{p_X}{1-p_X}\right)\left(\frac{c_2}{t_\lambda}\right)^{p_X}.
  \end{split}
\end{equation*}
\end{proof}

\begin{theorem}[Bernstein's inequality, see Theorem 26 in \cite{fischer2020sobolev}]\label{thm:bernstein}
Let $H$ be a separable Hilbert space and let $\xi_1,\dots,\xi_n$ be independent $H$-valued random
variables with $\mathbb E[\xi_i]=0$ for all $i$. Assume that there are constants
$\tilde\sigma,\tilde L>0$ such that
$$\sum_{i=1}^n\mathbb E\|\xi_i\|_H^q\le\frac{q!}{2}\,\tilde\sigma^2\tilde L^{q-2}\qquad\text{for all integers } q\ge2 .$$
Then, for every $\tau>0$, with probability at least $1-2e^{-\tau}$,
$$\Big\|\frac1n\sum_{i=1}^n\xi_i\Big\|_H\le\frac{\sqrt{2\tilde\sigma^2\tau}}{n}+\frac{2\tilde L\tau}{n}.$$
\end{theorem}

\begin{theorem}[Theorem 4, \cite{meunier2024optimal}]\label{th: meunier optimal}
  Let $g_\varepsilon$ be a filter with qualification greater than $\beta_Z/2$ and
  let Assumption \ref{ass:stage1} hold.
  Then, for $\varepsilon=O(m^{-1/(\beta_Z+p_z)})$, there exist constants $J,J'>0$ independent of $m$ and $\tau$, such that with $\mathbb{P}$-probability not less than $1-6e^{-\tau}$, we have
  $$\lVert M-\widehat{M}\rVert_{L^2(Z,\mathcal{H}_{K_X})}^2\leq J\tau^2 m^{-\beta_Z/(\beta_Z+p_Z)}\qquad \text{and}\qquad \lVert M-\widehat{M}\rVert_{\alpha_Z}\leq J' \tau^2 m^{(\alpha_Z-\beta_Z)/(\beta_Z+p_Z)}.$$
\end{theorem}

\end{document}